\documentclass[reqno]{amsart}
\usepackage{amsmath}
\usepackage{amssymb}
\usepackage{amsthm}
\usepackage{cite}
\usepackage{euscript}
\usepackage[arrow,curve,matrix,arc,2cell]{xy}
\usepackage[hyphens]{url}
\usepackage{hyperref}
\UseAllTwocells
\DeclareFontFamily{U}{rsfs}{} \DeclareFontShape{U}{rsfs}{n}{it}{<->
rsfs10}{} \DeclareSymbolFont{mscr}{U}{rsfs}{n}{it}
\DeclareSymbolFontAlphabet{\scr}{mscr}
\def\mathscr{\scr}
\begin{document}
%%%%%%%%%%%%%%%%%%%%%%%%%%%%%%%%%%%%%%%%%%%%%%%%%%%%%%%%%%%%%%%%%%%%%%%%
%%%%%%%%%%%%%%%%%%%%%%%%%%     Macros      %%%%%%%%%%%%%%%%%%%%%%%%%%%%%
%%%%%%%%%%%%%%%%%%%%%%%%%%%%%%%%%%%%%%%%%%%%%%%%%%%%%%%%%%%%%%%%%%%%%%%%
\def\e#1\e{\begin{equation}#1\end{equation}}
\def\ea#1\ea{\begin{align}#1\end{align}}
\def\eq#1{{\rm(\ref{#1})}}
\theoremstyle{plain}% default
\newtheorem{thm}{Theorem}[section]
\newtheorem{lem}[thm]{Lemma}
\newtheorem{prop}[thm]{Proposition}
\newtheorem{cor}[thm]{Corollary}
\newtheorem{princ}[thm]{Principle}
\newtheorem{claim}[thm]{Claim}
\theoremstyle{definition}
\newtheorem{dfn}[thm]{Definition}
\newtheorem{altdfn}[thm]{Alternative Definition}
\newtheorem{ex}[thm]{Example}
\newtheorem{rem}[thm]{Remark}
\numberwithin{figure}{section}
\numberwithin{equation}{section}
\def\dim{\mathop{\rm dim}\nolimits}
\def\codim{\mathop{\rm codim}\nolimits}
\def\vdim{\mathop{\rm vdim}\nolimits}
\def\depth{\mathop{\rm depth}\nolimits}
\def\Im{\mathop{\rm Im}\nolimits}
\def\Ker{\mathop{\rm Ker}}
\def\Coker{\mathop{\rm Coker}}
\def\ev{\mathop{\rm ev}\nolimits}
\def\Ho{\mathop{\rm Ho}}
\def\GL{\mathop{\rm GL}}
\def\Stab{\mathop{\rm Stab}\nolimits}
\def\supp{\mathop{\rm supp}}
\def\rank{\mathop{\rm rank}\nolimits}
\def\Sym{\mathop{\rm Sym}\nolimits}
\def\Aff{\mathop{\bf Aff}\nolimits}
\def\Sch{\mathop{\bf Sch}\nolimits}
\def\Sets{\mathop{\bf Sets}\nolimits}
\def\Groupoids{\mathop{\bf Groupoids}\nolimits}
\def\Hom{\mathop{\rm Hom}\nolimits}
\def\bHom{\mathop{\bf Hom}\nolimits}
\def\cHom{\mathop{\mathcal{H}om}\nolimits}
\def\cIso{\mathop{\mathcal{I}so}\nolimits}
\def\bcHom{\mathop{\bs{\mathcal{H}om}}\nolimits}
\def\bcEqu{\mathop{\bs{\mathcal{E}qu}}\nolimits}
\def\id{{\mathop{\rm id}\nolimits}}
\def\Obj{{\rm Obj}}
\def\CSch{{\mathop{\bf C^{\bs\iy}Sch}}}
\def\muKur{{\mathop{\bs\mu\bf Kur}\nolimits}}
\def\muKurb{{\mathop{\bs\mu\bf Kur^b}\nolimits}}
\def\muKurc{{\mathop{\bs\mu\bf Kur^c}\nolimits}}
\def\muKurcin{{\mathop{\bs\mu\bf Kur^c_{in}}\nolimits}}
\def\muKurcsi{{\mathop{\bs\mu\bf Kur^c_{si}}\nolimits}}
\def\muKurcis{{\mathop{\bs\mu\bf Kur^c_{is}}\nolimits}}
\def\muKurcst{{\mathop{\bs\mu\bf Kur^c_{st}}\nolimits}}
\def\muKurgc{{\mathop{\bs\mu\bf Kur^{gc}}}}
\def\muKurgcin{{\mathop{\bs\mu\bf Kur^{gc}_{in}}}}
\def\muKurgcsi{{\mathop{\bs\mu\bf Kur^{gc}_{si}}}}
\def\cmuKurc{{\mathop{\bs\mu\bf\check Kur^c}\nolimits}}
\def\cmuKurcin{{\mathop{\bs\mu\bf\check Kur^c_{in}}\nolimits}}
\def\cmuKurcsi{{\mathop{\bs\mu\bf\check Kur^c_{si}}\nolimits}}
\def\cmuKurcis{{\mathop{\bs\mu\bf\check Kur^c_{is}}\nolimits}}
\def\cmuKurcst{{\mathop{\bs\mu\bf\check Kur^c_{st}}\nolimits}}
\def\cmuKurgc{{\mathop{\bs\mu\bf\check Kur^{gc}}}}
\def\cmuKurgcin{{\mathop{\bs\mu\bf\check Kur^{gc}_{in}}}}
\def\cmuKurgcsi{{\mathop{\bs\mu\bf\check Kur^{gc}_{si}}}}
\def\mKur{{\mathop{\bf mKur}\nolimits}}
\def\mKurb{{\mathop{\bf mKur^b}\nolimits}}
\def\mKurc{{\mathop{\bf mKur^c}\nolimits}}
\def\mKurcin{{\mathop{\bf mKur^c_{in}}\nolimits}}
\def\mKurcsi{{\mathop{\bf mKur^c_{si}}\nolimits}}
\def\mKurcis{{\mathop{\bf mKur^c_{is}}\nolimits}}
\def\mKurcst{{\mathop{\bf mKur^c_{st}}\nolimits}}
\def\mKurgc{{\mathop{\bf mKur^{gc}}}}
\def\mKurgcin{{\mathop{\bf mKur^{gc}_{in}}}}
\def\mKurgcsi{{\mathop{\bf mKur^{gc}_{si}}}}
\def\cmKurc{{\mathop{\bf m\check{K}ur^c}\nolimits}}
\def\cmKurcin{{\mathop{\bf m\check{K}ur^c_{in}}\nolimits}}
\def\cmKurcsi{{\mathop{\bf m\check{K}ur^c_{si}}\nolimits}}
\def\cmKurcis{{\mathop{\bf m\check{K}ur^c_{is}}\nolimits}}
\def\cmKurcst{{\mathop{\bf m\check{K}ur^c_{st}}\nolimits}}
\def\cmKurgc{{\mathop{\bf m\check{K}ur^{gc}}}}
\def\cmKurgcin{{\mathop{\bf m\check{K}ur^{gc}_{in}}}}
\def\cmKurgcsi{{\mathop{\bf m\check{K}ur^{gc}_{si}}}}
\def\Kur{{\mathop{\bf Kur}\nolimits}}
\def\Kurb{{\mathop{\bf Kur^b}\nolimits}}
\def\Kurc{{\mathop{\bf Kur^c}\nolimits}}
\def\Kurcin{{\mathop{\bf Kur^c_{in}}\nolimits}}
\def\Kurcsi{{\mathop{\bf Kur^c_{si}}\nolimits}}
\def\Kurcis{{\mathop{\bf Kur^c_{is}}\nolimits}}
\def\Kurcst{{\mathop{\bf Kur^c_{st}}\nolimits}}
\def\Kurgc{{\mathop{\bf Kur^{gc}}}}
\def\Kurgcin{{\mathop{\bf Kur^{gc}_{in}}}}
\def\Kurgcsi{{\mathop{\bf Kur^{gc}_{si}}}}
\def\cKurcin{{\mathop{\bf \check{K}ur^c_{in}}\nolimits}}
\def\cKurcsi{{\mathop{\bf \check{K}ur^c_{si}}\nolimits}}
\def\cKurcis{{\mathop{\bf \check{K}ur^c_{is}}\nolimits}}
\def\cKurcst{{\mathop{\bf \check{K}ur^c_{st}}\nolimits}}
\def\cKurc{{\mathop{\bf \check{K}ur^c}\nolimits}}
\def\cKurgc{{\mathop{\bf \check{K}ur^{gc}}\nolimits}}
\def\cKurgcin{{\mathop{\bf \check{K}ur^{gc}_{in}}\nolimits}}
\def\cKurgcsi{{\mathop{\bf \check{K}ur^{gc}_{si}}\nolimits}}
\def\KurtrG{{\mathop{\bf Kur_{trG}}\nolimits}}
\def\KurtrGc{{\mathop{\bf Kur_{trG}^c}\nolimits}}
\def\KurtrGgc{{\mathop{\bf Kur_{trG}^{gc}}\nolimits}}
\def\KN{{\mathop{\bf KN}\nolimits}}
\def\mKN{{\mathop{\bf mKN}\nolimits}}
\def\GKN{{\mathop{\bf GKN}\nolimits}}
\def\GmKN{{\mathop{\bf GmKN}\nolimits}}
\def\Man{{\mathop{\bf Man}}}
\def\Manb{{\mathop{\bf Man^b}}}
\def\Manc{{\mathop{\bf Man^c}}}
\def\Mangc{{\mathop{\bf Man^{gc}}}}
\def\Mancst{{\mathop{\bf Man^c_{st}}}}
\def\Mancsi{{\mathop{\bf Man^c_{si}}}}
\def\Mancis{{\mathop{\bf Man^c_{is}}}}
\def\Mancin{{\mathop{\bf Man^c_{in}}}}
\def\Mangcin{{\mathop{\bf Man^{gc}_{in}}}}
\def\Mangcsi{{\mathop{\bf Man^{gc}_{si}}}}
\def\cManc{{\mathop{\bf\check{M}an^c}}}
\def\cMancst{{\mathop{\bf\check{M}an^c_{st}}}}
\def\cMancis{{\mathop{\bf\check{M}an^c_{is}}}}
\def\cMancsi{{\mathop{\bf\check{M}an^c_{si}}}}
\def\cMancin{{\mathop{\bf\check{M}an^c_{in}}}}
\def\cMangc{{\mathop{\bf\check{M}an^{gc}}}}
\def\cMangcin{{\mathop{\bf\check{M}an^{gc}_{in}}}}
\def\cMangcsi{{\mathop{\bf\check{M}an^{gc}_{si}}}}
\def\dMan{{\mathop{\bf dMan}}}
\def\dOrb{{\mathop{\bf dOrb}}}
\def\dSpa{{\mathop{\bf dSpa}}}
\def\SMod{{\mathop{\bf SMod}}}
\def\Orb{{\mathop{\bf Orb}}}
\def\Orbb{{\mathop{\bf Orb^b}}}
\def\Orbc{{\mathop{\bf Orb^c}}}
\def\OrbKur{{\mathop{\bf Orb_{\rm Kur}}}}
\def\Mon{{\mathop{\bf Mon}}}
\def\ul{\underline}
\def\bs{\boldsymbol}
\def\ge{\geqslant}
\def\le{\leqslant\nobreak}
\def\pr{\prec}
\def\O{{\mathcal O}}
\def\R{{\mathbin{\mathbb R}}}
\def\Z{{\mathbin{\mathbb Z}}}
\def\Q{{\mathbin{\mathbb Q}}}
\def\N{{\mathbin{\mathbb N}}}
\def\C{{\mathbin{\mathbb C}}}
\def\CP{{\mathbin{\mathbb{CP}}}}
\def\RP{{\mathbin{\mathbb{RP}}}}
\def\fC{{\mathbin{\mathfrak C}\kern.05em}}
\def\fD{{\mathbin{\mathfrak D}}}
\def\fE{{\mathbin{\mathfrak E}}}
\def\fF{{\mathbin{\mathfrak F}}}
\def\cA{{\mathbin{\mathcal A}}}
\def\cB{{\mathbin{\mathcal B}}}
\def\cC{{\mathbin{\mathcal C}}}
\def\cD{{\mathbin{\mathcal D}}}
\def\cE{{\mathbin{\mathcal E}}}
\def\cF{{\mathbin{\mathcal F}}}
\def\cG{{\mathbin{\mathcal G}}}
\def\cH{{\mathbin{\mathcal H}}}
\def\cI{{\mathbin{\mathcal I}}}
\def\cJ{{\mathbin{\mathcal J}}}
\def\cK{{\mathbin{\mathcal K}}}
\def\cL{{\mathbin{\mathcal L}}}
\def\cM{{\mathbin{\mathcal M}}}
\def\cN{{\mathbin{\mathcal N}}}
\def\cP{{\mathbin{\mathcal P}}}
\def\cQ{{\mathbin{\mathcal Q}}}
\def\cR{{\mathbin{\mathcal R}}}
\def\cS{{\mathbin{\mathcal S}}}
\def\cT{{\mathbin{\mathcal T}\kern -0.1em}}
\def\cW{{\mathbin{\mathcal W}}}
\def\fV{{\mathbin{\mathfrak V}}}
\def\fW{{\mathbin{\mathfrak W}}}
\def\fX{{\mathbin{\mathfrak X}}}
\def\fY{{\mathbin{\mathfrak Y}}}
\def\fZ{{\mathbin{\mathfrak Z}}}
\def\fe{{\mathfrak e}}
\def\ff{{\mathfrak f}}
\def\fg{{\mathfrak g}}
\def\fh{{\mathfrak h}}
\def\oM{{\mathbin{\smash{\,\,\overline{\!\!\mathcal M\!}\,}}}}
\def\bcE{{\mathbin{\bs{\mathcal E}}}}
\def\ur{{\underline{r\kern -0.15em}\kern 0.15em}{}}
\def\uU{{{\underline{U\kern -0.25em}\kern 0.2em}{}}{}}
\def\uX{{{\underline{X\!}\,}{}}{}}
\def\cU{{\mathcal U}}
\def\cV{{\mathcal V}}
\def\cW{{\mathcal W}}
\def\bA{{\bs A}}
\def\bB{{\bs B}}
\def\bC{{\bs C}}
\def\bD{{\bs D}}
\def\bE{{\bs E}}
\def\bF{{\bs F}}
\def\bG{{\bs G}}
\def\bH{{\bs H}}
\def\bM{{\bs M}}
\def\bN{{\bs N}}
\def\bO{{\bs O}}
\def\bP{{\bs P}}
\def\bQ{{\bs Q}}
\def\bS{{\bs S}}
\def\bT{{\bs T}}
\def\bU{{\bs U}}
\def\bV{{\bs V}}
\def\bW{{\bs W}\kern -0.1em}
\def\bX{{\bs X}}
\def\bY{{\bs Y}\kern -0.1em}
\def\bZ{{\bs Z}}
\def\al{\alpha}
\def\be{\beta}
\def\ga{\gamma}
\def\de{\delta}
\def\io{\iota}
\def\ep{\epsilon}
\def\la{\lambda}
\def\ka{\kappa}
\def\th{\theta}
\def\ze{\zeta}
\def\up{\upsilon}
\def\vp{\varphi}
\def\si{\sigma}
\def\om{\omega}
\def\De{\Delta}
\def\La{\Lambda}
\def\Om{\Omega}
\def\Io{{\rm I}}
\def\Ka{{\rm K}}
\def\Mu{{\rm M}}
\def\Nu{{\rm N}}
\def\Tau{{\rm T}}
\def\Up{\Upsilon}
\def\Al{{\rm A}}
\def\Be{{\rm B}}
\def\Ga{\Gamma}
\def\Si{\Sigma}
\def\Th{\Theta}
\def\utau{{\underline{\tau\kern -0.2em}\kern 0.2em}{}}
\def\uchi{{\underline{\chi\kern -0.1em}\kern 0.1em}{}}
\def\ueta{{\underline{\eta\kern -0.1em}\kern 0.1em}{}}
\def\pd{\partial}
\def\ts{\textstyle}
\def\st{\scriptstyle}
\def\sst{\scriptscriptstyle}
\def\w{\wedge}
\def\sm{\setminus}
\def\lt{\ltimes}
\def\bu{\bullet}
\def\sh{\sharp}
\def\op{\oplus}
\def\od{\odot}
\def\op{\oplus}
\def\ot{\otimes}
\def\ov{\overline}
\def\bigop{\bigoplus}
\def\bigot{\bigotimes}
\def\iy{\infty}
\def\es{\emptyset}
\def\ra{\rightarrow}
\def\rra{\rightrightarrows}
\def\Ra{\Rightarrow}
\def\RRa{\Rrightarrow}
\def\Longra{\Longrightarrow}
\def\ab{\allowbreak}
\def\longra{\longrightarrow}
\def\hookra{\hookrightarrow}
\def\dashra{\dashrightarrow}
\def\ha{{\ts\frac{1}{2}}}
\def\t{\times}
\def\ci{\circ}
\def\ti{\tilde}
\def\d{{\rm d}}
\def\md#1{\vert #1 \vert}
\def\bmd#1{\big\vert #1 \big\vert}
\def\an#1{\langle #1 \rangle}
\def\ban#1{\bigl\langle #1 \bigr\rangle}
%%%%%%%%%%%%%%%%%%%%%%%%%%%%%%%%%%%%%%%%%%%%%%%%%%%%%%%%%%%%%%%%%%%%%%%%
%%%%%%%%%%%%%%%%%%%%%%%%    Text of paper    %%%%%%%%%%%%%%%%%%%%%%%%%%%
%%%%%%%%%%%%%%%%%%%%%%%%%%%%%%%%%%%%%%%%%%%%%%%%%%%%%%%%%%%%%%%%%%%%%%%%
\title{Kuranishi spaces as a 2-category}
\author{Dominic Joyce}
\date{}

\begin{abstract} 
This is a survey of the author's paper \cite{Joyc5} and in-progress book \cite{Joyc6}. `Kuranishi spaces' were introduced in the work of Fukaya, Oh, Ohta and Ono \cite{Fuka,FOOO1,FOOO2,FOOO3,FOOO4,FOOO5,FOOO6,FOOO7,FOOO8,FOOO9,FuOn} in symplectic geometry, as the geometric structure on moduli spaces of $J$-holomorphic curves. We propose a new definition of Kuranishi space, which has the nice property that they form a 2-category~$\Kur$. 

Any Fukaya--Oh--Ohta--Ono (FOOO) Kuranishi space $\bX$ can be made into a Kuranishi space $\bX'$ uniquely up to equivalence in $\Kur$. The same holds for McDuff and Wehrheim's `Kuranishi atlases' \cite{McDu,McWe1,McWe2,McWe3}, and Hofer, Wysocki and Zehnder's `polyfold Fredholm structures'~\cite{Hofe,HWZ1,HWZ2,HWZ3,HWZ4,HWZ5,HWZ6,HWZ7}.

Our Kuranishi spaces are based on the author's theory of Derived Differential Geometry \cite{Joyc2,Joyc3,Joyc4}, the study of classes of derived manifolds and orbifolds that we call `d-manifolds' and `d-orbifolds'. There is an equivalence of 2-categories $\Kur\simeq\dOrb$, where $\dOrb$ is the 2-category of d-orbifolds. So Kuranishi spaces are really a form of derived orbifold. They have their own differential geometry, with notions of orientation, immersions, submersions, transverse fibre products, and so on.
\end{abstract}
\maketitle

\setcounter{tocdepth}{2}
\tableofcontents

\section{Introduction}
\label{kt1}

{\it Kuranishi spaces\/} were introduced in the work of Fukaya, Oh, Ohta and Ono \cite{Fuka,FOOO1,FOOO2,FOOO3,FOOO4,FOOO5,FOOO6,FOOO7,FOOO8,FOOO9,FuOn}, as the geometric structure on moduli spaces of $J$-holomorphic curves. We will refer to their most recent definition \cite[\S 4]{FOOO4} as {\it FOOO Kuranishi spaces}. They are used to define virtual cycles and virtual chains for such moduli spaces, for applications in symplectic geometry such as Gromov--Witten invariants \cite{FuOn} and Lagrangian Floer cohomology \cite{FOOO1}.

In related work, McDuff and Wehrheim \cite{McDu,McWe1,McWe2,McWe3} define {\it Kuranishi atlases}, which are quite similar to FOOO `good coordinate systems' in\cite{Fuka,FOOO1,FOOO2,FOOO3,FOOO4,FOOO5,FOOO6,FOOO7,FOOO8,FOOO9,FuOn}, and Dingyu Yang defines {\it Kuranishi structures\/} \cite{Yang1,Yang2,Yang3}, which are very similar to \cite{Fuka,FOOO1,FOOO2,FOOO3,FOOO4,FOOO5,FOOO6,FOOO7,FOOO8,FOOO9,FuOn}. An alternative theory, philosophically rather different, but which does essentially the same job, is the {\it polyfolds\/} of Hofer, Wysocki and Zehnder~\cite{Hofe,HWZ1,HWZ2,HWZ3,HWZ4,HWZ5,HWZ6,HWZ7}.

Although FOOO Kuranishi spaces are well defined, and are adequate for the applications in Fukaya et al.\ \cite{Fuka,FOOO1,FOOO2,FOOO3,FOOO4,FOOO5,FOOO6,FOOO7,FOOO8,FOOO9,FuOn}, they seem not very satisfactory as geometric spaces. For example, there is currently no good notion of morphism of FOOO Kuranishi spaces.

This is a partial survey of the author's paper \cite{Joyc5} and in-progress book \cite{Joyc6}. In it we present a new definition of Kuranishi spaces, which form a (weak) 2-category $\Kur$. That is, we have objects, the Kuranishi spaces $\bX,\bY,\ldots,$ and 1-morphisms $\bs f\colon \bX\ra\bY$, and also 2-morphisms $\bs\eta\colon \bs f\Ra\bs g$ between 1-morphisms~$\bs f,\bs g\colon \bX\ra\bY$.

If $\bX$ is a FOOO Kuranishi space, we can construct a Kuranishi space $\bX{}'$ in our sense with the same topological space, uniquely up to canonical equivalence in $\Kur$. The same holds for McDuff and Wehrheim's `Kuranishi atlases' \cite{McDu,McWe1,McWe2,McWe3}, and Hofer, Wysocki and Zehnder's `polyfold Fredholm structures'~\cite{Hofe,HWZ1,HWZ2,HWZ3,HWZ4,HWZ5,HWZ6,HWZ7}.

Our theory of Kuranishi spaces is based on the author's theory of Derived Differential Geometry \cite{Joyc2,Joyc3,Joyc4}, the study of derived manifolds and orbifolds, where `derived' is in the sense of the Derived Algebraic Geometry of Lurie \cite{Luri} and To\"en--Vezzosi \cite{Toen}. The author's derived manifolds and orbifolds are called `d-manifolds' and `d-orbifolds', and form strict 2-categories $\dMan,\dOrb$. For alternative definitions of derived manifolds, see Spivak \cite{Spiv} and Borisov--Noel \cite{BoNo}. The relations between \cite{Joyc2,Joyc3,Joyc4} and \cite{BoNo,Spiv} are explained by Borisov~\cite{Bori}.

D-manifolds and d-orbifolds \cite{Joyc2,Joyc3,Joyc4}, and also the derived manifolds of \cite{BoNo,Spiv}, are defined using $C^\iy$-algebraic geometry, as in \cite{Joyc1}: they are classes of derived schemes and stacks over $C^\iy$-rings. Their definition looks very different to that of Kuranishi spaces. Nonetheless, there is an equivalence of 2-categories $\Kur\simeq\dOrb$, so Kuranishi spaces and d-orbifolds are essentially the same thing. The author used results on d-orbifolds in \cite{Joyc4} to design the definition of Kuranishi spaces below, so as to arrange that~$\Kur\simeq\dOrb$.

One moral is that {\it Kuranishi spaces are really derived smooth orbifolds}, which does not seem to be well known amongst symplectic geometers. So we should use ideas from derived geometry to understand Kuranishi spaces. One such lesson is that `derived' objects generally form higher categories. In fact, even classical orbifolds are best defined to be a 2-category. Hence it should not be surprising that Kuranishi spaces are also a 2-category.

One can turn Kuranishi spaces into an ordinary category by passing to the {\it homotopy category\/} $\Ho(\Kur)$, with objects Kuranishi spaces $\bX,\bY,\ldots,$ and whose morphisms $[\bs f]\colon \bX\ra\bY$ are 2-isomorphism classes of 1-morphisms $\bs f\colon \bX\ab\ra\bY$ in $\Kur$. But the definition of $\Ho(\Kur)$ is not much easier than that of $\Kur$, and there are disadvantages to doing this:
\begin{itemize}
\setlength{\itemsep}{0pt}
\setlength{\parsep}{0pt}
\item Morphisms $[\bs f]\colon \bX\ra\bY$ in $\Ho(\Kur)$ are not local on $\bX$, that is, they do not form a sheaf on $\bX$.
\item Important constructions such as fibre products $\bX\t_{\bs g,\bZ,\bs h}\bY$ in $\Kur$  are characterized by universal properties involving 2-morphisms in $\Kur$. They do not satisfy any universal property in~$\Ho(\Kur)$.
\end{itemize}

All the Kuranishi-type structures that we consider are given by an `atlas of charts' on a topological space $X$, where the `charts', called {\it Kuranishi neighbourhoods}, are quintuples $(V,E,\Ga,s,\psi)$ for $V$ a manifold, $E\ra V$ a vector bundle, $\Ga$ a finite group acting smoothly on $V,E$, $s\colon V\ra E$ a $\Ga$-equivariant smooth section, and $\psi\colon s^{-1}(0)/\Ga\ra X$ a continuous map which is a homeomorphism with an open set~$\Im\psi\subseteq X$.

The important differences between FOOO Kuranishi spaces, and MW Kuranishi atlases, and our Kuranishi spaces, are in the `coordinate changes' between the charts in the atlas. Since our Kuranishi spaces are a 2-category, our Kuranishi neighbourhoods must be a 2-category too, as one can regard a Kuranishi neighbourhood as a Kuranishi space with only one chart. Thus, we define 1-{\it morphisms\/} $\Phi_{ij}\colon (V_i,E_i,\Ga_i,s_i,\psi_i)\ra (V_j,E_j,\Ga_j,s_j,\psi_j)$ between Kuranishi neighbourhoods on $X$, which include FOOO and MW coordinate changes as special cases, and also 2-{\it morphisms\/} $\La_{ij}\colon \Phi_{ij}\Ra\Phi_{ij}'$ between 1-morphisms. A {\it coordinate change\/} is a 1-morphism $\Phi_{ij}$ which is invertible up to 2-isomorphism.

Given Kuranishi neighbourhoods $(V_i,E_i,\Ga_i,s_i,\psi_i),(V_j,E_j,\Ga_j,s_j,\psi_j)$ on $X$, our 1- and 2-morphisms between them have the crucial property that they form a {\it stack}, or 2-{\it sheaf}, on $\Im\psi_i\cap\Im\psi_j$. This is a 2-categorical version of a sheaf, and means that 1- and 2-morphisms glue nicely on open covers. The stack property is essential in defining compositions $\bs g\ci\bs f\colon \bX\ra\bZ$ of 1-morphisms $\bs f\colon \bX\ra\bY$ and $\bs g\colon \bY\ra\bZ$ of Kuranishi spaces. The main reason why FOOO Kuranishi spaces and MW Kuranishi atlases do not have good notions of morphisms is that their coordinate changes lack such a sheaf/stack property. 

Our theory of Kuranishi spaces is part of a programme by the author \cite{Joyc6} to rewrite the foundations of areas of symplectic geometry involving moduli spaces of $J$-holomorphic curves. For such applications it will be an advantage that our Kuranishi spaces form a well-behaved 2-category. Here are three examples:
\begin{itemize}
\setlength{\itemsep}{0pt}
\setlength{\parsep}{0pt}
\item[(a)] Defining a Kuranishi structure on a moduli space $\cM$ involves making many arbitrary choices, and it is helpful to know how different choices are related. In our theory, we expect different choices of Kuranishi structure on $\cM$ to yield equivalent Kuranishi spaces $\bs\cM,\bs\cM'$ in the 2-category~$\Kur$.

Equivalence of Kuranishi spaces should be compared with McDuff and Wehrheim's `commensurate' Kuranishi atlases, as in \S\ref{kt23}, and Yang's `R-equivalent' Kuranishi structures, as in~\S\ref{kt25}.
\item[(b)] Fukaya, Oh, Ohta and Ono use finite group actions on Kuranishi spaces in \cite{Fuka}, \cite[\S 7]{FOOO8}, and $T^n$-actions on Kuranishi spaces in \cite{FOOO2,FOOO3}. In our theory it is easy to define and study actions of Lie groups on Kuranishi spaces, in a more flexible way than in~\cite{Fuka,FOOO2,FOOO3,FOOO8}.

\item[(c)] Fukaya \cite[\S 3, \S 5]{Fuka} (see also \cite[\S 4.2]{FOOO6}) works with a forgetful morphism $\mathfrak{forget}\colon \bs\cM_{l,1}(\be)\ra\bs\cM_{l,0}(\be)$ of $J$-holomorphic curve moduli spaces. This is a kind of morphism of Fukaya--Oh--Ohta--Ono Kuranishi spaces, but can only be defined when the Kuranishi structures on $\bs\cM_{l,1}(\be),\bs\cM_{l,0}(\be)$ are compatible in a strong sense, as they are in the examples of \cite{Fuka,FOOO6}. Our theory defines (1-)morphisms without such compatibility conditions.
\end{itemize}

In this survey we will for simplicity restrict attention almost entirely to Kuranishi spaces {\it without boundary}, which we just call Kuranishi spaces. Kuranishi spaces with boundary, and with corners, are studied in~\cite{Joyc5,Joyc6}.

Section \ref{kt2} reviews the definitions of Fukaya--Oh--Ohta--Ono's Kuranishi spaces \cite{Fuka,FOOO1,FOOO2,FOOO3,FOOO4,FOOO5,FOOO6,FOOO7,FOOO8,FOOO9,FuOn}, McDuff--Wehrheim's Kuranishi atlases \cite{McDu,McWe1,McWe2,McWe3}, and Dingyu Yang's Kuranishi structures \cite{Yang1,Yang2,Yang3}, and we also discuss Hofer--Wysocki--Zehnder's polyfolds \cite{Hofe,HWZ1,HWZ2,HWZ3,HWZ4,HWZ5,HWZ6,HWZ7}. Section \ref{kt3} explains Kuranishi neighbourhoods and their 1- and 2-morphisms and coordinate changes. Section \ref{kt4} defines the 2-category of Kuranishi spaces, and explains their relation to FOOO Kuranishi spaces, MW Kuranishi atlases, DY Kuranishi structures, polyfolds, and d-orbifolds. Section \ref{kt5} discusses the differential geometry of Kuranishi spaces. Finally, Appendix \ref{ktA} gives background on 2-categories, and stacks on topological spaces.

\medskip

\noindent{\it Acknowledgements.} I would like to thank Lino Amorim, Kenji Fukaya, Helmut Hofer, Dusa McDuff, and Dingyu Yang for helpful conversations. This research was supported by EPSRC grants EP/H035303/1 and EP/J016950/1.

\section{Previous definitions of Kuranishi space}
\label{kt2}

We explain {\it Kuranishi spaces\/} in the work of Fukaya, Oh, Ohta and Ono \cite{Fuka,FOOO1,FOOO2,FOOO3,FOOO4,FOOO5,FOOO6,FOOO7,FOOO8,FOOO9,FuOn}, and the similar idea of {\it Kuranishi atlases\/} in the work of McDuff and Wehrheim \cite{McDu,McWe1,McWe2,McWe3}. We also discuss  Dingyu Yang's {\it Kuranishi structures\/} \cite{Yang1,Yang2,Yang3}, a minor variation on Fukaya--Oh--Ohta--Ono Kuranishi spaces, and the rather different {\it polyfolds\/} of Hofer, Wysocki and Zehnder \cite{Hofe,HWZ1,HWZ2,HWZ3,HWZ4,HWZ5,HWZ6,HWZ7}. We have made some small changes in notation compared to our sources, for compatibility with \S\ref{kt3}--\S\ref{kt5}. We hope the authors concerned will not mind this.

\subsection{\texorpdfstring{Fukaya--Oh--Ohta--Ono's Kuranishi spaces}{Fukaya-Oh-Ohta-Ono\textquoteright s Kuranishi spaces}}
\label{kt21}

Kuranishi spaces are used in the work of Fukaya, Oh, Ohta and Ono \cite{Fuka,FOOO1,FOOO2,FOOO3,FOOO4,FOOO5,FOOO6,FOOO7,FOOO8,FOOO9,FuOn} as the geometric structure on moduli spaces of $J$-holomorphic curves. Initially introduced by Fukaya and Ono \cite[\S 5]{FuOn} in 1999, the definition has changed several times as their work has evolved \cite{Fuka,FOOO1,FOOO2,FOOO3,FOOO4,FOOO5,FOOO6,FOOO7,FOOO8,FOOO9,FuOn}. 

This section explains their most recent definition of Kuranishi space, taken from \cite[\S 4]{FOOO4}. As in the rest of our paper `Kuranishi neighbourhood', `coordinate change' and `Kuranishi space' have a different meaning, we will use the terms `FOOO Kuranishi neighbourhood', `FOOO coordinate change' and `FOOO Kuranishi space' below to refer to concepts from~\cite{FOOO4}.

For the next definitions, let $X$ be a compact, metrizable topological space. 

\begin{dfn} A {\it FOOO Kuranishi neighbourhood\/} on $X$ is a quintuple $(V,\ab E,\ab\Ga,\ab s,\ab\psi)$ such that:
\begin{itemize}
\setlength{\itemsep}{0pt}
\setlength{\parsep}{0pt}
\item[(a)] $V$ is a smooth manifold, which may or may not have boundary or corners.
\item[(b)] $E$ is a finite-dimensional real vector space.
\item[(c)] $\Ga$ is a finite group with a smooth, effective action on $V$, and a linear representation on $E$.
\item[(d)] $s\colon V\ra E$ is a $\Ga$-equivariant smooth map.
\item[(e)] $\psi$ is a homeomorphism from $s^{-1}(0)/\Ga$ to an open subset $\Im\psi$ in $X$, where $\Im\psi=\bigl\{\psi(x\Ga)\colon x\in s^{-1}(0)\bigr\}$ is the image of $\psi$, and is called the {\it footprint\/} of~$(V,E,\Ga,s,\psi)$.
\end{itemize}
We will write $\bar\psi\colon s^{-1}(0)\ra\Im\psi\subseteq X$ for the composition of $\psi$ with the projection~$s^{-1}(0)\ra s^{-1}(0)/\Ga$.

Now let $p\in X$. A {\it FOOO Kuranishi neighbourhood of\/ $p$ in\/} $X$ is a FOOO Kuranishi neighbourhood $(V_p,E_p,\Ga_p,s_p,\psi_p)$ with a distinguished point $o_p\in V_p$ such that $o_p$ is fixed by $\Ga_p$, and $s_p(o_p)=0$, and $\psi_p([o_p])=p$. Then $o_p$ is unique.

\label{kt2def1}
\end{dfn}

\begin{dfn} Let $(V_i,E_i,\Ga_i,s_i,\psi_i)$, $(V_j,E_j,\Ga_j,s_j,\psi_j)$ be FOOO Kuranishi nei\-gh\-bour\-hoods on $X$. Suppose $S\subseteq \Im\psi_i\cap\Im\psi_j\subseteq X$ is an open subset of the intersection of the footprints $\Im\psi_i,\Im\psi_j\subseteq X$. We say a quadruple $\Phi_{ij}=(V_{ij},h_{ij},\vp_{ij},\hat\vp_{ij})$ is a {\it FOOO coordinate change from\/ $(V_i,E_i,\Ga_i,s_i,\psi_i)$ to $(V_j,\ab E_j,\ab\Ga_j,\ab s_j,\ab\psi_j)$ over\/} $S$ if:
\begin{itemize}
\setlength{\itemsep}{0pt}
\setlength{\parsep}{0pt}
\item[(a)] $V_{ij}$ is a $\Ga_i$-invariant open neighbourhood of $\bar\psi_i^{-1}(S)$ in~$V_i$.
\item[(b)] $h_{ij}\colon \Ga_i\ra\Ga_j$ is an injective group homomorphism.
\item[(c)] $\vp_{ij}\colon V_{ij}\hookra V_j$ is an $h_{ij}$-equivariant smooth embedding, such that the induced map $(\vp_{ij})_*\colon  V_{ij}/\Ga_i\ra V_j/\Ga_j$ is injective. 
\item[(d)] $\hat\vp_{ij}\colon V_{ij}\t E_i\hookra V_j\t E_j$ is an $h_{ij}$-equivariant embedding of vector bundles over $\vp_{ij}\colon V_{ij}\hookra V_j$, viewing $V_{ij}\t E_i\ra V_{ij}$, $V_j\t E_j\ra V_j$ as trivial vector bundles.
\item[(e)] $\hat\vp_{ij}(s_i\vert_{V_{ij}})=\vp_{ij}^*(s_j)$, in sections of $\vp_{ij}^*(V_j\t E_j)\ra V_{ij}$.
\item[(f)] $\psi_i=\psi_j\ci(\vp_{ij})_*$ on $(s_i^{-1}(0)\cap V_{ij})/\Ga_i$.
\item[(g)] $h_{ij}$ restricts to an isomorphism $\Stab_{\Ga_i}(v)\ra \Stab_{\Ga_j}(\vp_{ij}(v))$ for all $v$ in $V_{ij}$, where $\Stab_{\Ga_i}(v)$ is the {\it stabilizer subgroup\/} $\bigl\{\ga\in\Ga_i\colon \ga(v)=v\bigr\}$.
\item[(h)] For each $v\in s_i^{-1}(0)\cap V_{ij}\subseteq V_{ij}\subseteq V_i$ we have a commutative diagram
\e
\begin{gathered}
\xymatrix@C=29pt@R=13pt{ 0 \ar[r] & T_vV_i \ar[rr]_{\d\vp_{ij}\vert_v} \ar[d]^{\d s_i\vert_v} && T_{\vp_{ij}(v)}V_j \ar[r] \ar[d]^{\d s_j\vert_{\vp_{ij}(v)}} & N_{ij}\vert_v \ar[r] \ar@{.>}[d]^{\d_{\rm fibre}s_j\vert_v} & 0 \\
0 \ar[r] & E_i\vert_v \ar[rr]^{\hat\vp_{ij}\vert_v} && E_j\vert_{\vp_{ij}(v)} \ar[r] & F_{ij}\vert_v \ar[r] & 0 
}
\end{gathered}
\label{kt2eq1}
\e
with exact rows, where $N_{ij}\ra V_{ij}$ is the normal bundle of $V_{ij}$ in $V_j$, and $F_{ij}=\vp_{ij}^*(E_j)/\hat\vp_{ij}(E_i\vert_{V_{ij}})$ the quotient bundle. We require that the induced morphism $\d_{\rm fibre}s_j\vert_v$ in \eq{kt2eq1} should be an isomorphism.
\end{itemize}

Note that $\d_{\rm fibre}s_j\vert_v$ an isomorphism in \eq{kt2eq1} is equivalent to the following complex being exact, which is how we write the analogous conditions in \S\ref{kt3}:
\e
\begin{gathered}
\xymatrix@C=17pt{ 0 \ar[r] & T_vV_i \ar[rrr]^(0.35){\d s_i\vert_v\op\d\vp_{ij}\vert_v} &&& E_i\vert_v \op T_{\vp_{ij}(v)}V_j 
\ar[rrr]^(0.56){\hat\vp_{ij}\vert_v\op -\d s_j\vert_{\vp_{ij}(v)}} &&& E_j\vert_{\vp_{ij}(v)} \ar[r] & 0. }
\end{gathered}
\label{kt2eq2}
\e

Now let $(V_p,E_p,\Ga_p,s_p,\psi_p)$, $(V_q,E_q,\Ga_q,s_q,\psi_q)$ be FOOO Kuranishi neighbourhoods of $p\in X$ and $q\in\Im\psi_p\subseteq X$, respectively. We say a quadruple $\Phi_{qp}=(V_{qp},h_{qp},\vp_{qp},\hat\vp_{qp})$ is a {\it FOOO coordinate change\/} if it is a FOOO coordinate change from $(V_q,E_q,\Ga_q,s_q,\psi_q)$ to $(V_p,E_p,\Ga_p,s_p,\psi_p)$ over $S_{qp}$, where $S_{qp}$ is any open neighbourhood of $q$ in $\Im\psi_q\cap\Im\psi_p$.
\label{kt2def2}
\end{dfn}

\begin{rem} Fukaya et al.\ \cite{FOOO4} only impose Definition \ref{kt2def2}(h) for Kuranishi spaces `with a tangent bundle'. As the author knows of no reason for considering Kuranishi spaces `without tangent bundles', and the notation appears to be merely historical, we will include `with a tangent bundle' in our definitions of FOOO coordinate changes and FOOO Kuranishi spaces.
\label{kt2rem1}
\end{rem}

\begin{dfn} A {\it FOOO Kuranishi structure\/} $\cK$ on $X$ of {\it virtual dimension\/} $n$ in $\Z$ in the sense of \cite[\S 4]{FOOO4}, including the `with a tangent bundle' condition, assigns a FOOO Kuranishi neighbourhood $(V_p,\ab E_p,\ab\Ga_p,s_p,\psi_p)$ for each $p\in X$ and a FOOO coordinate change $\Phi_{qp}=\ab(V_{qp},\ab h_{qp},\ab\vp_{qp},\ab\hat\vp_{qp})\colon (V_q,E_q,\Ga_q,s_q,\psi_q)\ra(V_p,E_p,\Ga_p,s_p,\psi_p)$ for each $q\in \Im\psi_p$ such that the following holds:
\begin{itemize}
\setlength{\itemsep}{0pt}
\setlength{\parsep}{0pt}
\item[(a)] $\dim V_p-\rank E_p=n$ for all $p\in X$.
\item[(b)] If $q\in\Im\psi_p$, $r\in\psi_q((V_{qp}\cap s_q^{-1}(0))/\Ga_q)$, then for each connected component $(\vp_{rq}^{-1}(V_{qp})\cap V_{rp})^\al$ of $\vp_{rq}^{-1}(V_{qp})\cap V_{rp}$ there exists $\ga_{rqp}^\al\in \Ga_p$ with
\e
\begin{gathered}
h_{qp}\ci h_{rq}=\ga_{rqp}^{\al}\cdot h_{rp} \cdot (\ga_{rqp}^{\al})^{-1}, \qquad
\vp_{qp} \ci \vp_{rq} = \ga_{rqp}^{\al}
\cdot \vp_{rp}, \\
\text{and}\qquad \vp_{rq}^*(\hat\vp_{qp}) \ci \hat\vp_{rq} = \ga_{rqp}^{\al}\cdot
\hat\vp_{rp},
\end{gathered}
\label{kt2eq3}
\e
where the second and third equations hold on $(\vp_{rq}^{-1}(V_{qp})\cap V_{rp})^\al$.
\end{itemize}

The pair $\bX=(X,\cK)$ is called a {\it FOOO Kuranishi space}, of {\it virtual dimension\/} $n\in\Z$, written $\vdim\bX=n$. If the $V_p$ for all $p\in X$ are manifolds without boundary, or with boundary, or with corners, then we call $\bX$ a {\it FOOO Kuranishi space without boundary}, or {\it with boundary}, or {\it with corners}, respectively.
\label{kt2def3}
\end{dfn}

FOOO Kuranishi spaces are differential-geometric spaces, generalizations of manifolds or orbifolds. Some of the differential geometry of manifolds extends to FOOO Kuranishi spaces. Here are analogues of the notions of orientations, smooth maps, and transverse fibre products of manifolds. 

\begin{dfn} Let $\bX$ be a FOOO Kuranishi space. Then for each $p\in X$, $q\in\Im\psi_p$ and $v\in s_q^{-1}(0)\cap V_{qp}$, we have an exact sequence \eq{kt2eq2}. Taking top exterior powers in \eq{kt2eq2} yields an isomorphism \begin{equation*}
\bigl(\det T_vV_q\bigr)\ot \det\bigl(E_p\vert_{\vp_{qp}(v)}\bigr)\cong
\bigl(\det E_q\vert_v\bigr)\ot\bigl(T_{\vp_{qp}(v)}V_p\bigr),
\end{equation*}
where $\det W$ means $\La^{\dim W}W$, or equivalently, a canonical isomorphism
\e
\bigl(\det T^*V_p\ot\det E_p\bigr)\vert_{\vp_{qp}(v)}\cong 
\bigl(\det T^*V_q\ot\det E_q\bigr)\vert_v.
\label{kt2eq4}
\e
Defining the isomorphism \eq{kt2eq4} requires a suitable sign convention. Sign conventions are discussed in Fukaya et al. \cite[\S 8.2]{FOOO1}. An {\it orientation\/} on $\bX$ is a choice of orientations on the line bundles 
\begin{equation*} 
\det T^*V_p\ot\det E_p\big\vert_{s_p^{-1}(0)} \longra s_p^{-1}(0)
\end{equation*}
for all $p\in X$, compatible with the isomorphisms~\eq{kt2eq4}.
\label{kt2def4}
\end{dfn}

\begin{dfn} Let $\bX$ be a FOOO Kuranishi space, and $Y$ a manifold. A {\it smooth map\/} $\bs f\colon \bX\ra Y$ is $\bs f=(f_p\colon p\in X)$ where $f_p\colon V_p\ra Y$ is a $\Ga_p$-invariant smooth map for all $p\in X$ (that is, $f_p$ factors via $V_p\ra V_p/\Ga_p\ra Y$), and $f_p\ci\vp_{qp}=f_q\vert_{V_{qp}}\colon V_{qp}\ra Y$ for all $q\in\Im\psi_p$. This induces a unique continuous map $f\colon X\ra Y$ with $f_p\vert_{s_p^{-1}(0)}=f\ci\bar\psi_p$ for all $p\in X$. We call $\bs f$ {\it weakly submersive\/} if each $f_p$ is a submersion.

Suppose $\bX,\bX'$ are FOOO Kuranishi spaces, $Y$ is a manifold, and $\bs f\colon \bX\ab\ra Y$, $\bs f'\colon \bX'\ra Y$ are weakly submersive. Then as in \cite[\S A1.2]{FOOO1} one can define a `fibre product' Kuranishi space $\bW=\bX\t_Y\bX'$, with topological space $W=\bigl\{(p,p')\in X\t X'\colon f(p)=f'(p')\bigr\}$, and FOOO Kuranishi neighbourhoods $(V_{p,p'},E_{p,p'},\Ga_{p,p'},s_{p,p'},\psi_{p,p'})$ for $(p,p')\in W$, where $V_{p,p'}=V_p\t_{f_p,Y,f_{\smash{p'}}'}V_{\smash{p'}}'$, $E_{p,p'}=\pi_{V_p}^*(E_p)\op \pi_{V_{\smash{p'}}'}^*(E'_{\smash{p'}})$, $\Ga_{p,p'}=\Ga_p\t\Ga_{\smash{p'}}'$, $s_{p,p'}=\pi_{V_p}^*(s_p)\op \pi_{V_{\smash{p'}}'}^*(s'_{\smash{p'}})$, and $\psi_{p,p'}=\psi_p\ci(\pi_{V_p})_*\t\psi_{\smash{p'}}'\ci(\pi_{V_{\smash{p'}}'})_*$. The weakly submersive condition ensures that $V_{p,p'}=V_p\t_YV_{\smash{p'}}'$ is well-defined.
\label{kt2def5}
\end{dfn}

\begin{rem}{\bf(i)} Fukaya et al.\ \cite{Fuka,FOOO1,FOOO2,FOOO3,FOOO4,FOOO5,FOOO6,FOOO7,FOOO8,FOOO9,FuOn} do not define morphisms between general FOOO Kuranishi spaces, so Kuranishi spaces in \cite{Fuka,FOOO1,FOOO2,FOOO3,FOOO4,FOOO5,FOOO6,FOOO7,FOOO8,FOOO9,FuOn} do not form a category. However, they do work with two special classes of morphisms. Smooth maps $\bs f\colon \bX\ra Y$ from Kuranishi spaces $\bX$ to manifolds $Y$ are important in their theory. Also Fukaya \cite[\S 3, \S 5]{Fuka} (see also \cite[\S 4.2]{FOOO6}) works with a forgetful morphism $\mathfrak{forget}\colon \bs\cM_{l,1}(\be)\ra\bs\cM_{l,0}(\be)$. This is an example of a kind of morphism $\bs f\colon \bX\ra\bY$ of FOOO Kuranishi spaces (though not given a formal definition in \cite{Fuka,FOOO6}), which can only be defined when the Kuranishi structures on $\bX,\bY$ are compatible in a strong sense, as they are by construction in~\cite{Fuka,FOOO6}. 
\smallskip

\noindent{\bf(ii)} The `fibre product' $\bX\t_Y\bX'$ in Definition \ref{kt2def5} is not a fibre product in the sense of category theory, characterized by a universal property, since Fukaya et al.\ in \cite{Fuka,FOOO1,FOOO2,FOOO3,FOOO4,FOOO5,FOOO6,FOOO7,FOOO8,FOOO9,FuOn} do not have a category of FOOO Kuranishi spaces in which to state such a universal property (though we do, see \S\ref{kt4} and \S\ref{kt52}). Their `fibre product' is extra data, an operation defined upon FOOO Kuranishi spaces.
\label{kt2rem2}
\end{rem}

\subsection{How FOOO Kuranishi spaces are used}
\label{kt22}

The next two theorems summarize some of the main results of \cite{FOOO1,FOOO3,FOOO4,FOOO7,FOOO5,FuOn}.

\begin{thm} The following classes of moduli spaces of\/ $J$-holomorphic curves in symplectic geometry may be given the structure of compact FOOO Kuranishi spaces of known dimension, after making many choices in the construction:
\begin{itemize}
\setlength{\itemsep}{0pt}
\setlength{\parsep}{0pt}
\item[{\bf(a)}] Moduli spaces $\bs\oM_{g,m}(J,\be)$ of stable maps $u\colon \Si\ra S$ from a prestable closed Riemann surface $\Si$ of genus $g$ with\/ $m$ marked points, to a compact symplectic manifold $(S,\om)$ with compatible almost complex structure $J,$ with homology class\/~$u_*([\Si])=\be\in H_2(S;\Z)$.

\item[{\bf(b)}] Moduli spaces $\bs\oM_k(L_1,\ldots,L_l,J,\be)$ of stable maps $u\colon \Si\ra S$ from a pre\-sta\-ble holomorphic disc $\Si$ with\/ $k$ boundary marked points, to a compact symplectic manifold $(S,\om)$ with compatible almost complex structure $J,$ such that $u(\pd\Si)$ lies in the union $L_1\cup\cdots\cup L_l$ of pairwise transversely intersecting, embedded, compact Lagrangians $L_1,\ldots,L_l$ in $S,$ with relative homology class\/~$u_*([\Si])=\be\in H_2(S,L_1\cup\cdots\cup L_l;\Z)$.
\end{itemize}
Here $\bs\oM_{g,m}(J,\be)$ is a FOOO Kuranishi space without boundary used to define Gromov--Witten invariants in\/ {\rm\cite{FuOn},} and\/ $\bs\oM_k(L_1,\ab\ldots,\ab L_l,J,\be)$ a FOOO Kuranishi space with corners used to define Lagrangian Floer cohomology in\/~{\rm\cite{FOOO1}}.

The `evaluation maps' at marked points $\bs\ev_i\colon \bs\oM_{g,m}(J,\be)\ra S$ and\/ $\bs\ev_j\ab\colon\ab \bs\oM_k(L_1,\ab\ldots,L_l,J,\be)\ra L_1\amalg\cdots\amalg L_l$ may be given the structure of weakly submersive smooth maps.

After choosing orientations and relative spin structures for $L_1,\ldots,L_l$ in {\bf(b)\rm,} there are canonical orientations on $\bs\oM_{g,m}(J,\be)$ and\/ $\bs\oM_k(L_1,\ldots,L_l,J,\be)$.
\label{kt2thm1}
\end{thm}

\begin{thm} Suppose $\bX$ is a compact, oriented FOOO Kuranishi space of virtual dimension $k,$ and\/ $Y$ is an oriented manifold of dimension $n,$ and\/ $\bs f\colon \bX\ra Y$ is a weakly submersive smooth map. Then, after making many choices, one can construct a \begin{bfseries}virtual (co)chain\end{bfseries} $[\bX]_{\rm virt}$ for $\bX,$ either:
\begin{itemize}
\setlength{\itemsep}{0pt}
\setlength{\parsep}{0pt}
\item[{\bf(i)}] in the smooth singular chains $C_k^{\rm ssi}(Y;\Q)$ of\/ $Y$ over $\Q;$ or 
\item[{\bf(ii)}] in the compactly-supported de Rham cochains $C^\iy_{\rm cs}(\La^{n-k}T^*Y)$.
\end{itemize}
If\/ $\pd\bX=\es$ then $[\bX]_{\rm virt}$ is a (co)cycle, and has a (co)homology class $[[\bX]_{\rm virt}],$ the \begin{bfseries}virtual cycle\end{bfseries}, in smooth singular homology $H_k^{\rm ssi}(Y;\Q),$ or in compactly-supported de Rham cohomology $H^{n-k}_{\rm dR,cs}(Y;\R)$. This $[[\bX]_{\rm virt}]$ is independent of choices in its construction, and depends on $(\bX,\bs f)$ only up to oriented bordism.

\label{kt2thm2}
\end{thm}

Thus, FOOO Kuranishi spaces are used as an intermediate stage in the construction of virtual chains or virtual cycles for moduli spaces of $J$-holomorphic curves. These virtual chains/cycles must also have other important properties which we have not stated --- basically, geometric relationships between moduli spaces should translate to algebraic relationships between their virtual chains/cycles. The virtual chains/cycles are used to define Gromov--Witten theory, Lagrangian Floer cohomology, Fukaya categories, and so on. 

The original proofs of Theorems \ref{kt2thm1} and \ref{kt2thm2} in \cite{FOOO1,FuOn} were widely criticized as not wholly correct/complete. Recently, Fukaya et al.\ \cite{FOOO4,FOOO5,FOOO7,FOOO9} have been working on providing more careful, detailed, and complete proofs. We should also acknowledge that the work of Fukaya, Oh, Ohta and Ono is very original, even visionary, and full of important ideas which other authors have used since. The work we discuss in \S\ref{kt23}--\S\ref{kt26} came 5 to 15 years after the first versions of \cite{FOOO1,FuOn}, partly as an attempt to resolve these initial problems.

\subsection{\texorpdfstring{McDuff--Wehrheim's Kuranishi atlases}{McDuff-Wehrheim\textquoteright s Kuranishi atlases}}
\label{kt23}

Next we discuss an approach to Kuranishi spaces developed by McDuff and Wehrheim \cite{McDu,McWe1,McWe2,McWe3}. Their main definition is that of a ({\it weak\/}) {\it Kuranishi atlas\/} on a topological space $X$. It is a variation on the notion of {\it good coordinate system\/} in the work of Fukaya, Oh, Ohta and Ono, as in \cite[Def.~6.1]{FuOn}, \cite[Lem.~A1.11]{FOOO1}, \cite[\S 15]{FOOO2}, and \cite[\S 5]{FOOO4}. Here are~\cite[Def.s 2.2.2 \& 2.2.8]{McWe2}.

\begin{dfn} An {\it MW Kuranishi neighbourhood\/} $(V,E,\Ga,s,\psi)$ on a topological space $X$ is the same as a FOOO Kuranishi neighbourhood in Definition \ref{kt2def1}, except that $\Ga$ is not required to act effectively on~$V$.
\label{kt2def6}
\end{dfn}

\begin{dfn} Suppose $(V_B,E_B,\Ga_B,s_B,\psi_B)$, $(V_C,E_C,\Ga_C,s_C,\psi_C)$ are MW Kuranishi neighbourhoods on a topological space $X$, and $S\subseteq \Im\psi_B\cap\Im\psi_C\subseteq X$ is open. We say a quadruple $\Phi_{BC}=(\ti V_{BC},\rho_{BC},\varpi_{BC},\hat\vp_{BC})$ is an {\it MW coordinate change from\/ $(V_B,E_B,\Ga_B,s_B,\psi_B)$ to $(V_C,\ab E_C,\ab\Ga_C,\ab s_C,\ab\psi_C)$ over\/} $S$ if:
\begin{itemize}
\setlength{\itemsep}{0pt}
\setlength{\parsep}{0pt}
\item[(a)] $\ti V_{BC}$ is a $\Ga_C$-invariant embedded submanifold of $V_C$ containing $\bar\psi_C^{-1}(S)$.
\item[(b)] $\rho_{BC}\colon \Ga_C\ra\Ga_B$ is a surjective group morphism, with kernel $\De_{BC}\subseteq\Ga_C$.

There should exist an isomorphism $\Ga_C\cong\Ga_B\t\De_{BC}$ identifying $\rho_{BC}$ with the projection $\Ga_B\t\De_{BC}\ra\Ga_B$.
\item[(c)] $\varpi_{BC}\colon \ti V_{BC}\ra V_B$ is a $\rho_{BC}$-equivariant \'etale map, with image $V_{BC}=\varpi_{BC}(\ti V_{BC})$ a $\Ga_B$-invariant open neighbourhood of $\bar\psi_B^{-1}(S)$ in $V_B$, such that $\varpi_{BC}\colon \ti V_{BC}\ra V_{BC}$ is a principal $\De_{BC}$-bundle.
\item[(d)] $\hat\vp_{BC}\colon E_B\ra E_C$ is an injective $\Ga_C$-equivariant linear map, where the $\Ga_C$-action on $E_B$ is induced from the $\Ga_B$-action by $\rho_{BC}$, so in particular $\De_{BC}$ acts trivially on~$E_B$.
\item[(e)] $\hat\vp_{BC}\ci s_B\ci\varpi_{BC}=s_C\vert_{\ti V_{BC}}\colon \ti V_{BC}\ra E_C$.
\item[(f)] $\psi_B\ci(\varpi_{BC})_*=\psi_C$ on $(s_C^{-1}(0)\cap\ti V_{BC})/\Ga_C$.
\item[(g)] For each $v\in\ti V_{BC}$ we have a commutative diagram
\e
\begin{gathered}
\xymatrix@C=25pt@R=13pt{ 0 \ar[r] & T_v\ti V_{BC} \ar[rr]_(0.6){\subset} \ar[d]^{\d(\varpi_{BC}^*(s_B))\vert_v} && T_vV_C \ar[r] \ar[d]^{\d s_C\vert_v} & N_{BC}\vert_v \ar[r] \ar@{.>}[d]^{\d_{\rm fibre}s_C\vert_v} & 0 \\
0 \ar[r] & E_B \ar[rr]^(0.65){\hat\vp_{BC}} && E_C \ar[r] & E_C/\hat\vp_{BC}(E_B) \ar[r] & 0 }
\end{gathered}
\label{kt2eq5}
\e
with exact rows, where $N_{BC}$ is the normal bundle of $\ti V_{BC}$ in $V_C$. We require the induced morphism $\d_{\rm fibre}s_C\vert_v$ in \eq{kt2eq5} to be an isomorphism.
\end{itemize}
\label{kt2def7}
\end{dfn}
 
\begin{dfn} Let $X$ be a compact, metrizable topological space. An {\it MW weak Kuranishi atlas\/ $\cK=\bigl(A,I,(V_B,E_B,\Ga_B,s_B,\psi_B)_{B\in I},\Phi_{BC,\;B,C\in I,\; B\subsetneq C}\bigr)$ on\/ $X$ of virtual dimension\/} $n\in\Z$, as in \cite[Def.~2.3.1]{McWe2}, consists of a finite indexing set $A$, a set $I$ of nonempty subsets of $A$, MW Kuranishi neighbourhoods $(V_B,E_B,\Ga_B,s_B,\psi_B)$ on $X$ for all $B\in I$ with $\dim V_B-\rank E_B=n$ and $X=\bigcup_{B\in I}\Im\psi_B$, and MW coordinate changes $\Phi_{BC}=(\ti V_{BC},\rho_{BC},\varpi_{BC},\hat\vp_{BC})$ from $(V_B,\ab E_B,\ab\Ga_B,\ab s_B,\ab\psi_B)$ to $(V_C,E_C,\Ga_C,s_C,\psi_C)$ on $S=\Im\psi_B\cap\Im\psi_C$ for all $B,C\in I$ with $B\subsetneq C$, satisfying the four conditions:
\begin{itemize}
\setlength{\itemsep}{0pt}
\setlength{\parsep}{0pt}
\item[(a)] We have $\{a\}\in I$ for all $a\in A$, and $I=\bigl\{\es\ne B\subseteq A\colon  \bigcap_{a\in B}\Im\psi_{\{a\}}\ne\es\bigr\}$. Also
$\Im\psi_B=\bigcap_{a\in B}\Im\psi_{\{a\}}$ for all $B\in I$.
\item[(b)] We have $\Ga_B=\prod_{a\in B}\Ga_{\{a\}}$ for all $B\in I$. If $B,C\in I$ with $B\subsetneq C$ then $\rho_{BC}\colon \Ga_C\ra\Ga_B$ is the obvious projection $\prod_{a\in C}\Ga_{\{a\}}\ra \prod_{a\in B}\Ga_{\{a\}}$, with kernel $\De_{BC}\cong \prod_{a\in C\sm B}\Ga_{\{a\}}$.
\item[(c)] We have $E_B=\prod_{a\in B}E_{\{a\}}$ for all $B\in I$, with the obvious representation of $\Ga_B=\prod_{a\in B}\Ga_{\{a\}}$. If $B\subsetneq C$ in $I$ then $\hat\vp_{BC}\colon E_B=\prod_{a\in B}E_{\{a\}}\ra E_C=\prod_{a\in C}E_{\{a\}}$ is $\id_{E_{\{a\}}}$ for $a\in B$, and maps to zero in $E_{\{a\}}$ for $a\in C\sm B$.
\item[(d)] If $B,C,D\in I$ with $B\subsetneq C\subsetneq D$ then $\varpi_{BC}\ci\varpi_{CD}=\varpi_{BD}$ on $\ti V_{BCD}:=\ti V_{BD}\cap \varpi_{CD}^{-1}(\ti V_{BC})$. One can show using (b),(c) and Definition \ref{kt2def7} that $\ti V_{BD}$ and $\varpi_{CD}^{-1}(\ti V_{BC})$ are both open subsets in $s_D^{-1}(\hat\vp_{BD}(E_B))$, which is a submanifold of $V_D$, so $\ti V_{BCD}$ is a submanifold of~$V_D$.
\end{itemize}

We call $\cK=\bigl(A,I,(V_B,E_B,\Ga_B,s_B,\psi_B)_{B\in I},\Phi_{BC,\;B\subsetneq C}\bigr)$ an {\it MW Kuranishi atlas on\/} $X$, as in \cite[Def.~2.3.1]{McWe2}, if it also satisfies:
\begin{itemize}
\setlength{\itemsep}{0pt}
\setlength{\parsep}{0pt}
\item[(e)] If $B,C,D\in I$ with $B\subsetneq C\subsetneq D$ then $\varpi_{CD}^{-1}(\ti V_{BC})\subseteq\ti V_{BD}$.
\end{itemize}

McDuff and Wehrheim also define {\it orientations\/} on MW weak Kuranishi atlases, in a very similar way to Definition~\ref{kt2def4}.

Two MW weak Kuranishi atlases $\cK,\cK'$ on $X$ are called {\it directly commensurate\/} if they are both contained in a third MW weak Kuranishi atlas $\cK''$. They are called {\it commensurate\/} if there exist MW weak Kuranishi atlases $\cK=\cK_0,\cK_1,\ldots,\cK_m=\cK'$ with $\cK_{i-1},\cK_i$ directly commensurate for $i=1,\ldots,m$. This is an equivalence relation on MW weak Kuranishi atlases on~$X$.
\label{kt2def8}
\end{dfn}

McDuff and Wehrheim argue that their concept of MW weak Kuranishi atlas is a more natural, or more basic, idea than a FOOO Kuranishi space, since in analytic moduli problems such as $J$-holomorphic curve moduli spaces, one has to construct an MW weak Kuranishi atlas (or something close to it) first, and then define the FOOO Kuranishi structure using this.

\subsection{How MW Kuranishi atlases are used}
\label{kt24}

McDuff and Wehrheim \cite{McWe1,McWe2,McWe3} prove analogues of Theorems \ref{kt2thm1} and \ref{kt2thm2}:

\begin{thm} Let\/ $(S,\om)$ be a symplectic manifold with tame almost complex structure $J,$ and\/ $\cM(\be,J)$ a compact moduli space of simple $J$-holomorphic maps $u\colon \CP^1\ra S$ in homology class\/ $\be\in H_2(S;\Z)$ with one marked point, modulo re\-para\-met\-riz\-a\-tion. Then one can construct an oriented MW weak Kuranishi atlas $\cK$ without boundary on $\cM(\be,J)$. 

The construction depends on many arbitrary choices, but any two such atlases $\cK,\cK'$ resulting from different choices are commensurate.
\label{kt2thm3}
\end{thm}

McDuff and Wehrheim also announce an extension of Theorem \ref{kt2thm3} to moduli spaces $\oM(\be,J)$ of genus zero prestable $J$-holomorphic maps $u\colon \Si\ra S$, without assuming $\Si$ nonsingular or $u$ simple.

\begin{thm} Let\/ $\cK$ be an oriented MW weak Kuranishi atlas without boundary of dimension $n$ on a compact, metrizable topological space $X$. Then $\cK$ determines a \begin{bfseries}virtual fundamental class\end{bfseries} $[[X]]_{\rm virt}$ in $\check H_n(X;\Q),$ where $\check H_*(-;\Q)$ is \v Cech homology over\/ $\Q$. Any two commensurate MW weak Kuranishi atlases $\cK,\cK'$ on $X$ yield the same virtual fundamental class.
\label{kt2thm4}
\end{thm}

One should also expect that many other moduli spaces of $J$-holomorphic curves admit MW weak Kuranishi atlases, and that if $X$ is compact with an oriented MW weak Kuranishi atlas $\cK$ with boundary or corners then $(X,\cK)$ has (choices of) virtual chains $[X]_{\rm virt}$, but these results are not yet available.

\subsection{\texorpdfstring{Dingyu Yang's Kuranishi structures, and Hofer--Wysocki--Zehnder's polyfolds}{Dingyu Yang\textquoteright s Kuranishi structures, and Hofer-Wysocki-Zehnder\textquoteright s polyfolds}}
\label{kt25}

As part of a project to define a truncation functor from polyfolds to Kuranishi spaces, Dingyu Yang \cite{Yang1,Yang2,Yang3} writes down his own theory of Kuranishi spaces:

\begin{dfn} Let $X$ be a compact, metrizable topological space. A {\it DY Kuranishi structure\/} $\cK$ on $X$ is a FOOO Kuranishi structure in the sense of Definition \ref{kt2def3}, satisfying two additional conditions \cite[Def.~1.11]{Yang2} the {\it maximality condition\/} and the {\it topological matching condition}, which are designed to ensure that the $V_p/\Ga_p$ in $\cK$ for all $p\in X$ can be glued nicely using the $h_{qp},\vp_{qp}$ in $\cK$ for all $p,q\in X$ to make a Hausdorff topological space. There are a few other small differences --- for instance, Yang does not require vector bundles $E_p$ in $(V_p,E_p,\Ga_p,s_p,\psi_p)$ to be trivial.
\label{kt2def9}
\end{dfn}

The next definition comes from Yang \cite[\S 1.6]{Yang1}, \cite[\S 5]{Yang2}, \cite[\S 2.4]{Yang3}.

\begin{dfn} Let $\cK,\cK'$ be DY Kuranishi structures on a compact topological space $X$. An {\it embedding\/} $\ep\colon \cK\hookra\cK'$ is a choice of FOOO coordinate change $\ep_p\colon (V_p,E_p,\Ga_p,s_p,\psi_p)\ra(V_p',E_p',\Ga_p',s_p',\psi_p')$ with domain $V_p$ for all $p\in X$, commuting with the FOOO coordinate changes $\Phi_{qp},\Phi'_{qp}$ in $\cK,\cK'$ up to elements of $\Ga_p'$. An embedding is a {\it chart refinement\/} if the $\ep_p$ come from inclusions of $\Ga_p$-invariant open sets~$V_p\hookra V_p'$.

DY Kuranishi structures $\cK,\cK'$ on $X$ are called {\it R-equivalent\/} (or {\it equivalent\/}) if there is a diagram of DY Kuranishi structures on $X$
\begin{equation*}
\xymatrix@C=30pt@R=15pt{\cK  & \cK_1 \ar[l]_\sim \ar@{=>}[r] & \cK_2 & \cK_3 \ar@{=>}[l] \ar[r]^\sim & \cK', }
\end{equation*}
where arrows $\Longra$ are embeddings, and ${\buildrel\sim\over\longra}$ are chart refinements. Using facts about existence of good coordinate systems, Yang proves \cite[Th.~1.6.17]{Yang1}, \cite[\S 11.2]{Yang2} that R-equivalence is an equivalence relation on DY Kuranishi structures.

\label{kt2def10}
\end{dfn}

Yang emphasizes the idea, which he calls {\it choice independence}, that when one constructs a (DY) Kuranishi structure $\cK$ on a moduli space $\oM$, it should be independent of choices up to R-equivalence. 

One goal of Yang's work is to relate the Kuranishi space theory of Fukaya, Oh, Ohta and Ono \cite{Fuka,FOOO1,FOOO2,FOOO3,FOOO4,FOOO5,FOOO6,FOOO7,FOOO8,FOOO9,FuOn} to the polyfold theory of Hofer, Wysocki and Zehnder \cite{Hofe,HWZ1,HWZ2,HWZ3,HWZ4,HWZ5,HWZ6,HWZ7}. Here is a very brief introduction to this:
\begin{itemize}
\setlength{\itemsep}{0pt}
\setlength{\parsep}{0pt}
\item An {\it sc-Banach space\/} $\cV$ is a sequence $\cV=(\cV_0\supset\cV_1\supset\cV_2\supset\cdots),$ where the $\cV_i$ are Banach spaces, the inclusions $\cV_{i+1}\hookra\cV_i$ are compact, bounded linear maps, and $\cV_\iy=\bigcap_{i\ge 0}\cV_i$ is dense in every~$\cV_i$. 

The {\it tangent space\/} $T\cV$ is $T\cV=(\cV_1\op\cV_0\supset\cV_2\op\cV_1\supset\cdots)$, an sc-Banach space. An {\it open set\/} $\cQ$ in $\cV$ is an open set $\cQ\subset\cV_0$, and we write $\cQ_i=\cQ\cap\cV_i$ for $i\ge 0$. Its {\it tangent space\/} is $T\cQ=\cQ_1\op\cV_0$, as an open set in~$T\cV$.

An example to bear in mind is if $M$ is a compact manifold, $E\ra M$ a smooth vector bundle, $\al\in(0,1)$, and $\cV_k=C^{k,\al}(E)$ for $k=0,1,\ldots.$ 
\item Let $\cV=(\cV_0\supset\cV_1\supset\cdots)$, $\cW=(\cW_0\supset\cW_1\supset\cdots)$ be sc-Banach spaces and $\cQ\subseteq\cV$, $\cR\subseteq\cW$ be open. A map $f\colon \cQ\ra\cR$ is called {\it sc}${}^0$ if $f(\cQ_i)\subseteq\cR_i$ and $f\vert_{\cQ_i}\colon \cQ_i\ra\cR_i$ is a continuous map of Banach manifolds for all~$i\ge 0$.

An sc${}^0$ map $f\colon \cQ\ra\cR$ is called {\it sc}${}^1$ if for each $q\in\cQ_1$ there exists a bounded linear map $Df_q\colon \cV_0\ra\cW_0$, such that $f\vert_{\cQ_1}\colon \cQ_1\ra\cR_0$ is a $C^1$ map of Banach manifolds with $\nabla f\vert_q=Df_q\vert_{\cV_1}\colon \cV_1\ra\cW_0$ for all $q\in\cQ_1$, and $Tf\colon T\cQ\ra T\cR$ mapping $Tf\colon (q,v)\mapsto(f(q),Df_q(v))$ is an sc${}^0$ map.

By induction on $k$, we call $f\colon \cQ\ra\cR$ an {\it sc}${}^k$ map for $k=2,3,\ldots$ if $f$ is sc${}^1$ and $Tf\colon T\cQ\ra T\cR$ is an sc${}^{k-1}$ map. We call $f\colon \cQ\ra\cR$ {\it sc-smooth}, or {\it sc}${}^\iy$, if it is sc${}^k$ for all $k=0,1,\ldots.$ This implies that $f\vert_{\cQ_{i+k}}\colon \cQ_{i+k}\ra\cR_i$ is a $C^k$-map of Banach manifolds for all~$i,k\ge 0$.
\item Let $\cV=(\cV_0\supset\cV_1\supset\cdots)$ be an sc-Banach space and $\cQ\subseteq\cV$ be open. An {\it sc\/$^\iy$-retraction\/} is an sc-smooth map $r\colon \cQ\ra\cQ$ with $r\ci r=r$. Set $\O=\Im r\subset\cV$. We call $(\O,\cV)$ a {\it local sc-model}.

If $\cV$ is finite-dimensional then $\O$ is just a smooth manifold. But in infinite dimensions, new phenomena occur, and the tangent spaces $T_x\O$ can vary discontinuously with $x\in\O$. This is important for `gluing'.
\item An {\it M-polyfold chart\/} $(\O,\cV,\psi)$ on a topological space $Z$ is a local sc-model $(\O,\cV)$ and a homeomorphism $\psi\colon \O\ra\Im\psi$ with an open set~$\Im\psi\subset Z$. 
\item M-polyfold charts $(\O,\cV,\psi),(\ti\O,\ti\cV,\ti\psi)$ on $Z$ are {\it compatible\/} if $\ti\psi^{-1}\ci\psi\ci r\colon \cQ\ra\ti\cV$ and $\psi^{-1}\ci\ti\psi\ci\ti r\colon \ti\cQ\ra\cV$ are sc-smooth, where $\cQ\subset\cV$, $\ti\cQ\subset\ti\cV$ are open and $r\colon \cQ\ra\cQ$, $\ti r\colon \ti\cQ\ra\ti\cQ$ are sc-smooth with $r\ci r=r$, $\ti r\ci\ti r=\ti r$ and $\Im r=\psi^{-1}(\Im\ti\psi)\subseteq\O$, $\Im\ti r=\ti\psi^{-1}(\Im\psi)\subseteq\ti\O$.
\item An {\it M-polyfold\/} is roughly a metrizable topological space $Z$ with a maximal atlas of pairwise compatible M-polyfold charts.
\item {\it Polyfolds\/} are the orbifold version of M-polyfolds, proper \'etale groupoids in M-polyfolds.
\item A {\it polyfold Fredholm structure\/} $\cP$ on a metrizable topological space $X$ writes $X$ as the zeroes of an sc-Fredholm section ${\mathfrak s}\colon \fV\ra\fE$ of a strong polyfold vector bundle $\fE\ra\fV$ over a polyfold $\fV$. 
\end{itemize}
This is deep mathematics, and all rather complicated. The motivation for local sc-models $(\O,\cV)$ is that they can be used to describe functional-analytic problems involving `gluing', `bubbling', and `neck-stretching', including moduli spaces of $J$-holomorphic curves with singularities of various kinds.

Yang proves \cite[Th.~3.1.7]{Yang1} (see also \cite[\S 2.6]{Yang3}):

\begin{thm} Suppose we are given a `polyfold Fredholm structure'\/ $\cP$ on a compact metrizable topological space $X,$ that is, we write $X$ as the zeroes of an sc-Fredholm section ${\mathfrak s}\colon \fV\ra\fE$ of a strong polyfold vector bundle $\fE\ra\fV$ over a polyfold\/ $\fV,$ where $\mathfrak s$ has constant Fredholm index $n\in\Z$. Then we can construct a DY Kuranishi structure $\cK$ on $X,$ of virtual dimension $n,$ which is independent of choices up to R-equivalence.
\label{kt2thm5}
\end{thm}

\subsection{How polyfolds are used}
\label{kt26}

Hofer, Wysocki and Zehnder's polyfold programme \cite{Hofe,HWZ1,HWZ2,HWZ3,HWZ4,HWZ5,HWZ6,HWZ7} aims to show that moduli spaces of $J$-holomorphic curves in symplectic geometry may be given a polyfold Fredholm structure, and that compact spaces with oriented polyfold Fredholm structures have virtual chains and virtual classes. They will then use these virtual chains/classes to define big theories in symplectic geometry, such as Gromov--Witten invariants \cite{FuOn,HWZ6} or Symplectic Field Theory \cite{EGH}. 

Here is an analogue of Theorems \ref{kt2thm1}(a) and \ref{kt2thm3} in polyfold theory, proved by Hofer, Wysocki and Zehnder~\cite{HWZ6}.

\begin{thm} Moduli spaces $\oM_{g,m}(J,\be)$ of stable maps $u\colon \Si\ra S$ from a prestable closed Riemann surface $\Si$ of genus $g$ with\/ $m$ marked points, to a compact symplectic manifold $(S,\om)$ of dimension $2n$ with compatible almost complex structure $J,$ with homology class\/ $u_*([\Si])=\be\in H_2(S;\Z),$ may be given a `polyfold Fredholm structure'. That is, we may write $\oM_{g,m}(J,\be)$ as the zeroes of an sc-Fredholm section ${\mathfrak s}\colon \fV\ra\fE$ of a strong polyfold vector bundle $\fE\ra\fV$ over a polyfold $\fV,$ where $\mathfrak s$ has Fredholm index $2\bigl(c_1(S)\cdot\be+(n-3)(1-g)+m\bigr)$. 

The `evaluation maps' at marked points $\ev_i\colon \oM_{g,m}(J,\be)\ra S$ lift to sc-smooth maps of polyfolds $\mathfrak{ev}_i\colon \fV\ra S$. The polyfold Fredholm structure $(\fV,\fE,{\mathfrak s})$ has a natural `orientation'.
\label{kt2thm6}
\end{thm}

Hofer, Wysocki and Zehnder also announce the proofs of existence of polyfold Fredholm structures on other moduli spaces of $J$-holomorphic curves, in particular those relevant to Symplectic Field Theory \cite{EGH}. Constructing a polyfold Fredholm structure on a moduli space of $J$-holomorphic curves involves far fewer arbitrary choices than defining a Kuranishi structure.

Combining Theorems \ref{kt2thm5}--\ref{kt2thm6} gives a DY Kuranishi structure on $\oM_{g,m}(J,\be)$, uniquely up to R-equivalence. This is an example of a FOOO Kuranishi structure, giving an alternative proof of Theorem~\ref{kt2thm1}(a). 

For the polyfold analogue of Theorems \ref{kt2thm2} and \ref{kt2thm4}, Hofer et al.~\cite{HWZ4} prove:

\begin{thm} Let\/ $(\fV,\fE,{\mathfrak s})$ be an oriented polyfold Fredholm structure with corners of virtual dimension $k$ on a compact topological space $X$. Then there exist small perturbations $\ti{\mathfrak s}$ of $\mathfrak s$ as an `sc-smooth multisection', such that\/ $\ti X=\ti{\mathfrak s}^{-1}(0)$ is a compact, oriented, branched, $\Q$-weighted\/ $k$-orbifold with corners.  

We may triangulate $\ti X$ by $\Q$-weighted\/ $k$-simplices to define a \begin{bfseries}virtual class\end{bfseries} $[[X]]_{\rm virt}$ in the relative singular homology $H^{\rm si}_k(\fV,\pd\fV;\Q)$. Alternatively, using relative de Rham cohomology we may define $[[X]]^{\rm virt}\in H_{\rm dR}^k(\fV,\pd\fV;\R)^*$ mapping $[\om]\mapsto\int_{\ti X}\om$. Both of\/ $[[X]]_{\rm virt},[[X]]^{\rm virt}$ are independent of the choice of\/ $\ti{\mathfrak s},\ti X$. If\/ $\pd\fV=\es$ then $\pd\ti X=\es,$ $[[X]]_{\rm virt}\in H^{\rm si}_k(\fV;\Q),$ and\/ $[[X]]^{\rm virt}\in H_{\rm dR}^k(\fV;\R)^*$.
\label{kt2thm7}
\end{thm}

\section{Kuranishi neighbourhoods as a 2-category}
\label{kt3}

FOOO Kuranishi spaces and MW Kuranishi atlases in \S\ref{kt2} are each built out of Kuranishi neighbourhoods $(V,E,\Ga,s,\psi)$ on a topological space $X$, and `coordinate changes' $\Phi_{ij}\colon (V_i,E_i,\Ga_i,s_i,\psi_i)\ra (V_j,E_j,\Ga_j,s_j,\psi_j)$ between them.

Our Kuranishi spaces $\bX$ will be built out of Kuranishi neighbourhoods $(V,\ab E,\ab\Ga,\ab s,\ab\psi)$ on $X$, and 1-morphisms $\Phi_{ij}$ between them (which include FOOO and MW coordinate changes as special cases), and 2-morphisms $\La_{ij}\colon \Phi_{ij}\Ra\Phi_{ij}'$ between 1-morphisms, so that Kuranishi neighbourhoods form a 2-category.

This section explains our notions of 1- and 2-morphisms of Kuranishi neighbourhoods, and their properties. Section \ref{kt4} will then use these to define the 2-category of Kuranishi spaces. With the exception of \S\ref{kt36}, all of this section, including proofs of quoted results, comes from~\cite{Joyc5,Joyc6}.

\subsection{Kuranishi neighbourhoods, 1-morphisms, and 2-morphisms}
\label{kt31}

We will use following `$O(s)$' and `$O(s^2)$' notation very often:

\begin{dfn} Let $V$ be a manifold, $E\ra V$ a vector bundle, and $s\in C^\iy(E)$ a smooth section.
\begin{itemize}
\setlength{\itemsep}{0pt}
\setlength{\parsep}{0pt}
\item[(i)] If $F\ra V$ is another vector bundle and $t_1,t_2\in C^\iy(F)$ are smooth sections, we write $t_1=t_2+O(s)$ if there exists $\al\in C^\iy(E^*\ot
F)$ such that $t_1=t_2+\al\cdot s$ in $C^\iy(F)$, where the
contraction $\al\cdot s$ is formed using the natural pairing of
vector bundles $(E^*\ot F)\t E\ra F$ over $V$.
\item[(ii)] We write $t_1=t_2+O(s^2)$ if there exists $\al\in
C^\iy(E^*\ot E^*\ot F)$ such that $t_1=t_2+\al\cdot(s\ot s)$ in
$C^\iy(F)$, where $\al\cdot(s\ot s)$ uses the pairing~$(E^*\ot
E^*\ot F)\t(E\ot E)\ra F$.
\end{itemize}

Now let $W$ be another manifold, and $f,g\colon V\ra W$ be smooth maps. 
\begin{itemize}
\setlength{\itemsep}{0pt}
\setlength{\parsep}{0pt}
\item[(iii)] We write $f=g+O(s)$ if whenever $h\colon W\ra\R$ is a smooth map, there
exists $\al\in C^\iy(E^*)$ such that $h\ci f=h\ci g+\al\cdot s$.
\item[(iv)] We write $f=g+O(s^2)$ if whenever $h\colon W\ra\R$ is a smooth
map, there exists $\al\in C^\iy(E^*\ot E^*)$ such that $h\ci f=h\ci
g+\al\cdot (s\ot s)$.
\item[(v)] If $\La\in C^\iy\bigl(E^*\ot f^*(TW)\bigr)$, we write $f=g+\La\cdot s+O(s^2)$ if whenever $h\colon W\ra\R$ is a smooth map, there exists $\al\in C^\iy(E^*\ot E^*)$ such that $h\ci f=h\ci g+\La\cdot(s\ot f^*(\d h))+\al\cdot (s\ot s)$. Here $s\ot f^*(\d h)$ lies in $C^\iy\bigl(E\ot f^*(T^*W)\bigr)$, and so has an obvious pairing with $\La$.
\end{itemize}

Next suppose $f,g\colon V\ra W$ with $f=g+O(s)$, and $F\ra V$, $G\ra W$ are vector bundles, and $t_1\in C^\iy(F\ot f^*(G))$, $t_2\in C^\iy(F\ot g^*(G))$. We wish to compare $t_1,t_2$, even though they are sections of different vector bundles.
\begin{itemize}
\setlength{\itemsep}{0pt}
\setlength{\parsep}{0pt}
\item[(vi)] We write $t_1=t_2+O(s)$ if for all $\be\in C^\iy(G^*)$ we have $t_1\cdot f^*(\be)=t_2\cdot g^*(\be)+O(s)$ in sections of $F\ra V$, as in (i).
\end{itemize}
Given any $t_1\in C^\iy(F\ot f^*(G))$, there exists $t_2\in C^\iy(F\ot g^*(G))$ with $t_1=t_2+O(s)$ in the sense of (vi), and if $t_2'$ is an alternative choice then $t_2'=t_2+O(s)$ in the sense of (i). To prove this, first suppose $G$ is trivial, so that $f^*(G)=g^*(G)$, and $t_2=t_1$ is a possible choice. In general we can trivialize $G$ locally on $Y$, and combine the corresponding local choices for $t_2$ on $X$ by a partition of unity. 

The moral is that if $f=g+O(s)$ and we are interested in smooth sections $t$ of $F\ot f^*(G)$ up to $O(s)$, then we can treat the vector bundles $F\ot f^*(G)$ and $F\ot g^*(G)$ as essentially the same.

If instead $f,g\colon V\ra W$ with $f=g+O(s^2)$ and $F,G,t_1,t_2$ are as above\begin{itemize}
\setlength{\itemsep}{0pt}
\setlength{\parsep}{0pt}
\item[(vii)] We write $t_1=t_2+O(s^2)$ if for all $\be\in C^\iy(G^*)$ we have $t_1\cdot f^*(\be)=t_2\cdot g^*(\be)+O(s^2)$ in sections of $F\ra V$, as in (ii).
\end{itemize}
Given any $t_1\in C^\iy(F\ot f^*(G))$, there exists $t_2\in C^\iy(F\ot g^*(G))$ with $t_1=t_2+O(s^2)$ in the sense of (vii), and if $t_2'$ is an alternative choice then $t_2'=t_2+O(s^2)$ in the sense of~(ii). 
\label{kt3def1}
\end{dfn}

Here is our notion of Kuranishi neighbourhood. It is very similar to FOOO and MW Kuranishi neighbourhoods in Definitions \ref{kt2def1} and \ref{kt2def6}, but it is more general, in that we allow $E$ to be a nontrivial vector bundle.

\begin{dfn} Let $X$ be a topological space. A {\it Kuranishi neighbourhood\/} on $X$ is a quintuple $(V,E,\Ga,s,\psi)$ such that:
\begin{itemize}
\setlength{\itemsep}{0pt}
\setlength{\parsep}{0pt}
\item[(a)] $V$ is a smooth manifold.
\item[(b)] $\pi\colon E\ra V$ is a real vector bundle over $V$, called the {\it obstruction bundle}.
\item[(c)] $\Ga$ is a finite group with a smooth action on $V$ (not necessarily effective), and a compatible action on $E$ preserving the vector bundle structure.
\item[(d)] $s\colon V\ra E$ is a $\Ga$-equivariant smooth section of $E$, called the {\it Kuranishi section}.
\item[(e)] $\psi$ is a homeomorphism from $s^{-1}(0)/\Ga$ to an open subset $\Im\psi=\bigl\{\psi(\Ga v)\colon\ab v\in s^{-1}(0)\bigr\}$ in $X$, called the {\it footprint\/} of~$(V,E,\Ga,s,\psi)$.
\end{itemize}
We will write $\bar\psi\colon s^{-1}(0)\ra\Im\psi\subseteq X$ for the composition of $\psi$ with the projection $s^{-1}(0)\ra s^{-1}(0)/\Ga$.

We call $(V,E,\Ga,s,\psi)$ a {\it global\/} Kuranishi neighbourhood on $X$ if~$\Im\psi=X$.
\label{kt3def2}
\end{dfn}

The next two definitions are crucial to our programme.

\begin{dfn} Let $X,Y$ be topological spaces, $f\colon X\ra Y$ a continuous map, $(V_i,E_i,\Ga_i,s_i,\psi_i)$, $(V_j,E_j,\Ga_j,s_j,\psi_j)$ be Kuranishi neighbourhoods on $X,Y$ respectively, and $S\subseteq\Im\psi_i\cap f^{-1}(\Im\psi_j)\subseteq X$ be an open set. A 1-{\it morphism $\Phi_{ij}=(P_{ij},\pi_{ij},\phi_{ij},\hat\phi_{ij})\colon (V_i,E_i,\Ga_i,s_i,\psi_i)\ra (V_j,E_j,\Ga_j,s_j,\psi_j)$ of Kuranishi neighbourhoods over\/} $(S,f)$ is a quadruple $(P_{ij},\pi_{ij},\phi_{ij},\hat\phi_{ij})$ satisfying 
\begin{itemize}
\setlength{\itemsep}{0pt}
\setlength{\parsep}{0pt}
\item[(a)] $P_{ij}$ is a manifold, with commuting smooth actions of $\Ga_i,\Ga_j$ (that is, with a smooth action of $\Ga_i\t\Ga_j$), with the $\Ga_j$-action free.
\item[(b)] $\pi_{ij}\colon P_{ij}\ra V_i$ is a smooth map which is $\Ga_i$-equivariant, $\Ga_j$-invariant, and \'etale (a local diffeomorphism). The image $V_{ij}:=\pi_{ij}(P_{ij})$ is a $\Ga_i$-invariant open neighbourhood of $\bar\psi_i^{-1}(S)$ in $V_i$, and the fibres $\pi_{ij}^{-1}(v)$ of $\pi_{ij}$ for $v\in V_{ij}$ are $\Ga_j$-orbits, so that $\pi_{ij}\colon P_{ij}\ra V_{ij}$ is a principal $\Ga_j$-bundle.

We do not require $\bar\psi_i^{-1}(S)=V_{ij}\cap s_i^{-1}(0)$, only that~$\bar\psi_i^{-1}(S)\subseteq V_{ij}\cap s_i^{-1}(0)$.
\item[(c)] $\phi_{ij}\colon P_{ij}\ra V_j$ is a $\Ga_i$-invariant and $\Ga_j$-equivariant smooth map, that is, $\phi_{ij}(\ga_i\cdot p)=\phi_{ij}(p)$, $\phi_{ij}(\ga_j\cdot p)=\ga_j\cdot\phi_{ij}(p)$ for all $\ga_i\in\Ga_i$, $\ga_j\in\Ga_j$, $p\in P_{ij}$.
\item[(d)] $\hat\phi_{ij}\colon \pi_{ij}^*(E_i)\ra\phi_{ij}^*(E_j)$ is a $\Ga_i$- and $\Ga_j$-equivariant morphism of vector bundles on $P_{ij}$, where the $\Ga_i,\Ga_j$-actions are induced by the given $\Ga_i$-action and the trivial $\Ga_j$-action on $E_i$, and vice versa for $E_j$.
\item[(e)] $\hat\phi_{ij}(\pi_{ij}^*(s_i))=\phi_{ij}^*(s_j)+O(\pi_{ij}^*(s_i)^2)$, in the sense of Definition~\ref{kt3def1}.
\item[(f)] $f\ci\bar\psi_i\ci\pi_{ij}=\bar\psi_j\ci\phi_{ij}$ on~$\pi_{ij}^{-1}(s_i^{-1}(0))\subseteq P_{ij}$.
\end{itemize}

If $Y=X$ and $f=\id_X$ then we call $\Phi_{ij}$ a 1-{\it morphism of Kuranishi neighbourhoods over\/} $S$, or just a 1-{\it morphism over\/}~$S$.

Let $(V_i,E_i,\Ga_i,s_i,\psi_i)$ be a Kuranishi neighbourhood on $X$, and $S\subseteq\Im\psi_i$ be open. We will define the {\it identity\/ $1$-morphism\/} over $S$
\e
\id_{(V_i,E_i,\Ga_i,s_i,\psi_i)}=(P_{ii},\pi_{ii},\phi_{ii},\hat\phi_{ii})\colon (V_i,E_i,\Ga_i,s_i,\psi_i)\ra (V_i,E_i,\Ga_i,s_i,\psi_i).
\label{kt3eq1}
\e
Since $P_{ii}$ must have two different actions of $\Ga_i$, for clarity we write $\Ga_i^d=\Ga_i^t=\Ga_i$, where $\Ga_i^d$ and $\Ga_i^t$ mean the copies of $\Ga_i$ acting on the domain and target of the 1-morphism in \eq{kt3eq1}, respectively.
 
Define $P_{ii}=V_i\t\Ga_i$, and let $\Ga_i^d$ act on $P_{ii}$ by $\ga^d\colon (v,\ga)\mapsto (\ga^d\cdot v,\ga(\ga^d)^{-1})$ and $\Ga_i^t$ act on $P_{ii}$ by $\ga^t\colon (v,\ga)\mapsto(v,\ga^t\ga)$. Define $\pi_{ii},\phi_{ii}\colon P_{ii}\ra V_i$ by $\pi_{ii}\colon (v,\ga)\mapsto v$ and $\phi_{ii}\colon (v,\ga)\mapsto\ga\cdot v$. Then $\pi_{ii}$ is $\Ga_i^d$-equivariant and $\Ga_i^t$-invariant, and is a $\Ga_i^t$-principal bundle, and $\phi_{ii}$ is $\Ga_i^d$-invariant and $\Ga_i^t$-equivariant.

At $(v,\ga)\in P_{ii}$, the morphism $\hat\phi_{ii}\colon \pi_{ii}^*(E_i)\ra\phi_{ii}^*(E_i)$ must map $E_i\vert_v\ra E_i\vert_{\ga\cdot v}$. We have such a map, the lift of the $\ga$-action on $V_i$ to $E_i$. So we define $\hat\phi_{ii}$ on $V_i\t\{\ga\}\subseteq P_{ii}$ to be the lift to $E_i$ of the $\ga$-action on $V_i$, for each $\ga\in\Ga$. 
\label{kt3def3}
\end{dfn}

\begin{rem}{\bf(i)} We will use 1-morphisms of Kuranishi neighbourhoods for two purposes in \S\ref{kt41}: Kuranishi spaces $\bX$ involve `coordinate changes', which are 1-morphisms of Kuranishi neighbourhoods  over $\id_X\colon X\ra X$ satisfying an invertibility condition. And 1-morphisms of Kuranishi spaces $\bs f\colon \bX\ra\bY$ involve 1-morphisms of Kuranishi neighbourhoods over~$f\colon X\ra Y$.
\smallskip

\noindent{\bf(ii)} As we will explain in \S\ref{kt34}, FOOO coordinate changes in \S\ref{kt21} and MW coordinate changes in \S\ref{kt23} give examples of 1-morphisms of Kuranishi neighbourhoods over~$\id_X\colon X\ra X$.
\smallskip

\noindent{\bf(iii)} Definition \ref{kt3def3} does not include the analogue of \eq{kt2eq2} being exact (the `with a tangent bundle' condition in the Fukaya--Oh--Ohta--Ono theory). But as Theorem \ref{kt3thm3} shows, this will be a consequence of the definition of `coordinate change' in Definition \ref{kt3def8}, which are special 1-morphisms over~$\id_X:X\ra X$.
\smallskip

\noindent{\bf(iv)} In a 1-morphism $\Phi_{ij}=(P_{ij},\pi_{ij},\phi_{ij},\hat\phi_{ij})$, the data $P_{ij},\pi_{ij},\phi_{ij}$ are actually a standard way to write `smooth maps' between orbifolds, using `bibundles':
\begin{equation*}
(P_{ij},\pi_{ij},\phi_{ij})\colon [V_{ij}/\Ga_i]\longra[V_j/\Ga_j]
\end{equation*}
is a {\it Hilsum--Skandalis morphism\/} between the quotient orbifolds $[V_{ij}/\Ga_i],[V_j/\Ga_j]$, as in Lerman \cite[\S 3.3]{Lerm} or Henriques and Metzler~\cite{HeMe}.

Thus, we can interpret $\Phi_{ij}$ as giving the following data:
\begin{itemize}
\setlength{\itemsep}{0pt}
\setlength{\parsep}{0pt}
\item An open neighbourhood $[V_{ij}/\Ga_i]$ of $\psi_i^{-1}(S)$ in the orbifold $[V_i/\Ga_i]$;
\item A `smooth map' of orbifolds $\Phi_{ij}'=(P_{ij},\pi_{ij},\phi_{ij})\colon [V_{ij}/\Ga_i]\ra[V_j/\Ga_j]$; and
\item A morphism of orbifold vector bundles $\hat\phi_{ij}\colon E_i\vert_{[V_{ij}/\Ga_i]}\ra (\Phi_{ij}')^*(E_j)$.
\end{itemize}
\label{kt3rem1}
\end{rem}

In our theory, 2-morphisms between 1-morphisms will also be important.

\begin{dfn} Suppose $X,Y$ are topological spaces, $f\colon X\ra Y$ is a continuous map, $(V_i,E_i,\Ga_i,s_i,\psi_i)$, $(V_j,E_j,\Ga_j,s_j,\psi_j)$ are Kuranishi neighbourhoods on $X,Y$ respectively, $S\subseteq\Im\psi_i\cap f^{-1}(\Im\psi_j)\subseteq X$ is open, and $\Phi_{ij},\Phi_{ij}'\colon (V_i,\ab E_i,\ab\Ga_i,\ab s_i,\ab\psi_i)\ra (V_j,\ab E_j,\ab\Ga_j,\ab s_j,\ab\psi_j)$ are two 1-morphisms over $(S,f)$, with $\Phi_{ij}=(P_{ij},\pi_{ij},\phi_{ij},\hat\phi_{ij})$ and~$\Phi_{ij}'=(P_{ij}',\pi_{ij}',\phi_{ij}',\hat\phi_{ij}')$.

Consider triples $(\dot P_{ij},\la_{ij},\hat\la_{ij})$ satisfying:
\begin{itemize}
\setlength{\itemsep}{0pt}
\setlength{\parsep}{0pt}
\item[(a)] $\dot P_{ij}$ is a $\Ga_i$- and $\Ga_j$-invariant open neighbourhood of $\pi_{ij}^{-1}(\bar\psi_i^{-1}(S))$ in $P_{ij}$.
\item[(b)] $\la_{ij}\colon \dot P_{ij}\ra P_{ij}'$ is a $\Ga_i$- and $\Ga_j$-equivariant smooth map with $\pi_{ij}'\ci\la_{ij}=\pi_{ij}\vert_{\dot P_{ij}}$. This implies that $\la_{ij}$ is an isomorphism of principal $\Ga_j$-bundles over $\dot V_{ij}:=\pi_{ij}(\dot P_{ij})$, so $\la_{ij}$ is a diffeomorphism with a $\Ga_i$- and $\Ga_j$-invariant open set $\la_{ij}(\dot P_{ij})$ in~$P_{ij}'$.
\item[(c)] $\hat\la_{ij}\colon \pi_{ij}^*(E_i)\vert_{\dot P_{ij}}\ra \phi_{ij}^*(TV_j)\vert_{\dot P_{ij}}$ is a $\Ga_i$- and $\Ga_j$-invariant smooth morphism of vector bundles on $\dot P_{ij}$, satisfying
\e
\begin{split}
\phi_{ij}'\ci\la_{ij}&=\phi_{ij}\vert_{\dot P_{ij}}+\hat\la_{ij}\cdot \pi_{ij}^*(s_i)+O\bigl(\pi_{ij}^*(s_i)^2\bigr)\;\>\text{and}\\ 
\la_{ij}^*(\hat\phi_{ij}')&=\hat\phi_{ij}\vert_{\dot P_{ij}}+\hat\la_{ij}\cdot \phi_{ij}^*(\d s_j)+O\bigl(\pi_{ij}^*(s_i)\bigr)\;\> \text{on $\dot P_{ij}$,}
\end{split}
\label{kt3eq2}
\e
in the notation of Definition \ref{kt3def1}. Here $\d s_j$ is a shorthand for the derivative $\nabla s_j$ in $C^\iy(T^*V_j\ot E_j)$ with respect to any connection $\nabla$ on $E_j$, and \eq{kt3eq2} is independent of this choice.
\end{itemize}

Define a binary relation $\approx$ on such triples by $(\dot P_{ij},\la_{ij},\hat\la_{ij})\approx(\dot P_{ij}',\la_{ij}',\hat\la_{ij}')$ if there exists an open neighbourhood $\ddot P_{ij}$ of $\pi_{ij}^{-1}(\bar\psi_i^{-1}(S))$ in $\dot P_{ij}\cap \dot P_{ij}'$ with
\begin{equation*}
\la_{ij}\vert_{\ddot P_{ij}}=\la_{ij}'\vert_{\ddot P_{ij}}\quad\text{and}\quad \hat\la_{ij}\vert_{\ddot P_{ij}}=\hat\la_{ij}'\vert_{\ddot P_{ij}}+O\bigl(\pi_{ij}^*(s_i)\bigr)\quad\text{on $\ddot P_{ij}$.}
\end{equation*}
Then $\approx$ is an equivalence relation. We also write $\approx_S$ for $\approx$ if we wish to stress the open set $S$. Write $[\dot P_{ij},\la_{ij},\hat\la_{ij}]$ for the $\approx$-equivalence class of $(\dot P_{ij},\la_{ij},\hat\la_{ij})$. We say that $[\dot P_{ij},\la_{ij},\hat\la_{ij}]\colon \Phi_{ij}\Ra\Phi_{ij}'$ is a 2-{\it morphism of\/ $1$-morphisms of Kuranishi neighbourhoods on\/ $X$ over\/} $(S,f)$, or just a 2-{\it morphism over\/} $(S,f)$. We often write~$\La_{ij}=[\dot P_{ij},\la_{ij},\hat\la_{ij}]$.

If $Y=X$ and $f=\id_X$ then we call $\La_{ij}$ a 2-{\it morphism of Kuranishi neighbourhoods over\/} $S$, or just a 2-{\it morphism over\/}~$S$.

For a 1-morphism $\Phi_{ij}=(P_{ij},\pi_{ij},\phi_{ij},\hat\phi_{ij})$, define the {\it identity\/ $2$-morphism\/}
\begin{equation*}
\id_{\Phi_{ij}}=[P_{ij},\id_{P_{ij}},0]\colon \Phi_{ij}\Longra\Phi_{ij}.
\end{equation*}

\label{kt3def4}
\end{dfn}

\begin{rem} This definition of 2-morphism probably seems rather arbitrary. The best way the author knows to motivate it is via d-orbifolds, as we explain in \S\ref{kt36}. It is also justified by the useful properties of 2-morphisms, in particular Theorems \ref{kt3thm1}, \ref{kt3thm2} and \ref{kt3thm3} below.
\label{kt3rem2}
\end{rem}

\subsection{Making Kuranishi neighbourhoods into a 2-category}
\label{kt32}

Readers unfamiliar with 2-categories are advised to look at \S\ref{ktA1} at this point. We wish to make Kuranishi neighbourhoods on a topological space $X$ into a weak 2-category. We have already defined the objects (Kuranishi neighbourhoods), 1-morphisms, 2-morphisms, and identity 1- and 2-morphisms, in \S\ref{kt31}. But there remains quite a lot of structure to define, which we do in this section:
\begin{itemize}
\setlength{\itemsep}{0pt}
\setlength{\parsep}{0pt}
\item Composition of 1-morphisms;
\item Horizontal composition of 2-morphisms;
\item Vertical composition of 2-morphisms;
\item Coherence 2-isomorphisms $\al_{g,f,e}\colon (g\ci f)\ci e\Ra g\ci (f\ci e)$, $\be_f\colon f\ci\id_X\Ra f$, and $\ga_f\colon \id_Y\ci f\Ra f$, as in \eq{ktAeq5} and~\eq{ktAeq7}.
\end{itemize}

\begin{dfn} Let $X,Y,Z$ be topological spaces, $f\colon X\ra Y$, $g\colon Y\ra Z$ be continuous maps, $(V_i,E_i,\Ga_i,s_i,\psi_i),(V_j,E_j,\Ga_j,s_j,\psi_j),(V_k,E_k,\Ga_k,s_k,\psi_k)$ be Kuranishi neighbourhoods on $X,Y,Z$ respectively, and $T\subseteq \Im\psi_j\cap g^{-1}(\Im\psi_k)\ab\subseteq Y$ and $S\subseteq\Im\psi_i\cap f^{-1}(T)\subseteq X$ be open. Suppose $\Phi_{ij}=(P_{ij},\pi_{ij},\phi_{ij},\hat\phi_{ij})\colon\ab (V_i,E_i,\Ga_i,s_i,\psi_i)\ra (V_j,E_j,\Ga_j,s_j,\psi_j)$ is a 1-morphism of Kuranishi neighbourhoods over $(S,f)$, and $\Phi_{jk}=(P_{jk},\pi_{jk},\phi_{jk},\hat\phi_{jk})\colon (V_j,E_j,\Ga_j,s_j,\psi_j)\ra (V_k,E_k,\ab\Ga_k,\ab s_k,\ab\psi_k)$ is a 1-morphism of Kuranishi neighbourhoods over~$(T,g)$.

Consider the diagram of manifolds and smooth maps:
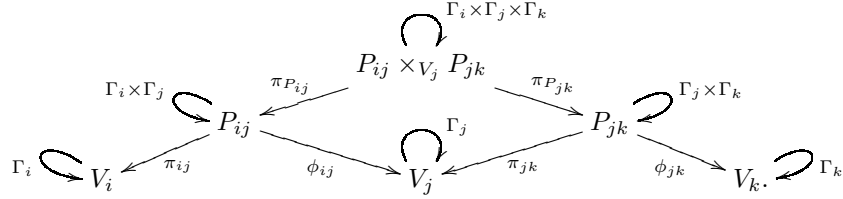
\begin{equation*}
\xymatrix@C=33pt@R=6pt{ 
&& P_{ij}\t_{V_j}P_{jk} \ar@(ul,ur)[]^(0.7){\Ga_i\t\Ga_j\t\Ga_k} \ar[dl]_(0.65){\pi_{P_{ij}}} \ar[dr]^(0.65){\pi_{P_{jk}}} \\
& P_{ij} \ar@(ul,l)[]_(0.5){\Ga_i\t\Ga_j} \ar[dl]^{\pi_{ij}} \ar[dr]_{\phi_{ij}} && P_{jk} \ar@(ur,r)[]^(0.5){\Ga_j\t\Ga_k} \ar[dl]^{\pi_{jk}} \ar[dr]_{\phi_{jk}} \\
V_i \ar@(ul,l)[]_(0.6){\Ga_i} && V_j \ar@(ul,ur)[]^(0.7){\Ga_j} && V_k. \ar@(ur,r)[]^(0.6){\Ga_k} }
\end{equation*}
Here the fibre product $P_{ij}\t_{V_j}P_{jk}$ is transverse, and so exists, as $\pi_{jk}$ is \'etale. We have shown the actions of various combinations of $\Ga_i,\Ga_j,\Ga_k$ on each space. In fact $\Ga_i\t\Ga_j\t\Ga_k$ acts on the whole diagram, with all maps equivariant, but we have omitted the trivial actions (for instance, $\Ga_j,\Ga_k$ act trivially on~$V_i$).

As $\Ga_j$ acts freely on $P_{ij}$, it also acts freely on $P_{ij}\t_{V_j}P_{jk}$, so $P_{ik}:=(P_{ij}\t_{V_j}P_{jk})/\Ga_j$ is a manifold, with projection $\Pi\colon P_{ij}\t_{V_j}P_{jk}\ra P_{ik}$. The commuting actions of $\Ga_i,\Ga_k$ on $P_{ij}\t_{V_j}P_{jk}$ descend to commuting actions of $\Ga_i,\Ga_k$ on $P_{ik}$ such that $\Pi$ is $\Ga_i$- and $\Ga_k$-equivariant. As $\pi_{ij}\ci\pi_{P_{ij}}\colon P_{ij}\t_{V_j}P_{jk}\ra V_i$ and $\phi_{jk}\ci\pi_{P_{jk}}\colon P_{ij}\t_{V_j}P_{jk}\ra V_k$ are $\Ga_j$-invariant, they factor through $\Pi$, so there are unique smooth maps $\pi_{ik}\colon P_{ik}\ra V_i$ and $\phi_{ik}\colon P_{ik}\ra V_k$ such that $\pi_{ij}\ci\pi_{P_{ij}}=\pi_{ik}\ci\Pi$ and~$\phi_{jk}\ci\pi_{P_{jk}}=\phi_{ik}\ci\Pi$. 

Consider the diagram of vector bundles on $P_{ij}\t_{V_j}P_{jk}$:
\begin{equation*}
\xymatrix@C=15pt@R=13pt{
*+[r]{\Pi^*\ci\pi_{ik}^*(E_i)} \ar@{.>}[rrrrrrr]_{\Pi^*(\hat\phi_{ik})} \ar@{=}[d] &&&&&&& *+[l]{\Pi^*\ci \phi_{ik}^*(E_k)} \ar@{=}[d] \\
*+[r]{\pi_{P_{ij}}^*\ci\pi_{ij}^*(E_i)} \ar[rrr]^(0.6){\pi_{P_{ij}}^*(\hat\phi_{ij})} 
&&& \pi_{P_{ij}}^*\ci\phi_{ij}^*(E_j) \ar@{=}[r] &
\pi_{P_{jk}}^*\ci\pi_{jk}^*(E_j) \ar[rrr]^(0.4){\pi_{P_{jk}}^*(\hat\phi_{jk})} &&&
*+[l]{\pi_{P_{jk}}^*\ci\phi_{jk}^*(E_k).} }
\end{equation*}
There is a unique morphism on the top line making the diagram commute. As $\hat\phi_{ij},\hat\phi_{jk}$ are $\Ga_j$-equivariant, this is $\Ga_j$-equivariant, so it is the pullback under $\Pi^*$ of a unique morphism $\hat\phi_{ik}\colon \pi_{ik}^*(E_i)\ra \phi_{ik}^*(E_k)$, as shown. It is now easy to check that $(P_{ik},\pi_{ik},\phi_{ik},\hat\phi_{ik})$ satisfies Definition \ref{kt3def3}(a)--(f), and is a 1-morphism $\Phi_{ik}=(P_{ik},\pi_{ik},\phi_{ik},\hat\phi_{ik})\colon (V_i,E_i,\Ga_i,s_i,\psi_i)\ra(V_k,E_k,\Ga_k,s_k,\psi_k)$ over $(S,g\ci f)$. We write $\Phi_{jk}\ci\Phi_{ij}=\Phi_{ik}$, and call it the {\it composition of\/ $1$-morphisms}. 

If we have three such 1-morphisms $\Phi_{ij},\Phi_{jk},\Phi_{kl}$, define
\e
\la_{ijkl}\colon \bigl[P_{ij}\t_{V_j}\bigl((P_{jk}\t_{V_k}P_{kl})/\Ga_k\bigr)\bigr]/\Ga_j\ra\bigl[\bigl((P_{ij}\t_{V_j}P_{jk})/\Ga_j\bigr)\t_{V_k}P_{kl}\bigr]/\Ga_k
\label{kt3eq3}
\e
to be the natural identification. Then we have a 2-isomorphism
\e
\begin{split}
\bs\al_{\Phi_{kl},\Phi_{jk},\Phi_{ij}}:=\bigl[[P_{ij}\t_{V_j}((P_{jk}&\t_{V_k}P_{kl})/\Ga_k)]/\Ga_j,\la_{ijkl},0\bigr]\colon \\
&(\Phi_{kl}\ci\Phi_{jk})\ci\Phi_{ij}\Longra\Phi_{kl}\ci(\Phi_{jk}\ci\Phi_{ij}).
\end{split}
\label{kt3eq4}
\e
That is, composition of $1$-morphisms is associative up to canonical 2-iso\-mor\-ph\-ism, as for weak 2-categories in~\S\ref{ktA1}.

For $\Phi_{ij}\colon (V_i,E_i,\Ga_i,s_i,\psi_i)\ra (V_j,E_j,\Ga_j,s_j,\psi_j)$ as above, define 
\begin{align*}
\mu_{ij}&\colon ((V_i\t\Ga_i)\t_{V_i}P_{ij})/\Ga_i\longra P_{ij},\\
\nu_{ij}&\colon (P_{ij}\t_{V_j}(V_j\t\Ga_j))/\Ga_j\longra P_{ij},
\end{align*}
to be the natural identifications. Then we have 2-isomorphisms
\e
\begin{split}
\bs\be_{\Phi_{ij}}:=\bigl[((V_i\t\Ga_i)\t_{V_i}P_{ij})/\Ga_i,\mu_{ij},0\bigr]&\colon \Phi_{ij}\ci\id_{(V_i,E_i,\Ga_i,s_i,\psi_i)}\Longra\Phi_{ij},\\
\bs\ga_{\Phi_{ij}}:=\bigl[(P_{ij}\t_{V_j}(V_j\t\Ga_j))/\Ga_j,\nu_{ij},0\bigr]&\colon \id_{(V_j,E_j,\Ga_j,s_j,\psi_j)}\ci\Phi_{ij}\Longra\Phi_{ij},
\end{split}
\label{kt3eq5}
\e
so identity 1-morphisms behave as they should up to canonical 2-isomorphism, as for weak 2-categories in~\S\ref{ktA1}.
\label{kt3def5}
\end{dfn}

\begin{dfn} Let $X,Y$ be topological spaces, $f\colon X\ra Y$ be continuous, $(V_i,\ab E_i,\ab\Ga_i,\ab s_i,\ab\psi_i)$, $(V_j,E_j,\Ga_j,s_j,\psi_j)$ be Kuranishi neighbourhoods on $X,Y$, $S\subseteq \Im\psi_i\cap f^{-1}(\Im\psi_j)\subseteq X$ be open, and $\Phi_{ij},\Phi_{ij}',\Phi_{ij}''\colon (V_i,\ab E_i,\ab\Ga_i,\ab s_i,\ab\psi_i)\ra (V_j,\ab E_j,\ab\Ga_j,\ab s_j,\ab\psi_j)$ be 1-morphisms over $(S,f)$ with $\Phi_{ij}=(P_{ij},\pi_{ij},\phi_{ij},\hat\phi_{ij})$, $\Phi_{ij}'=(P_{ij}',\pi_{ij}',\phi_{ij}',\hat\phi_{ij}')$, $\Phi_{ij}''=(P_{ij}'',\pi_{ij}'',\phi_{ij}'',\hat\phi_{ij}'')$. Suppose $\La_{ij}=[\dot P_{ij},\la_{ij},\hat\la_{ij}]\colon \Phi_{ij}\Ra\Phi_{ij}'$ and $\La_{ij}'=[\dot P_{ij}',\la_{ij}',\hat\la_{ij}']\colon \Phi_{ij}'\Ra\Phi_{ij}''$ are 2-morphisms over $(S,f)$. We will define the {\it vertical composition of\/ $2$-morphisms over\/} $(S,f)$, written
\begin{equation*}
\La_{ij}'\od\La_{ij}=[\dot P_{ij}',\la_{ij}',\hat\la_{ij}']\od [\dot P_{ij},\la_{ij},\hat\la_{ij}]\colon \Phi_{ij}\Longra\Phi_{ij}''.
\end{equation*}

Choose representatives $(\dot P_{ij},\la_{ij},\hat\la_{ij}),(\dot P_{ij}',\la_{ij}',\hat\la_{ij}')$ in the $\approx$-equivalence cla\-sses $[\dot P_{ij},\la_{ij},\hat\la_{ij}],[\dot P_{ij}',\la_{ij}',\hat\la_{ij}']$. Define $\dot P_{ij}''=\la_{ij}^{-1}(\dot P_{ij})\subseteq \dot P_{ij}\subseteq P_{ij}$, and $\la_{ij}''=\la_{ij}'\ci\la_{ij}\vert_{\dot P_{ij}''}$. Consider the morphism of vector bundles
\begin{equation*}
\la_{ij}^*(\hat\la_{ij}')\colon \pi_{ij}^*(E_i)\vert_{\dot P_{ij}''}=\la_{ij}^*\ci\pi_{ij}^{\prime *}(E_i)\vert_{\dot P_{ij}''}\longra\la_{ij}^*\ci\phi_{ij}^{\prime *}(TV_j)=(\phi_{ij}'\ci\la_{ij})^*(TV_j)\vert_{\dot P_{ij}''}.
\end{equation*}
Since $\phi_{ij}'\ci\la_{ij}\vert_{\dot P_{ij}''}=\phi_{ij}\vert_{\dot P_{ij}''}+O(\pi_{ij}^*(s_i))$ by \eq{kt3eq2}, the discussion after Definition \ref{kt3def1}(vi) shows that there exists $\check\la'_{ij}\colon \pi_{ij}^*(E_i)\vert_{\dot P_{ij}''}\ra \phi_{ij}^*(TV_j)\vert_{\dot P_{ij}''}$ with
\e
\check\la_{ij}'=\la_{ij}\vert_{\dot P_{ij}''}^*(\hat\la_{ij}')+O(\pi_{ij}^*(s_i)),
\label{kt3eq6}
\e
as in Definition \ref{kt3def1}(vi), and $\check\la_{ij}'$ is unique up to $O(\pi_{ij}^*(s_i))$. By averaging over the $\Ga_i\t\Ga_j$-action we can suppose $\check\la_{ij}'$ is $\Ga_i$- and $\Ga_j$-equivariant, as $\hat\la_{ij}'$ is.

Define $\hat\la_{ij}''\colon \pi_{ij}^*(E_i)\vert_{\dot P_{ij}''}\ra \phi_{ij}^*(TV_j)\vert_{\dot P_{ij}''}$ by $\hat\la_{ij}''=\hat\la_{ij}\vert_{\dot P_{ij}''}+\check\la_{ij}'$. It is now easy to see that $(\dot P_{ij}'',\la_{ij}'',\hat\la_{ij}'')$ satisfies Definition \ref{kt3def4}(a)--(c) for $\Phi_{ij},\Phi_{ij}''$, using \eq{kt3eq2} for $\hat\la_{ij},\hat\la_{ij}'$ and \eq{kt3eq6} to prove \eq{kt3eq2} for $\hat\la_{ij}''$. Hence $\La_{ij}''=[\dot P_{ij}'',\la_{ij}'',\hat\la_{ij}'']\colon \Phi_{ij}\Ra\Phi_{ij}''$ is a 2-morphism over $(S,f)$. It is independent of choices. We define $[\dot P_{ij}',\la_{ij}',\hat\la_{ij}']\od [\dot P_{ij},\la_{ij},\hat\la_{ij}]=[\dot P_{ij}'',\la_{ij}'',\hat\la_{ij}'']$, or $\La_{ij}'\od\La_{ij}=\La_{ij}''$. All $2$-morphisms over $(S,f)$ are invertible under vertical composition, that is, they are 2-isomorphisms.

Write $\bHom_{S,f}\bigl((V_i,E_i,\Ga_i,s_i,\psi_i),(V_j,E_j,\Ga_j,s_j,\psi_j)\bigr)$ for the groupoid with objects 1-morphisms $\Phi_{ij}\colon (V_i,E_i,\Ga_i,s_i,\psi_i)\ab\ra (V_j,\ab E_j,\ab\Ga_j,\ab s_j,\ab\psi_j)$ over $(S,f)$, and morphisms 2-morphisms $\La_{ij}\colon \Phi_{ij}\Ra\Phi_{ij}'$ over $(S,f)$. If $X=Y$ and $f=\id_X$, we write this as~$\bHom_S\bigl((V_i,E_i,\Ga_i,s_i,\psi_i),(V_j,E_j,\Ga_j,s_j,\psi_j)\bigr)$.
\label{kt3def6}
\end{dfn}

\begin{dfn} Let $X,Y,Z$ be topological spaces, $f\colon X\ra Y$, $g\colon Y\ra Z$ be continuous maps, $(V_i,E_i,\Ga_i,s_i,\psi_i),(V_j,E_j,\Ga_j,s_j,\psi_j),(V_k,E_k,\Ga_k,s_k,\psi_k)$ be Kuranishi neighbourhoods on $X,Y,Z$, and $T\subseteq \Im\psi_j\cap g^{-1}(\Im\psi_k)\ab\subseteq Y$ and $S\subseteq\Im\psi_i\cap f^{-1}(T)\subseteq X$ be open. Suppose $\Phi_{ij},\Phi_{ij}'\colon (V_i,E_i,\Ga_i,s_i,\psi_i)\ra (V_j,E_j,\Ga_j,s_j,\psi_j)$ are 1-morphisms of Kuranishi neighbourhoods over $(S,f)$, and $\La_{ij}\colon \Phi_{ij}\Ra\Phi_{ij}'$ is a 2-morphism over $(S,f)$, and $\Phi_{jk},\Phi_{jk}'\colon (V_j,E_j,\Ga_j,\ab s_j,\ab\psi_j)\ab\ra (V_k,E_k,\ab\Ga_k,\ab s_k,\ab\psi_k)$ are 1-morphisms of Kuranishi neighbourhoods over $(T,g)$, and $\La_{jk}\colon \Phi_{jk}\Ra\Phi_{jk}'$ is a 2-morphism over~$(T,g)$.

We will define the {\it horizontal composition of\/ $2$-morphisms}, written
\e
\La_{jk}*\La_{ij}\colon \Phi_{jk}\ci\Phi_{ij}\Longra\Phi_{jk}'\ci\Phi_{ij}'\qquad\text{over $(S,g\ci f)$.}
\label{kt3eq7}
\e
Use our usual notation for $\Phi_{ij},\ldots,\La_{jk}$, and write $(P_{ik},\pi_{ik},\phi_{ik},\hat\phi_{ik})=\Phi_{jk}\ci\Phi_{ij}$, $(P_{ik}',\pi_{ik}',\phi_{ik}',\hat\phi_{ik}')=\Phi_{jk}'\ci\Phi_{ij}'$, as in Definition \ref{kt3def5}. Choose representatives $(\dot P_{ij},\la_{ij},\hat\la_{ij})$, $(\dot P_{jk},\la_{jk},\hat\la_{jk})$ for $\La_{ij}=[\dot P_{ij},\la_{ij},\hat\la_{ij}]$ and~$\La_{jk}=[\dot P_{jk},\la_{jk},\hat\la_{jk}]$. 

Then $P_{ik}=(P_{ij}\t_{V_j}P_{jk})/\Ga_j$, and $\dot P_{ij}\subseteq P_{ij}$, $\dot P_{jk}\subseteq P_{jk}$ are open and $\Ga_j$-invariant, so $\dot P_{ij}\t_{V_j}\dot P_{jk}$ is open and $\Ga_j$-invariant in $P_{ij}\t_{V_j}P_{jk}$. Define $\dot P_{ik}=(\dot P_{ij}\t_{V_j}\dot P_{jk})/\Ga_j$, as an open subset of $P_{ik}$. It is $\Ga_i$- and $\Ga_k$-invariant, as $\dot P_{ij}$, $\dot P_{jk}$ are $\Ga_i$- and $\Ga_k$-invariant, respectively.

The maps $\la_{ij}\colon \dot P_{ij}\ra P_{ij}'$, $\la_{jk}\colon \dot P_{jk}\ra P_{jk}'$ satisfy $\phi_{ij}'\ci\la_{ij}=\phi_{ij}\vert_{\dot P_{ij}}\colon \dot P_{ij}\ra V_j$ and $\pi_{jk}'\ci\la_{jk}=\pi_{jk}\vert_{\dot P_{jk}}\colon \dot P_{jk}\ra V_j$. Hence by properties of fibre products they induce a unique smooth map $\ti\la_{ik}\colon \dot P_{ij}\t_{\phi_{ij},V_j,\pi_{jk}}\dot P_{jk}\ra P_{ij}'\t_{\phi_{ij}',V_j,\pi_{jk}'}P_{jk}'$ with $\pi_{P_{ij}'}\ci\ti\la_{ik}=\la_{ij}\ci\pi_{\dot P_{ij}}$ and $\pi_{P_{jk}'}\ci\ti\la_{ik}=\la_{jk}\ci\pi_{\dot P_{jk}}$. As everything is $\Ga_j$-equivariant, $\ti\la_{ik}$ descends to the quotients by $\Ga_j$. Thus we obtain a unique smooth map
\begin{equation*}
\la_{ik}\colon \dot P_{ik}=(\dot P_{ij}\t_{V_j}\dot P_{jk})/\Ga_j\longra 
(P_{ij}'\t_{V_j}P_{jk}')/\Ga_j=P_{ik}'
\end{equation*}
with $\la_{ik}\ci\Pi=\Pi'\ci\ti\la_{ik}$, for $\Pi\colon \dot P_{ij}\t_{V_j}\dot P_{jk}\ra (\dot P_{ij}\t_{V_j}\dot P_{jk})/\Ga_j$, $\Pi'\colon P_{ij}'\t_{V_j}P_{jk}'\ra (P_{ij}'\t_{V_j}P_{jk}')/\Ga_j$ the projections.

Define a morphism of vector bundles on $\dot P_{ij}\t_{V_j}\dot P_{jk}$
\begin{gather*}
\check\la_{ik}\colon \Pi^*\ci\pi_{ik}^*(E_i)=(\pi_{ij}\ci\pi_{\dot P_{ij}})^*(E_i)\longra (\phi_{jk}\ci\pi_{\dot P_{jk}})^*(TV_k)=\Pi^*\ci\phi_{ik}^*(TV_k)\\
\text{by}\;\>
\check\la_{ik}=\pi_{\dot P_{jk}}^*(\d\phi_{jk}\ci (\d\pi_{jk})^{-1})\ci\pi_{\dot P_{ij}}^*(\hat\la_{ij})+
\pi_{\dot P_{jk}}^*(\hat\la_{jk})\ci \pi_{\dot P_{ij}}^*(\hat\phi_{ij}),
\end{gather*}
where the morphisms are given in the diagram
\begin{equation*}
\xymatrix@C=83pt@R=2pt{
*+[r]{(\pi_{ij}\ci\pi_{\dot P_{ij}})^*(E_i)} 
\ar[dd]^{\pi_{\dot P_{ij}}^*(\hat\la_{ij})}
\ar[r]_(0.6){\pi_{\dot P_{ij}}^*(\hat\phi_{ij})} & (\phi_{ij}\ci\pi_{\dot P_{ij}})^*(E_j) \ar@{=}[r] & *+[l]{(\pi_{jk}\ci\pi_{\dot P_{jk}})^*(E_j)} \ar[ddd]_{\pi_{\dot P_{jk}}^*(\hat\la_{jk})} 
\\ \\
*+[r]{(\phi_{ij}\ci\pi_{\dot P_{ij}})^*(TV_j)} \ar@{=}[d]
\\
*+[r]{(\pi_{jk}\ci\pi_{\dot P_{jk}})^*(TV_j)} 
\ar@<1ex>[r]^(0.65){\pi_{\dot P_{jk}}^*((\d\pi_{jk})^{-1})} & \pi_{\dot P_{jk}}^*(T\dot P_{jk}) \ar[r]^(0.35){\pi_{\dot P_{jk}}^*(\d\phi_{jk})} 
\ar[l]^(0.35){\pi_{\dot P_{jk}}^*(\d\pi_{jk})} & 
*+[l]{(\phi_{jk}\ci\pi_{\dot P_{jk}})^*(TV_k).} }
\end{equation*}
Here $\d\pi_{jk}\colon T\dot P_{jk}\ra \pi_{jk}^*(TV_j)$ is invertible as $\pi_{jk}$ is \'etale. As all the ingredients are $\Ga_i,\Ga_j,\Ga_k$-invariant or equivariant, $\check\la_{ik}$ is $\Ga_j$-invariant, and so descends to $\dot P_{ik}=(\dot P_{ij}\t_{V_j}\dot P_{jk})/\Ga_j$. That is, there is a unique morphism $\hat\la_{ik}\colon \pi_{ik}\vert_{\dot P_{ik}}^*(E_i)\ra \phi_{ik}\vert_{\dot P_{ik}}^*(TV_k)$ of vector bundles on $\dot P_{ik}$ with $\Pi^*(\hat\la_{ik})=\check\la_{ik}$. As $\check\la_{ik}$ is $\Ga_i$- and $\Ga_k$-equivariant, so is $\hat\la_{ik}$.

One can now check that $(\dot P_{ik},\la_{ik},\hat\la_{ik})$ satisfies Definition \ref{kt3def4}(a)--(c), where \eq{kt3eq2} for $\hat\la_{ik}$ follows from adding the pullbacks to $\dot P_{ij}\t_{V_j}\dot P_{jk}$ of \eq{kt3eq2} for $\hat\la_{ij},\hat\la_{jk}$, so $\La_{ik}=[\dot P_{ik},\la_{ik},\hat\la_{ik}]$ is a 2-morphism as in \eq{kt3eq7}, which is independent of choices of $(\dot P_{ij},\la_{ij},\hat\la_{ij})$, $(\dot P_{jk},\la_{jk},\hat\la_{jk})$. We define $\La_{jk}*\La_{ij}=\La_{ik}$ in~\eq{kt3eq7}.
\label{kt3def7}
\end{dfn}

We have now defined all the structures of a weak 2-category, except that our 1- and 2-morphisms are defined over an open set $S\subseteq X$, which is not part of the 2-category structure. As in \cite[\S 4.1]{Joyc5} and \cite[\S 6.1]{Joyc6}, there are two ways to make a genuine 2-category out of this, either (a) working on a fixed topological space $X$ and open $S\subseteq X$, or (b) allowing $X$ to vary but requiring $S=X$ throughout.

\begin{thm}{\bf(a)} Let\/ $X$ be a topological space and\/ $S\subseteq X$ an open set. The definitions above give a \begin{bfseries}weak\/ $2$-category $\KN_S(X)$ of Kuranishi neighbourhoods over\/\end{bfseries} $S,$ with objects Kuranishi neighbourhoods $(V_i,E_i,\Ga_i,s_i,\psi_i)$ on $X$ with $S\subseteq\Im\psi_i$, and $1$- and $2$-morphisms $1$-morphisms $\Phi_{ij}$ of Kuranishi neighbourhoods over $S,$ and $2$-morphisms $\La_{ij}$ over $S$.
\smallskip

\noindent{\bf(b)} The definitions above also give a weak\/ $2$-category\/ $\GKN$ of \begin{bfseries}global Kuranishi neighbourhoods\end{bfseries}, such that:
\begin{itemize}
\setlength{\itemsep}{0pt}
\setlength{\parsep}{0pt}
\item[{\bf(i)}] The objects $\bigl(X,(V,E,\Ga,s,\psi)\bigr)$ in\/ $\GKN$ are pairs of a topological space $X$ and a Kuranishi neighbourhood\/ $(V,E,\Ga,s,\psi)$ on $X$ with\/~$\Im\psi=X$.
\item[{\bf(ii)}] The\/ $1$-morphisms $(f,\Phi)\!:\! \bigl(X,(V,E,\Ga,s,\psi)\bigr)\!\ra\!\bigl(Y,(W,F,\De,t,\chi)\bigr)$ in\/ $\GKN$ are pairs of a continuous map $f\colon X\ra Y$ and a $1$-morphism $\Phi\colon (V,\ab E,\ab\Ga,\ab s,\ab\psi)\ab\ra(W,F,\De,t,\chi)$ of Kuranishi neighbourhoods over $(X,f)$.
\item[{\bf(iii)}] For $1$-morphisms $(f,\Phi),(g,\Psi)\colon \bigl(X,(V,E,\Ga,s,\psi)\bigr)\ra\bigl(Y,(W,F,\De,t,\chi)\bigr),$ the $2$-morphisms $\La\colon (f,\Phi)\Ra(g,\Psi)$ in $\GKN$ exist only if\/ $f=g,$ and are $2$-morphisms $\La\colon \Phi\Ra\Psi$ of Kuranishi neighbourhoods over\/~$(X,f)$.
\end{itemize}
All\/ $2$-morphisms in $\KN_S(X),\GKN$ are $2$-isomorphisms, that is, $\KN_S(X),\ab\GKN$ are $(2,1)$-categories.
\label{kt3thm1}
\end{thm}

\begin{dfn} Recall from \S\ref{ktA1} that an {\it equivalence\/} in a 2-category $\bs\cC$ is a 1-morphism $f\colon A\ra B$ in $\bs\cC$ such that there exist a 1-morphism $g\colon B\ra A$ (called a {\it quasi-inverse\/}) and 2-isomorphisms $\eta\colon g\ci f\Ra\id_A$ and $\ze\colon f\ci g\Ra\id_B$. A 1-morphism $f\colon A\ra B$ is an equivalence if and only if $[f]\colon A\ra B$ is an isomorphism (is invertible) in the homotopy category~$\Ho(\bs\cC)$.

A 1-morphism $\Phi_{ij}\colon (V_i,E_i,\Ga_i,s_i,\psi_i)\ra (V_j,E_j,\Ga_j,s_j,\psi_j)$ on $X$ over $S$ is a {\it coordinate change over\/} $S$ if $\Phi_{ij}$ is an equivalence in the 2-category~$\KN_S(X)$.

\label{kt3def8}
\end{dfn}

\begin{dfn} Let $T\subseteq S\subseteq X$ be open. Define the {\it restriction\/ $2$-functor\/} $\vert_T\colon \KN_S(X)\ra\KN_T(X)$ to map objects $(V_i,E_i,\Ga_i,s_i,\psi_i)$ to exactly the same objects, and 1-morphisms $\Phi_{ij}$ to exactly the same 1-morphisms but regarded as 1-morphisms over $T$, and 2-morphisms $\La_{ij}=[\dot P_{ij},\la_{ij},\hat\la_{ij}]$ over $S$ to $\La_{ij}\vert_T=[\dot P_{ij},\la_{ij},\hat\la_{ij}]\vert_T$, where $[\dot P_{ij},\la_{ij},\hat\la_{ij}]\vert_T$ is the $\approx_T$-equivalence class of any representative $(\dot P_{ij},\ab\la_{ij},\ab\hat\la_{ij})$ for the $\approx_S$-equivalence class~$[\dot P_{ij},\la_{ij},\hat\la_{ij}]$. 

Then $\vert_T\colon \KN_S(X)\ra\KN_T(X)$ commutes with all the structure, so we may we regard $\vert_T$ as a weak 2-functor as in \S\ref{ktA2} for which the additional 2-isomorphisms $F_{g,f}$ are identities.

If $U\subseteq T\subseteq S\subseteq X$ are open then $\vert_U\ci\vert_T=\vert_U\colon \KN_S(X)\ra\KN_U(X)$. Also $\vert_S$ gives a functor $\vert_T\colon \bHom_S\bigl((V_i,E_i,\Ga_i,s_i,\psi_i),(V_j,E_j,\Ga_j,s_j,\psi_j)\bigr)\ra\ab \bHom_T\bigl((V_i,\ab E_i,\ab\Ga_i,s_i,\psi_i),(V_j,E_j,\Ga_j,s_j,\psi_j)\bigr)$ when $T\subseteq S\subseteq \Im\psi_i\cap\Im\psi_j$, in the notation of Definition~\ref{kt3def6}.
\label{kt3def9}
\end{dfn}

\begin{dfn} So far we have discussed 1- and 2-morphisms of Kuranishi neighbourhoods, and coordinate changes, {\it over a specified open set\/} $S\subseteq X$, or over $(S,f)$. We now make the convention that {\it when we do not specify a domain\/ $S$ for a $1$-morphism, $2$-morphism, or coordinate change, the domain should be as large as possible}. For example, if we say that $\Phi_{ij}\colon (V_i,E_i,\Ga_i,s_i,\psi_i)\ra (V_j,E_j,\Ga_j,s_j,\psi_j)$ is a 1-morphism (or a 1-morphism over $f\colon X\ra Y$) without specifying $S$, we mean that $S=\Im\psi_i\cap\Im\psi_j$ (or~$S=\Im\psi_i\cap f^{-1}(\Im\psi_j)$).

Similarly, if we write a formula involving several 2-morphisms (possibly defined on different domains), without specifying the domain $S$, we make the convention {\it that the domain where the formula holds should be as large as possible}. That is, the domain $S$ is taken to be the intersection of the domains of each 2-morphism in the formula, and we implicitly restrict each morphism in the formula to $S$ as in Definition \ref{kt3def9}, to make it make sense.
\label{kt3def10}
\end{dfn}

\subsection{Properties of 1- and 2-morphisms}
\label{kt33}

The next result \cite[Th.~4.14]{Joyc5}, \cite[\S 6.1]{Joyc6} is very important in our theory. We will call Theorem \ref{kt3thm2} the {\it stack property}. We will use it in \S\ref{kt41} to construct compositions of 1- and 2-morphisms of Kuranishi spaces. {\it Stacks on topological spaces\/} are the 2-category analogue of sheaves, and are defined in~\S\ref{ktA4}.

\begin{thm}{\bf(a)} Let\/ $(V_i,E_i,\Ga_i,s_i,\psi_i),(V_j,E_j,\Ga_j,s_j,\psi_j)$ be Kuranishi nei\-ghbourhoods on a topological space $X$. For each open $S\subseteq\Im\psi_i\cap\Im\psi_j,$ set
\begin{align*}
\bcHom\bigl((V_i&,E_i,\Ga_i,s_i,\psi_i),(V_j,E_j,\Ga_j,s_j,\psi_j)\bigr)(S)\\
&=\bHom_S\bigl((V_i,E_i,\Ga_i,s_i,\psi_i),(V_j,E_j,\Ga_j,s_j,\psi_j)\bigr),
\end{align*}
as in Definition\/ {\rm\ref{kt3def6},} for all open $T\subseteq S\subseteq \Im\psi_i\cap\Im\psi_j$ define a functor
\begin{align*}
\rho_{ST}\colon \,&\bcHom\bigl((V_i,E_i,\Ga_i,s_i,\psi_i),(V_j,E_j,\Ga_j,s_j,\psi_j)\bigr)(S)\longra\\
&\bcHom\bigl((V_i,E_i,\Ga_i,s_i,\psi_i),(V_j,E_j,\Ga_j,s_j,\psi_j)\bigr)(T)\end{align*}
by $\rho_{ST}=\vert_T,$ as in Definition\/ {\rm\ref{kt3def9},} and for all open $U\subseteq T\subseteq S\subseteq \Im\psi_i\cap\Im\psi_j$ take the obvious isomorphism $\eta_{STU}=\id_{\rho_{SU}}\colon \rho_{TU}\ci\rho_{ST}\Ra \rho_{SU}$.
Then $\bcHom\bigl((V_i,E_i,\Ga_i,s_i,\psi_i),(V_j,E_j,\Ga_j,s_j,\psi_j)\bigr)$ is a stack on $\Im\psi_i\cap\Im\psi_j$. 

Coordinate changes $(V_i,E_i,\Ga_i,s_i,\psi_i)\ra(V_j,E_j,\Ga_j,s_j,\psi_j)$ also form a stack\/ $\bcEqu\bigl((V_i,E_i,\Ga_i,s_i,\psi_i),(V_j,E_j,\Ga_j,s_j,\psi_j)\bigr)$ on $\Im\psi_i\cap\Im\psi_j,$ which is a substack of\/~$\bcHom\bigl((V_i,E_i,\Ga_i,s_i,\psi_i),(V_j,E_j,\Ga_j,s_j,\psi_j)\bigr)$.
\smallskip

\noindent{\bf(b)} Let\/ $f\colon X\ra Y$ be continuous, and\/ $(U_i,D_i,\Be_i,r_i,\chi_i),(V_j,E_j,\ab\Ga_j,\ab s_j,\ab\psi_j)$ be Kuranishi neighbourhoods on $X,Y,$ respectively. Then $1$- and\/ $2$-morphisms from $(U_i,D_i,\Be_i,r_i,\chi_i)$ to $(V_j,E_j,\Ga_j,s_j,\psi_j)$ over $f$ form a stack\/ $\bcHom_f\bigl((U_i,\ab D_i,\ab\Be_i,\ab r_i,\ab\chi_i),\ab (V_j,E_j,\Ga_j,s_j,\psi_j)\bigr)$ on\/~$\Im\chi_i\cap f^{-1}(\Im\psi_j)\subseteq X$.
\label{kt3thm2}
\end{thm}

Here \cite[Th.~4.16]{Joyc5}, \cite[\S 10.5]{Joyc6} is a criterion for recognizing when a 1-morphism of Kuranishi neighbourhoods is a coordinate change, that is, is invertible up to 2-isomorphism. Note the similarity of \eq{kt3eq8} to equation~\eq{kt2eq2}.

\begin{thm} Let\/ $\Phi_{ij}=(P_{ij},\pi_{ij},\phi_{ij},\hat\phi_{ij})\colon (V_i,E_i,\Ga_i,s_i,\psi_i)\ra (V_j,E_j,\ab\Ga_j,\ab s_j,\ab\psi_j)$ be a $1$-morphism of Kuranishi neighbourhoods over\/ $S\subseteq X$. Let\/ $p\in\pi_{ij}^{-1}(\bar\psi_i^{-1}(S))\subseteq P_{ij},$ and set\/ $v_i=\pi_{ij}(p)\in V_i$ and\/ $v_j=\phi_{ij}(p)\in V_j$. Consider the sequence of real vector spaces:
\e
\xymatrix@C=14pt{ 0 \ar[r] & T_{v_i}V_i \ar[rrrrr]^(0.43){\d s_i\vert_{v_i}\op(\d\phi_{ij}\vert_p\ci\d\pi_{ij}\vert_p^{-1})} &&&&& E_i\vert_{v_i} \op T_{v_j}V_j 
\ar[rrr]^(0.56){-\hat\phi_{ij}\vert_p\op \d s_j\vert_{v_j}} &&& E_j\vert_{v_j} \ar[r] & 0. }{}
\label{kt3eq8}
\e
Here $\d\pi_{ij}\vert_p\colon T_pP_{ij}\ra T_{v_i}V_i$ is invertible as $\pi_{ij}$ is \'etale. Definition\/ {\rm\ref{kt3def3}(e)} implies that\/ \eq{kt3eq8} is a complex. Also consider the morphism of finite groups
\e
\begin{split}
&\rho_p\colon \bigl\{(\ga_i,\ga_j)\in\Ga_i\t\Ga_j\colon (\ga_i,\ga_j)\cdot p=p\bigr\}
\longra\bigl\{\ga_j\in\Ga_j\colon \ga_j\cdot v_j=v_j\bigr\},\\
&\rho_p\colon (\ga_i,\ga_j)\longmapsto \ga_j.
\end{split}
\label{kt3eq9}
\e

Then $\Phi_{ij}$ is a coordinate change over $S,$ in the sense of Definition\/ {\rm\ref{kt3def8},} if and only if\/ \eq{kt3eq8} is exact and\/ \eq{kt3eq9} is an isomorphism for all\/~$p\in\pi_{ij}^{-1}(\bar\psi_i^{-1}(S))$.

\label{kt3thm3}
\end{thm}

\subsection{\texorpdfstring{Relation to Fukaya--Oh--Ohta--Ono's work}{Relation to Fukaya-Oh-Ohta-Ono\textquoteright s work}}
\label{kt34}

We now relate our definitions in \S\ref{kt31}--\S\ref{kt33} to FOOO Kuranishi neighbourhoods and FOOO coordinate changes from \S\ref{kt21}.

\begin{ex} Compare Definition \ref{kt2def1} of FOOO Kuranishi neighbourhoods $(V,\ab E,\ab\Ga,\ab s,\ab\psi)$, and Definition \ref{kt3def2} of our Kuranishi neighbourhoods $(\ti V,\ti E,\ti\Ga,\ti s,\ab\ti\psi)$. The differences are that $E$ is a vector space, but $\ti E\ra\ti V$ is a vector bundle over $\ti V$. Also $\Ga$ must act effectively on $V$, but $\ti\Ga$ need not act effectively on $\ti V$.

To make a FOOO Kuranishi neighbourhood $(V,E,\Ga,s,\psi)$ into one of our Kuranishi neighbourhoods $(\ti V,\ti E,\ti\Ga,\ti s,\ti\psi)$, take $\ti V=V$, $\ti\Ga=\Ga,$ $\ti\psi=\psi$, let $\ti\pi\colon \ti E\ra \ti V$ be the trivial vector bundle $\pi_V\colon V\t E\ra V$ with fibre $E$, and $\ti s=\id\t s\colon V\ra V\t E$. Thus, FOOO Kuranishi neighbourhoods correspond to special examples of our Kuranishi neighbourhoods $(\ti V,\ti E,\ti\Ga,\ti s,\ti\psi)$, in which $\ti\pi\colon \ti E\ra \ti V$ is a trivial vector bundle, and $\ti\Ga$ acts effectively on~$\ti V$.

By an abuse of notation, we will sometimes identify FOOO Kuranishi neighbourhoods with the corresponding Kuranishi neighbourhoods in \S\ref{kt31}. That is, we will use $E$ to denote both a vector space, and the corresponding trivial vector bundle over $V$, and $s$ to denote both a map, and a section of a trivial bundle. Fukaya et al.\ \cite[Def.~4.3(4)]{FOOO4} also make the same abuse of notation.
\label{kt3ex1}
\end{ex}

\begin{ex} Let $\Phi_{ij}=(V_{ij},h_{ij},\vp_{ij},\hat\vp_{ij})\colon (V_i,E_i,\Ga_i,s_i,\psi_i)\ra(V_j,\ab E_j,\ab\Ga_j,\ab s_j,\ab\psi_j)$ be a FOOO coordinate change over $S$, as in Definition \ref{kt2def2}. As in Example \ref{kt3ex1}, regard the FOOO Kuranishi neighbourhoods $(V_i,E_i,\Ga_i,s_i,\psi_i),(V_j,\ab E_j,\ab\Ga_j,\ab s_j,\ab\psi_j)$ as examples of Kuranishi neighbourhoods in the sense of~\S\ref{kt31}.

Set $P_{ij}=V_{ij}\t\Ga_j$. Let $\Ga_i$ act on $P_{ij}$ by $\ga_i\colon (v,\ga)\mapsto (\ga_i\cdot v,\ga h_{ij}(\ga_i)^{-1})$. Let $\Ga_j$ act on $P_{ij}$ by $\ga_j\colon (v,\ga)\mapsto (v,\ga_j\ga)$. Define $\pi_{ij}\colon P_{ij}\ra V_i$ and $\phi_{ij}\colon P_{ij}\ra V_j$ by $\pi_{ij}\colon (v,\ga)\mapsto v$ and $\phi_{ij}\colon (v,\ga)\mapsto\ga\cdot\vp_{ij}(v)$. Then $\pi_{ij}$ is $\Ga_i$-equivariant and $\Ga_j$-invariant. Since $\vp_{ij}$ is $h_{ij}$-equivariant, $\phi_{ij}$ is $\Ga_i$-invariant, and $\Ga_j$-equivariant.

We will define a vector bundle morphism $\hat\phi_{ij}\colon \pi_{ij}^*(E_i)\ra\phi_{ij}^*(E_j)$. At $(v,\ga)\in P_{ij}$, this $\hat\phi_{ij}$ must map $E_i\vert_v\ra E_j\vert_{\ga\cdot \vp_{ij}(v)}$. We define $\hat\phi_{ij}\vert_{(v,\ga)}$ to be the composition of $\hat\vp_{ij}\vert_v\colon E_i\vert_v\ra E_j\vert_{\vp_{ij}(v)}$ with $\ga\cdot \colon E_j\vert_{\vp_{ij}(v)}\ra E_j\vert_{\ga\cdot\vp_{ij}(v)}$ from the $\Ga_j$-action on $E_j$. That is, $\hat\phi_{ij}\vert_{V_{ij}\t\{\ga\}}=\ga\cdot \hat\vp_{ij}$ for each $\ga\in\Ga_j$.

It is now easy to see that $\ti\Phi_{ij}=(P_{ij},\pi_{ij},\phi_{ij},\hat\phi_{ij})\colon (V_i,E_i,\Ga_i,s_i,\psi_i)\ra (V_j,\ab E_j,\ab\Ga_j,\ab s_j,\ab\psi_j)$ is a 1-morphism over $S$, in the sense of \S\ref{kt31}. Theorem \ref{kt3thm3} shows that $\ti\Phi_{ij}$ is a coordinate change over $S$, in the sense of~\S\ref{kt32}.
\label{kt3ex2}
\end{ex}

We will show that the elements $\ga_{rqp}^\al\in\Ga_p$ in Definition \ref{kt2def3}(b) correspond in the setting of \S\ref{kt3} to a 2-morphism~$\La_{rqp}\colon \ti\Phi_{qp}\ci\ti\Phi_{rq}\Ra\ti\Phi_{rp}$.

\begin{ex}{\bf(i)} In the Fukaya--Oh--Ohta--Ono theory \cite{Fuka,FOOO1,FOOO2,FOOO3,FOOO4,FOOO5,FOOO6,FOOO7,FOOO8,FOOO9,FuOn}, one often relates two FOOO coordinate changes in the following way. Let $\Phi_{ij}=(V_{ij},\ab h_{ij},\ab \vp_{ij},\ab\hat\vp_{ij}),\ab\Phi_{ij}'=(V_{ij}',\ab h_{ij}',\ab\vp_{ij}',\hat\vp_{ij}')\colon (V_i,E_i,\Ga_i,s_i,\psi_i)\ra(V_j,\ab E_j,\ab\Ga_j,\ab s_j,\ab\psi_j)$ be FOOO coordinate changes over $S$. Suppose there exists $\ga\in\Ga_j$ such that
\e
h_{ij}=\ga\cdot h_{ij}'\cdot\ga^{-1}, \quad
\phi_{ij}=\ga\cdot \phi_{ij}',\quad\text{and}\quad 
\hat\phi_{ij}=\ga\cdot\hat\phi_{ij}',
\label{kt3eq10}
\e
where the second and third equations hold on~$\dot V_{ij}:=V_{ij}\cap V_{ij}'$.

Let $\ti\Phi_{ij},\ti\Phi_{ij}'\colon (V_i,E_i,\Ga_i,s_i,\psi_i)\ra (V_j,\ab E_j,\ab\Ga_j,\ab s_j,\ab\psi_j)$ be the 1-morphisms in the sense of \S\ref{kt31} corresponding to $\Phi_{ij},\Phi_{ij}'$ in Example \ref{kt3ex2}. Set $\dot P_{ij}=\dot V_{ij}\t\Ga_j\subseteq P_{ij}$. Define $\la_{ij}\colon \dot P_{ij}=\dot V_{ij}\t\Ga_j\ra V_{ij}'\t\Ga_j=P_{ij}'$ by $\la_{ij}\colon (v,\ga')\mapsto(v,\ga'\ga)$, and $\hat\la_{ij}=0$. Then $(\dot P_{ij},\la_{ij},\hat\la_{ij})$ satisfies Definition \ref{kt3def4}(a)--(c), so we have defined a 2-morphism $\La_{ij}=[\dot P_{ij},\la_{ij},\hat\la_{ij}]\colon \ti\Phi_{ij}\Ra\ti\Phi_{ij}'$, in the sense of~\S\ref{kt31}.
\smallskip

\noindent{\bf(ii)} This enables us to interpret Definition \ref{kt2def3}(b) in terms of a 2-morphism. In the situation of Definition \ref{kt2def3}(b), the composition of the FOOO coordinate changes $\Phi_{rq},\Phi_{qp}$ is $\Phi_{qp}\ci\Phi_{qp}=\bigl(\vp_{rq}^{-1}(V_{qp}),h_{qp}\ci h_{rq},\vp_{qp} \ci \vp_{rq}\vert_{\vp_{rq}^{-1}(V_{qp})},\vp_{rq}^*(\hat\vp_{qp}) \ci \hat\vp_{rq}\vert_{\vp_{rq}^{-1}(V_{qp})}\bigr)$. Thus, \eq{kt2eq3} relates $\Phi_{qp}\ci\Phi_{rq}$ to $\Phi_{rp}$ in the same way that \eq{kt3eq10} relates $\Phi_{ij}$ to $\Phi_{ij}'$, except for allowing $\ga_{rqp}$ to vary on different connected components. Hence, if $\ti\Phi_{rq},\ti\Phi_{qp},\ti\Phi_{rp}$ are the coordinate changes in the sense of \S\ref{kt32} associated to $\Phi_{rq},\ab\Phi_{qp},\ab\Phi_{rp}$ in Example \ref{kt3ex2}, then the method of {\bf(i)} defines a 2-morphism $\La_{pqr}\colon \ti\Phi_{qp}\ci\ti\Phi_{rq}\Ra\ti\Phi_{rp}$, in the sense of~\S\ref{kt31}.
\smallskip

\noindent{\bf(iii)} In the situation of Definition \ref{kt2def3}(b), suppose $v\in (\vp_{rq}^{-1}(V_{qp})\cap V_{rp})^\al$ is generic. Then $\Stab_{\Ga_r}(v)=\{1\}$, as $\Ga_r$ acts (locally) effectively on $V_r$ by Definition \ref{kt2def1}(c). Hence $\Stab_{\Ga_p}(\vp_{rp}(v))=\{1\}$ by Definition \ref{kt2def2}(g). Therefore the point $\ga_{rqp}^\al\cdot \vp_{rp}(v)=\vp_{qp}\ci\vp_{rq}(v)$ in $V_p$ determines $\ga_{rqp}^\al$ in $\Ga_p$. So the second equation of \eq{kt2eq3} determines $\ga_{rqp}^\al\in\Ga_p$ uniquely, provided it exists. Thus the 2-morphism $\La_{pqr}\colon \ti\Phi_{qp}\ci\ti\Phi_{rq}\Ra\ti\Phi_{rp}$ in {\bf(ii)} is also determined uniquely.
\label{kt3ex3}
\end{ex}

\subsection{\texorpdfstring{Relation to McDuff and Wehrheim's work}{Relation to McDuff and Wehrheim\textquoteright s work}}
\label{kt35}

We can also connect our definitions in \S\ref{kt31}--\S\ref{kt33} to MW Kuranishi neighbourhoods and MW coordinate changes from \S\ref{kt23}.
As in Example \ref{kt3ex1}, by an abuse of notation we will regard MW Kuranishi neighbourhoods as examples of our Kuranishi neighbourhoods in~\S\ref{kt31}. We relate MW coordinate changes to ours:

\begin{ex}  Let $\Phi_{BC}=(\ti V_{BC},\rho_{BC},\varpi_{BC},\hat\vp_{BC})\colon (V_B,E_B,\Ga_B,s_B,\psi_B)\ra(V_C,\ab E_C,\ab\Ga_C,\ab s_C,\ab\psi_C)$ be an MW coordinate change over $S$, as in Definition \ref{kt2def7}. Regard $(V_B,E_B,\Ga_B,s_B,\psi_B),(V_C,E_C,\Ga_C,s_C,\psi_C)$ as Kuranishi neighbourhoods in the sense of~\S\ref{kt31}. 

Set $P_{BC}=\ti V_{BC}\t\Ga_B$. Let $\Ga_B$ act on $P_{BC}$ by $\ga_B\colon (v,\ga)\mapsto (v,\ga_B\ga)$. Let $\Ga_C$ act on $P_{BC}$ by $\ga_C\colon (v,\ga)\mapsto (\ga_C\cdot v,\ga\rho_{BC}(\ga_C)^{-1})$. Define $\pi_{BC}\colon P_{BC}\ra V_B$ and $\phi_{BC}\colon P_{BC}\ra V_C$ by $\pi_{BC}\colon (v,\ga)\mapsto\ga\cdot\varpi_{BC}(v)$ and $\phi_{BC}\colon (v,\ga)\mapsto v$. Then $\pi_{BC}$ is $\Ga_B$-equivariant and $\Ga_C$-invariant, and $\phi_{BC}$ is $\Ga_B$-invariant and $\Ga_C$-equivariant.

Define $\hat\phi_{BC}:\pi_{BC}^*(V_B\t E_B)\ra\phi_{BC}^*(V_C\t E_C)$, as a morphism of trivial vector bundles with fibres $E_B,E_C$ on $P_{BC}=\ti V_{BC}\t\Ga_B$, by $\hat\phi_{BC}\vert_{\ti V_{BC}\t\{\ga\}}=\hat\vp_{BC}\ci(\ga^{-1}\cdot -)$ for each $\ga\in\Ga_B$. It is easy to see that $\ti\Phi_{BC}=(P_{BC},\pi_{BC},\phi_{BC},\hat\phi_{BC})\colon \ab(V_B,E_B,\Ga_B,s_B,\psi_B)\ra (V_C,\ab E_C,\ab\Ga_C,\ab s_C,\ab\psi_C)$ is a 1-morphism over $S$, in the sense of \S\ref{kt31}. Combining Definition \ref{kt2def7}(g) and Theorem \ref{kt3thm3} shows that $\ti\Phi_{BC}$ is a coordinate change over $S$, in the sense of~\S\ref{kt32}.
\label{kt3ex4}
\end{ex}

We relate Definition \ref{kt2def8}(d) to 2-morphisms in~\S\ref{kt31}:

\begin{ex}  In the situation of Definition \ref{kt2def8}(d), let $\ti\Phi_{BC},\ti\Phi_{BD},\ti\Phi_{CD}$ be the coordinate changes in the sense of \S\ref{kt32} associated to the MW coordinate changes $\Phi_{BC},\Phi_{BD},\Phi_{CD}$ in Example \ref{kt3ex4}. The composition coordinate change $\ti\Phi_{CD}\ci\ti\Phi_{BC}=(P_{BCD},\pi_{BCD},\phi_{BCD},\hat\phi_{BCD})$ from Definition \ref{kt3def5} has
\e
\begin{split}
P_{BCD}&=\bigl[(\ti V_{BC}\t\Ga_B)\t_{V_C}(\ti V_{CD}\t\Ga_C)\bigr]\big/\Ga_C\\
&\cong (\ti V_{BC}\t_{V_C}\ti V_{CD})\t\Ga_B \cong \varpi_{CD}^{-1}(\ti V_{BC})\t\Ga_B.
\end{split}
\label{kt3eq11}
\e
Define $\dot P_{BCD}$ to be the open subset of $P_{BCD}$ identified with $\ti V_{BCD}\t\Ga_B$ by \eq{kt3eq11}, and $\la_{BCD}\colon \dot P_{BCD}\ra P_{BD}=\ti V_{BD}\t\Ga_B$ to be the map identified by \eq{kt3eq11} with the inclusion $\ti V_{BCD}\t\Ga_B\hookra \ti V_{BD}\t\Ga_B$, and $\hat\la_{BCD}=0$. Then as in Example \ref{kt3ex3}(i), we can show that $(\dot P_{BCD},\la_{BCD},\hat\la_{BCD})$ satisfies Definition \ref{kt3def4}(a)--(c), so we have defined a 2-morphism $\La_{BCD}=[\dot P_{BCD},\la_{BCD},\hat\la_{BCD}]\colon \ti\Phi_{CD}\ci\ti\Phi_{BC}\Ra\ti\Phi_{BD}$ on $S_{BCD}=\Im\psi_B\cap\Im\psi_C\cap\Im\psi_D$, in the sense of~\S\ref{kt31}.
\label{kt3ex5}
\end{ex}

\subsection{Relation to d-orbifolds}
\label{kt36}

The author has developed a theory of Derived Differential Geometry \cite{Joyc2,Joyc3,Joyc4}, involving `d-manifolds' and `d-orbifolds', derived versions of smooth manifolds and orbifolds, where `derived' is in the sense of Derived Algebraic Geometry. D-manifolds and d-orbifolds form strict 2-categories $\dMan,\dOrb$. For other definitions of derived manifolds, which form $\iy$-categories, see Spivak \cite{Spiv} and Borisov--Noel \cite{BoNo}. Borisov \cite{Bori} explains how the d-manifolds of \cite{Joyc2,Joyc3,Joyc4} are related to the derived manifolds of Spivak \cite{Spiv} and Borisov--Noel~\cite{BoNo}.

D-manifolds are a full 2-subcategory of the strict 2-category of `d-spaces' $\dSpa$, where a d-space $\bX=(X,\O_X^\bu)$ is a topological space $\bX$ with a sheaf $\O_X^\bu$ of `square zero dg $C^\iy$-rings', satisfying some conditions. Here a square zero dg $C^\iy$-ring is a 2-categorical, differential-geometric analogue of the commutative differential graded algebras (cdgas) used as basic objects in Derived Algebraic Geometry. D-orbifolds are a Deligne--Mumford stack version of d-manifolds.

The initial definitions of d-manifolds and d-orbifolds broadly follow those of schemes and stacks in algebraic geometry, and look nothing like \S\ref{kt31}--\S\ref{kt32}. But in order to describe d-manifolds and d-orbifolds in classical differential-geometric language, the author \cite[\S 3 \& \S 10]{Joyc4} proves results on `standard model' d-manifolds and d-orbifolds, and `standard model' 1- and 2-morphisms. 

In terms of the ideas of \S\ref{kt31}, a `standard model' d-manifold $\bS_{V,E,s}$ or d-orbifold $\bS_{V,E,\Ga,s}$ is constructed from a Kuranishi neighbourhood $(V,E,\Ga,s,\psi)$, with $\Ga=\{1\}$ in the d-manifold case, and `standard model' 1- and 2-morphisms are constructed from the data in 1- and 2-morphisms of Kuranishi neighbourhoods. By comparing \S\ref{kt31}--\S\ref{kt32} with results in \cite[\S 3 \& \S 10]{Joyc4}, we can deduce:

\begin{thm}{\bf(a)} Suppose $(V,E,\Ga,s,\psi)$ is a global Kuranishi neighbourhood on a topological space $X,$ where `global' means $\Im\psi=X$. Then we can define a `standard model' d-orbifold\/ $\bS_{V,E,\Ga,s},$ with underlying topological space\/ $X$.
\smallskip

\noindent{\bf(b)} Let\/ $(V_i,E_i,\Ga_i,s_i,\psi_i),(V_j,E_j,\Ga_j,s_j,\psi_j)$ be global Kuranishi neighbourhoods on $X,Y,$ and\/ $f\colon X\ra Y$ be continuous, and\/ $\Phi_{ij}\colon (V_i,\ab E_i,\ab\Ga_i,\ab s_i,\ab\psi_i)\ra (V_j,\ab E_j,\ab\Ga_j,\ab s_j,\ab\psi_j)$ be a $1$-morphism of Kuranishi neighbourhoods over\/ $(X,f)$. Then we can define a  `standard model'\/ $1$-morphism of d-orbifolds in $\dOrb$
\begin{equation*}
\smash{\bS_{\Phi_{ij}}\colon \bS_{V_i,E_i,\Ga_i,s_i}\longra \bS_{V_j,E_j,\Ga_j,s_j},}
\end{equation*}
with underlying continuous map\/~$f$.
\smallskip

\noindent{\bf(c)} Let\/ $\bs f\colon \bS_{V_i,E_i,\Ga_i,s_i}\ra \bS_{V_j,E_j,\Ga_j,s_j}$ be a $1$-morphism in $\dOrb$ between standard model d-orbifolds. Then $\bs f$ is $2$-isomorphic in $\dOrb$ to some standard model\/ $1$-morphism $\bS_{\Phi_{ij}},$ as in\/~{\bf(b)}.
\smallskip

\noindent{\bf(d)} Suppose $\Phi_{ij},\Phi_{ij}'\colon (V_i,E_i,\Ga_i,s_i,\psi_i)\ra (V_j,E_j,\Ga_j,s_j,\psi_j)$ are $1$-morphisms of Kuranishi neighbourhoods over\/ $(X,f),$ with associated\/ $1$-morphisms of d-orbifolds
$\bS_{\Phi_{ij}},\bS_{\Phi_{ij}'}$. Then there is a canonical\/ $1$-$1$ correspondence between $2$-morphisms $\La\colon \Phi_{ij}\Ra\Phi_{ij}'$ of Kuranishi neighbourhoods over $(X,f),$ and\/ $2$-morphisms $\bs\la\colon \bS_{\Phi_{ij}}\ab\Ra\bS_{\Phi_{ij}'}$ of d-orbifolds.
\smallskip

\noindent These maps from Kuranishi neighbourhoods, $1$- and\/ $2$-morphisms to d-orbifolds, $1$- and\/ $2$-morphisms are compatible with the rest of the $2$-category structures.
\label{kt3thm4}
\end{thm}

This implies that the full strict 2-subcategory $\SMod\subset\dOrb$ with objects `standard model' d-orbifolds is equivalent to the weak 2-category $\GKN$ of global Kuranishi neighbourhoods in Theorem~\ref{kt3thm1}.

In fact the author deliberately wrote \S\ref{kt31}--\S\ref{kt32} by translating facts about d-orbifolds into Kuranishi-style, differential-geometric language, which makes Theorem \ref{kt3thm4} almost a tautology. But the facts about `standard model' d-manifolds and d-orbifolds on which Theorem \ref{kt3thm4} is based are genuine theorems, not just simple consequences of the definition of d-orbifolds.

\section{The weak 2-category of Kuranishi spaces}
\label{kt4}

We now define our 2-category of Kuranishi spaces $\Kur$. The material of this section, and the proofs of the theorems, have a different character to those of \S\ref{kt3}. Section \ref{kt3} was largely about differential geometry: in Kuranishi neighbourhoods $(V_i,E_i,\Ga_i,s_i,\psi_i)$ it really mattered that $V_i$ is a manifold, and so on, and the proofs of Theorems \ref{kt3thm1}, \ref{kt3thm2} and \ref{kt3thm3} use many properties of manifolds.

This section involves a lot of stack theory from algebraic geometry, but very little differential geometry. When we define Kuranishi structures in \S\ref{kt41}, in Kuranishi neighbourhoods $(V_i,E_i,\Ga_i,s_i,\psi_i)$ we do not care that $V_i$ is a manifold, it only matters that $(V_i,E_i,\Ga_i,s_i,\psi_i)$ is an object in a 2-category $\KN_S(X)$ for open $S\subseteq\Im\psi_i\subseteq X$, which has a stack property as in Theorem~\ref{kt3thm2}.

Because of this, the framework below works whenever we have some class of `charts' on a topological space satisfying the analogues of Theorems \ref{kt3thm1} and \ref{kt3thm2}, so for instance we can easily define Kuranishi spaces with boundary and corners, or Kuranishi spaces associated to other categories of manifolds.

All of this section, including proofs of quoted results, comes from~\cite{Joyc5,Joyc6}.

\subsection{Kuranishi spaces, 1-morphisms, and 2-morphisms}
\label{kt41}

Here is one of the main definitions of the paper:

\begin{dfn} Let $X$ be a Hausdorff, second countable topological space (not necessarily compact), and $n\in\Z$. A {\it Kuranishi structure\/ $\cK$ on $X$ of virtual dimension\/} $n$ is data $\cK=\bigl(I,(V_i,E_i,\Ga_i,s_i,\psi_i)_{i\in I}$, $\Phi_{ij,\;i,j\in I}$, $\La_{ijk,\; i,j,k\in I}\bigr)$, where:
\begin{itemize}
\setlength{\itemsep}{0pt}
\setlength{\parsep}{0pt}
\item[(a)] $I$ is an indexing set (not necessarily finite).
\item[(b)] $(V_i,E_i,\Ga_i,s_i,\psi_i)$ is a Kuranishi neighbourhood on $X$ for each $i\in I$, with $\dim V_i-\rank E_i=n$.
\item[(c)] $\Phi_{ij}=(P_{ij},\pi_{ij},\phi_{ij},\hat\phi_{ij})\colon (V_i,E_i,\ab\Ga_i,\ab s_i,\ab\psi_i)\ab\ra (V_j,E_j,\Ga_j,s_j,\psi_j)$ is a coordinate change for all $i,j\in I$ (as usual, defined over $S=\Im\psi_i\cap\Im\psi_j$).
\item[(d)] $\La_{ijk}=[\dot P_{ijk},\la_{ijk},\hat\la_{ijk}]\colon \Phi_{jk}\ci\Phi_{ij}\Ra\Phi_{ik}$ is a 2-morphism for all $i,j,k\in I$ (as usual, defined over $S=\Im\psi_i\cap\Im\psi_j\cap\Im\psi_k$).
\item[(e)] $\bigcup_{i\in I}\Im\psi_i=X$. 
\item[(f)] $\Phi_{ii}=\id_{(V_i,E_i,\Ga_i,s_i,\psi_i)}$ for all $i\in I$.
\item[(g)] $\La_{iij}=\bs\be_{\Phi_{ij}}$ and $\La_{ijj}=\bs\ga_{\Phi_{ij}}$ for all $i,j\in I$, for $\bs\be_{\Phi_{ij}},\bs\ga_{\Phi_{ij}}$ as in \eq{kt3eq5}.
\item[(h)] The following diagram of 2-morphisms over $S=\Im\psi_i\ab\cap\ab\Im\psi_j\ab\cap\ab\Im\psi_k\ab\cap\ab\Im\psi_l$ commutes for all $i,j,k,l\in I$, for $\bs\al_{\Phi_{kl},\Phi_{jk},\Phi_{ij}}$ as in \eq{kt3eq4}:
\begin{equation*}
\xymatrix@C=90pt@R=15pt{
*+[r]{(\Phi_{kl}\ci\Phi_{jk})\ci\Phi_{ij}} \ar@{=>}[d]^{\bs\al_{\Phi_{kl},\Phi_{jk},\Phi_{ij}}}
\ar@{=>}[rr]_(0.53){\La_{jkl}*\id_{\Phi_{ij}}} && *+[l]{\Phi_{jl}\ci\Phi_{ij}} \ar@{=>}[d]_{\La_{ijl}}  
\\
*+[r]{\Phi_{kl}\ci(\Phi_{jk}\ci\Phi_{ij})} \ar@{=>}[r]^(0.65){\id_{\Phi_{kl}}*\La_{ijk}}
& \Phi_{kl}\ci\Phi_{ik} \ar@{=>}[r]^{\La_{ikl} } & *+[l]{\Phi_{il}.\!} }
\end{equation*}
\end{itemize}
We call $\bX=(X,\cK)$ a {\it Kuranishi space}, of {\it virtual dimension\/}~$\vdim\bX=n$.

When we write $x\in\bX$, we mean that $x\in X$. We can define {\it orientations\/} on Kuranishi spaces in a very similar way to Definition~\ref{kt2def4}.
\label{kt4def1}
\end{dfn}

\begin{ex} Let $V$ be a manifold, $E\ra V$ a vector bundle, $\Ga$ a finite group with a smooth action on $V$ and a compatible action on $E$ preserving the vector bundle structure, and $s\colon V\ra E$ a $\Ga$-equivariant smooth section. Set $X=s^{-1}(0)/\Ga$, with the quotient topology induced from the closed subset $s^{-1}(0)\subseteq V$. Then $X$ is Hausdorff and second countable, as $V$ is and $\Ga$ is finite. 

Define a Kuranishi structure $\cK=\bigl(\{0\},(V_0,E_0,\Ga_0,s_0,\psi_0),\Phi_{00},\La_{000}\bigr)$ on $X$ with indexing set $I=\{0\}$, one Kuranishi neighbourhood $(V_0,E_0,\Ga_0,s_0,\psi_0)$ with $V_0=V$, $E_0=E$, $\Ga_0=\Ga$, $s_0=s$ and $\psi_0=\id_X$, one coordinate change $\Phi_{00}=\id_{(V_0,E_0,\Ga_0,s_0,\psi_0)}$, and one 2-morphism $\La_{000}=\id_{\Phi_{00}}$. Then $\bX=(X,\cK)$ is a Kuranishi space, with $\vdim\bX=\dim V-\rank E$. We write~$\bS_{V,E,\Ga,s}=\bX$.
\label{kt4ex1}
\end{ex}

We will need notation to distinguish Kuranishi neighbourhoods, coordinate chan\-ges, and 2-morphisms on different Kuranishi spaces. We will often use the following notation for Kuranishi spaces $\bW,\bX,\bY,\bZ$:
\ea
\bW&=(W,\cH),& \cH&=\bigl(H,(T_h,C_h,\Al_i,q_h,\vp_h)_{h\in H},\; \Si_{hh'}=(O_{hh'},\pi_{hh'},\si_{hh'},
\nonumber\\
&& {}\hskip -60pt \hat\si_{hh'})_{h,h'\in H}&,\; \Io_{hh'h''}=[\dot O_{hh'h''},\io_{hh'h''},\hat\io_{hh'h''}]_{h,h',h''\in H}\bigr),
\label{kt4eq1}\\
\bX&=(X,\cI),& \cI&=\bigl(I,(U_i,D_i,\Be_i,r_i,\chi_i)_{i\in I},\;
\Tau_{ii'}=(P_{ii'},\pi_{ii'},\tau_{ii'},
\nonumber\\
&&{}\hskip -60pt \hat\tau_{ii'})_{i,i'\in I}&,\; \Ka_{ii'i''}=[\dot P_{ii'i''},\ka_{ii'i''},\hat\ka_{ii'i''}]_{i,i',i''\in I}\bigr),
\label{kt4eq2}\\
\bY&=(Y,\cJ),& \cJ&=\bigl(J,(V_j,E_j,\Ga_j,s_j,\psi_j)_{j\in J},\; \Up_{jj'}=(Q_{jj'},\pi_{jj'},\up_{jj'},
\nonumber\\
&&{}\hskip -60pt \hat\up_{jj'})_{j,j'\in J}&,\; \La_{jj'j''}=[\dot Q_{jj'j''},\la_{jj'j''},\hat\la_{jj'j''}]_{j,j',j''\in J}\bigr),
\label{kt4eq3}\\
\bZ&=(Z,\cK),& \cK&=\bigl(K,(W_k,F_k,\De_k,t_k,\om_k)_{k\in K},\; \Phi_{kk'}=(R_{kk'},\pi_{kk'},\phi_{kk'},
\nonumber\\
&&{}\hskip -60pt \hat\phi_{kk'})_{k,k'\in K}&,\; \Mu_{kk'k''}=[\dot R_{kk'k''},\mu_{kk'k''},\hat\mu_{kk'k''}]_{k,k',k''\in K}\bigr).
\label{kt4eq4}
\ea

Next we define 1- and 2-morphisms of Kuranishi spaces. Note a possible confusion: we will be defining 1-morphisms of Kuranishi spaces $\bs f,\bs g\colon \bX\ra\bY$ and 2-morphisms of Kuranishi spaces $\bs\eta\colon \bs f\Ra\bs g$, but these will be built out of 1-morphisms of Kuranishi neighbourhoods $\bs f_{ij},\bs g_{ij}\colon (U_i,D_i,\Be_i,r_i,\chi_i)\ra (V_j,E_j,\Ga_j,s_j,\psi_j)$ and 2-morphisms of Kuranishi neighbourhoods $\bs\eta_{ij}\colon \bs f_{ij}\Ra\bs g_{ij}$ in the sense of \S\ref{kt31}, so `1-morphism' and `2-morphism' can mean two things.

\begin{dfn} Let $\bX=(X,\cI)$ and $\bY=(Y,\cJ)$ be Kuranishi spaces, with notation \eq{kt4eq2}--\eq{kt4eq3}. A 1-{\it morphism of Kuranishi spaces\/} $\bs f\colon \bX\ra\bY$ is data
\begin{equation*}
\bs f=\bigl(f,\bs f_{ij,\;i\in I,\; j\in J},\; \bs F_{ii',\;i,i'\in I}^{j,\; j\in J},\; \bs F_{i,\;i\in I}^{jj',\; j,j'\in J}\bigr),
\end{equation*}
satisfying the conditions:
\begin{itemize}
\setlength{\itemsep}{0pt}
\setlength{\parsep}{0pt}
\item[(a)] $f\colon X\ra Y$ is a continuous map.
\item[(b)] $\bs f_{ij}=(P_{ij},\pi_{ij},f_{ij},\hat f_{ij})\colon (U_i,D_i,\Be_i,r_i,\chi_i)\ra (V_j,E_j,\Ga_j,s_j,\psi_j)$ is a 1-mor\-ph\-ism of Kuranishi neighbourhoods over $f$ for all $i\in I$, $j\in J$ (defined over $S=\Im\chi_i\cap f^{-1}(\Im\psi_j)$, as usual).
\item[(c)] $\bs F_{ii'}^j=[\dot P_{ii'}^j,F_{ii'}^j,\hat F_{ii'}^j]\colon \bs f_{i'j}\ci\Tau_{ii'}\Ra \bs f_{ij}$ is a 2-morphism over $f$ for all $i,i'\in I$ and $j\in J$ (defined over $S=\Im\chi_i\cap\Im\chi_{i'}\cap f^{-1}(\Im\psi_j)$).
\item[(d)] $\bs F_i^{jj'}=[\dot P_i^{jj'},F_i^{jj'},\hat F_i^{jj'}]\colon \Up_{jj'}\ci\bs f_{ij}\Ra \bs f_{ij'}$ is a 2-morphism over $f$ for all $i\in I$ and $j,j'\in J$ (defined over $S=\Im\chi_i\cap f^{-1}(\Im\psi_j\cap\Im\psi_{j'})$).
\item[(e)] $\bs F_{ii}^j=\bs\be_{\bs f_{ij}}$ and $\bs F_i^{jj}=\bs\ga_{\bs f_{ij}}$ for all $i\in I$, $j\in J$, for $\bs\be_{\bs f_{ij}},\bs\ga_{\bs f_{ij}}$ as in~\eq{kt3eq5}.
\item[(f)] The following commutes for all $i,i',i''\in I$ and $j\in J$: 
\begin{equation*}
\xymatrix@C=91pt@R=15pt{
*+[r]{(\bs f_{i''j}\ci\Tau_{i'i''})\ci\Tau_{ii'}} \ar@{=>}[d]^{\bs\al_{\bs f_{i''j},\Tau_{i'i''},\Tau_{ii'}}}
\ar@{=>}[rr]_(0.53){\bs F_{i'i''}^j*\id_{\Tau_{ii'}}} && *+[l]{\bs f_{i'j}\ci\Tau_{ii'}} \ar@{=>}[d]_{\bs F_{ii'}^j}  
\\
*+[r]{\bs f_{i''j}\ci(\Tau_{i'i''}\ci\Tau_{ii'})} \ar@{=>}[r]^(0.65){\id_{\bs f_{i''j}}*\Ka_{ii'i''}}
& \bs f_{i''j}\ci \Tau_{ii''} \ar@{=>}[r]^{\bs F_{ii''}^j} & *+[l]{\bs f_{ij}.\!} }
\end{equation*}
\item[(g)] The following commutes for all $i,i'\in I$ and $j,j'\in J$:
\begin{equation*}
\xymatrix@C=91pt@R=15pt{
*+[r]{(\Up_{jj'}\ci\bs f_{i'j})\ci\Tau_{ii'}} \ar@{=>}[d]^{\bs\al_{\Up_{jj'},\bs f_{i'j},\Tau_{ii'}}}
\ar@{=>}[rr]_(0.53){\bs F_{i'}^{jj'} *\id_{\Tau_{ii'}}} && *+[l]{\bs f_{i'j'}\ci\Tau_{ii'}} \ar@{=>}[d]_{\bs F_{ii'}^{j'}}  
\\
*+[r]{\Up_{jj'}\ci(\bs f_{i'j}\ci\Tau_{ii'})} \ar@{=>}[r]^(0.65){\id_{\Up_{jj'}}*\bs F_{ii'}^j}
& \Up_{jj'}\ci\bs f_{ij} \ar@{=>}[r]^{\bs F_i^{jj'}} & *+[l]{\bs f_{ij'}.\!} }
\end{equation*}
\item[(h)] The following commutes for all $i\in I$ and $j,j',j''\in J$:
\begin{equation*}
\xymatrix@C=89pt@R=15pt{
*+[r]{(\Up_{j'j''}\ci\Up_{jj'})\ci\bs f_{ij}} \ar@{=>}[d]^{\bs\al_{\Up_{j'j''},\Up_{jj'},\bs f_{ij}}}
\ar@{=>}[rr]_(0.53){\La_{jj'j''}*\id_{\bs f_{ij}}} && *+[l]{\Up_{jj''}\ci\bs f_{ij}} \ar@{=>}[d]_{\bs F_i^{jj''}}  
\\
*+[r]{\Up_{j'j''}\ci(\Up_{jj'}\ci\bs f_{ij})} \ar@{=>}[r]^(0.68){\id_{\Up_{j'j''}}*\bs F_i^{jj'}}
& \Up_{j'j''}\ci\bs f_{ij'} \ar@{=>}[r]^{\bs F_i^{j'j''}} & *+[l]{\bs f_{ij''}.\!} }
\end{equation*}
\end{itemize}

If $x\in\bX$ (i.e. $x\in X$), we will write $\bs f(x)=f(x)\in\bY$.

When $\bY=\bX$, define the {\it identity\/ $1$-morphism\/} $\bs\id_\bX\colon \bX\ra\bX$ by
\begin{equation*}
\bs\id_\bX=\bigl(\id_X,\Tau_{ij,\; i,j\in I},\; \Ka_{ii'j,\;i,i'\in I}^{\;\;\;\;\;\;\;\; j\in I},\; \Ka_{ijj',\;i\in I}^{\;\;\; j,j'\in I}\bigr),
\end{equation*}
Then Definition \ref{kt4def1}(h) implies that (f)--(h) above hold.
\label{kt4def2}
\end{dfn}

\begin{dfn} Let $\bX=(X,\cI)$ and $\bY=(Y,\cJ)$ be Kuranishi spaces, with notation as in \eq{kt4eq2}--\eq{kt4eq3}, and $\bs f,\bs g\colon \bX\ra\bY$ be 1-morphisms. Suppose the continuous maps $f,g\colon X\ra Y$ in $\bs f,\bs g$ satisfy $f=g$. A 2-{\it morphism of Kuranishi spaces\/} $\bs\eta\colon \bs f\Ra\bs g$ is data $\bs\eta=\bigl(\bs\eta_{ij,\; i\in I,\; j\in J}\bigr)$, where $\bs\eta_{ij}=[\dot P_{ij},\eta_{ij},\hat\eta_{ij}]\colon \bs f_{ij}\Ra\bs g_{ij}$ is a 2-morphism of Kuranishi neighbourhoods over $f=g$ (defined over $S=\Im\chi_i\cap f^{-1}(\Im\psi_j)$, as usual), satisfying the conditions:
\begin{itemize}
\setlength{\itemsep}{0pt}
\setlength{\parsep}{0pt}
\item[(a)] $\bs G_{ii'}^j\od(\bs\eta_{i'j}*\id_{\Tau_{ii'}})=\bs\eta_{ij}\od\bs F_{ii'}^j\colon \bs f_{i'j}\ci\Tau_{ii'}\Ra\bs g_{ij}$ for all $i,i'\in I$, $j\in J$. 
\item[(b)] $\bs G_i^{jj'}\od(\id_{\Up_{jj'}}*\bs\eta_{ij})=\bs\eta_{ij'}\od\bs F_i^{jj'}\colon \Up_{jj'}\ci\bs f_{ij}\Ra\bs g_{ij'}$ for all $i\in I$, $j,j'\in J$.\end{itemize}
Note that by definition, 2-morphisms $\bs\eta\colon \bs f\Ra\bs g$ only exist if $f=g$.

If $\bs f=\bs g$, the {\it identity\/ $2$-morphism\/} is $\bs\id_{\bs f}=\bigl(\id_{\bs f_{ij},\; i\in I,\; j\in J}\bigr)\colon \bs f\Ra\bs f$.
\label{kt4def3}
\end{dfn}

\subsection{Making Kuranishi spaces into a 2-category}
\label{kt42}

We will make Kuranishi spaces into a weak 2-category. We have already defined the objects, 1-morphisms, 2-morphisms, and identity 1- and 2-morphisms in \S\ref{kt41}. As for \S\ref{kt32}, it remains to explain:
\begin{itemize}
\setlength{\itemsep}{0pt}
\setlength{\parsep}{0pt}
\item Composition of 1-morphisms;
\item Horizontal composition of 2-morphisms;
\item Vertical composition of 2-morphisms; and
\item Coherence 2-isomorphisms $\bs\al_{\bs g,\bs f,\bs e}\colon (\bs g\ci\bs f)\ci\bs e\Ra\bs g\ci (\bs f\ci\bs e)$, $\bs\be_{\bs f}\colon \bs f\ci\bs\id_\bX\Ra\bs f$, and $\bs\ga_{\bs f}\colon \bs\id_\bY\ci\bs f\Ra\bs f$, as in \eq{ktAeq5} and~\eq{ktAeq7}.
\end{itemize}

We will define composition of 1-morphisms in Proposition \ref{kt4prop1} and Definition \ref{kt4def4}. The definition involves some interesting issues, so we discuss these first. Suppose $\bs f\colon \bX\ra\bY$, $\bs g\colon \bY\ra\bZ$ are 1-morphisms of Kuranishi spaces. We want to define the composition $\bs g\ci\bs f\colon \bX\ra\bZ$. Use notation \eq{kt4eq2}--\eq{kt4eq4} for $\bX,\bY,\bZ$. Then Kuranishi neighbourhoods on $X$ are parametrized by $i\in I$, where $\{\Im\chi_i\colon i\in I\}$ is an open cover of $X$, and similarly $\{\Im\psi_j\colon j\in J\}$ is an open cover of $Y$, and $\{\Im\om_k\colon k\in K\}$ is an open cover of~$Z$.

To make $\bs g\ci\bs f$, for each $i\in I$ and $k\in K$ we need a 1-morphism of Kuranishi neighbourhoods $(\bs g\ci\bs f)_{ik}$ defined on $\Im\chi_i\cap (g\ci f)^{-1}(\Im\om_k)\subseteq X$. But we actually have 1-morphisms $\bs g_{jk}\ci\bs f_{ij}$ defined on $\Im\chi_i\cap f^{-1}(\psi_j)\cap (g\ci f)^{-1}(\Im\om_k)\subseteq X$ for all $j\in J$. We will construct $(\bs g\ci\bs f)_{ik}$ by gluing together the $\bs g_{jk}\ci\bs f_{ij}$ for all $j\in J$ using the stack property of Kuranishi neighbourhoods, Theorem~\ref{kt3thm2}.

The stack property is crucial, as without it we could not define composition of 1-morphisms of Kuranishi spaces, and could not make Kuranishi spaces into a 2-category. For FOOO Kuranishi spaces in \S\ref{kt21}, or MW weak Kuranishi atlases in \S\ref{kt23}, one could write down notions of morphisms similar to Definition \ref{kt4def2}. However, as FOOO/MW coordinate changes do not form sheaves or stacks, there seems no obvious way to define compositions of such morphisms.

In the next proposition, taken from \cite[\S 4.3]{Joyc5}, \cite[\S 6.2]{Joyc6}, part (a) constructs candidates $\bs h$ for $\bs g\ci\bs f$, part (b) shows such $\bs h$ are unique up to canonical 2-isomorphism, and part (c) that $\bs f,\bs g$ are allowed candidates for~$\bs\id_\bY\ci\bs f,\bs g\ci\bs\id_\bY$.

\begin{prop}{\bf(a)} Let\/ $\bX=(X,\cI),\bY=(Y,\cJ),\bZ=(Z,\cK)$ be Kuranishi spaces with notation {\rm\eq{kt4eq2}--\eq{kt4eq4},} and\/ $\bs f\colon \bX\ra\bY,$ $\bs g\colon \bY\ra\bZ$ be $1$-morphisms, with\/ $\bs f=\bigl(f,\bs f_{ij},\bs F_{ii'}^j,\bs F_i^{jj'}\bigr),$ $\bs g=\bigl(g,\bs g_{jk},\bs G_{jj'}^k,\bs G_j^{kk'}\bigr)$. Then there exists a $1$-morphism $\bs h\colon \bX\ra\bZ$ with\/ $\bs h=\bigl(h,\bs h_{ik},\bs H_{ii'}^k,\bs H_i^{kk'}\bigr),$ such that\/ $h=g\ci f\colon X\ra Z,$ and for all\/ $i\in I,$ $j\in J,$ $k\in K$ we have $2$-morphisms over $h$
\e
\Th_{ijk}\colon \bs g_{jk}\ci\bs f_{ij}\Longra\bs h_{ik},
\label{kt4eq5}
\e
where as usual\/ \eq{kt4eq5} holds over\/ $S=\Im\chi_i\cap f^{-1}(\Im\psi_j)\cap h^{-1}(\Im\om_k),$ and for all\/ $i,i'\in I,$ $j,j'\in J,$ $k,k'\in K$ the following commute:
\begin{gather*}
\xymatrix@C=90pt@R=15pt{
*+[r]{(\bs g_{jk}\ci\bs f_{i'j})\ci\Tau_{ii'}} \ar@{=>}[d]^{\bs\al_{\bs g_{jk},\bs f_{i'j},\Tau_{ii'}}}
\ar@{=>}[rr]_(0.53){\Th_{i'jk}*\id_{\Tau_{ii'}}} && *+[l]{\bs h_{i'k}\ci\Tau_{ii'}} \ar@{=>}[d]_{\bs H_{ii'}^k}  
\\
*+[r]{\bs g_{jk}\ci(\bs f_{i'j}\ci\Tau_{ii'})} \ar@{=>}[r]^(0.65){\id_{\bs g_{jk}}*\bs F_{ii'}^j}
& \bs g_{jk}\ci\bs f_{ij} \ar@{=>}[r]^{\Th_{ijk}} & *+[l]{\bs h_{ik},\!} }
\allowdisplaybreaks\\
\xymatrix@C=90pt@R=15pt{
*+[r]{(\bs g_{j'k}\ci\Up_{jj'})\ci\bs f_{ij}} \ar@{=>}[d]^{\bs\al_{\bs g_{j'k},\Up_{jj'},\bs f_{ij}}}
\ar@{=>}[rr]_(0.53){\bs G_{jj'}^k*\id_{\bs f_{ij}}} && *+[l]{\bs g_{jk}\ci\bs f_{ij}} \ar@{=>}[d]_{\Th_{ijk}}  
\\
*+[r]{\bs g_{j'k}\ci(\Up_{jj'}\ci\bs f_{ij})} \ar@{=>}[r]^(0.65){\id_{\bs g_{j'k}}*\bs F_i^{jj'}}
& \bs g_{j'k}\ci\bs f_{ij'} \ar@{=>}[r]^{\Th_{ij'k}} & *+[l]{\bs h_{ik},\!} }
\allowdisplaybreaks\\
\xymatrix@C=90pt@R=15pt{
*+[r]{(\Phi_{kk'}\ci\bs g_{jk})\ci\bs f_{ij}} \ar@{=>}[d]^{\bs\al_{\Phi_{kk'},\bs g_{jk},\bs f_{ij}}}
\ar@{=>}[rr]_(0.53){\bs G_j^{kk'}*\id_{\bs f_{ij}}} && *+[l]{\bs g_{jk'}\ci\bs f_{ij}} \ar@{=>}[d]_{\Th_{ijk'}}  
\\
*+[r]{\Phi_{kk'}\ci(\bs g_{jk}\ci\bs f_{ij})} \ar@{=>}[r]^(0.65){\id_{\Phi_{kk'}}*\Th_{ijk}}
& \Phi_{kk'}\ci\bs h_{ik} \ar@{=>}[r]^{\bs H_i^{kk'}} & *+[l]{\bs h_{ik'}.\!} }
\end{gather*}

\noindent{\bf(b)} If\/ $\bs{\ti h}=\bigl(h,\bs{\ti h}_{ik},\bs{\ti H}{}_{ii'}^k,\bs{\ti H}{}_i^{kk'}\bigr),\ti\Th_{ijk}$ are alternative choices for $\bs h,\Th_{ijk}$ in {\bf(a)\rm,} then there is a unique $2$-morphism of Kuranishi spaces $\bs\eta=(\bs\eta_{ik})\colon \bs h\Ra\bs{\ti h}$ satisfying $\bs\eta_{ik}\od\Th_{ijk}=\ti\Th_{ijk}\colon \bs g_{jk}\ci\bs f_{ij}\Ra\bs{\ti h}_{ik}$ for all\/ $i\in I,$ $j\in J,$ $k\in K$.
\smallskip

\noindent{\bf(c)} If\/ $\bX=\bY$ and\/ $\bs f=\bs\id_\bY$ in {\bf(a)\rm,} so that\/ $I=J,$ then a possible choice for $\bs h,\Th_{ijk}$ in {\bf(a)} is $\bs h=\bs g$ and\/~$\Th_{ijk}=\bs G_{ij}^k$.

Similarly, if\/ $\bZ=\bY$ and\/ $\bs g=\bs\id_\bY$ in {\bf(a)\rm,} so that\/ $K=J,$ then a possible choice for $\bs h,\Th_{ijk}$ in {\bf(a)} is $\bs h=\bs f$ and\/~$\Th_{ijk}=\bs F_i^{jk}$.
\label{kt4prop1}
\end{prop}

Proposition \ref{kt4prop1}(a) gives possible values $\bs h$ for the composition $\bs g\ci\bs f\colon \bX\ra\bZ$. Since there is no distinguished choice, we choose $\bs g\ci\bs f$ arbitrarily.

\begin{dfn} For all pairs of 1-morphisms of Kuranishi spaces $\bs f\colon \bX\ra\bY$ and $\bs g\colon \bY\ra\bZ$, use the Axiom of Global Choice (a strong form of the Axiom of Choice working for classes as well as sets, see Shulman \cite[\S 7]{Shul} or Herrlick and Strecker \cite[\S 1.2]{HeSt}) to choose possible values of $\bs h\colon \bX\ra\bZ$ and $\Th_{ijk}$ in Proposition \ref{kt4prop1}(a), and write $\bs g\ci\bs f=\bs h$, and for $i\in I$, $j\in J$, $k\in K$
\begin{equation*}
\Th_{ijk}^{\bs g,\bs f}=\Th_{ijk}\colon \bs g_{jk}\ci\bs f_{ij}\Longra(\bs g\ci\bs f)_{ik}.
\end{equation*}
We call $\bs g\ci\bs f$ the {\it composition of\/ $1$-morphisms of Kuranishi spaces}.

For general $\bs f,\bs g$ we make these choices arbitrarily. However, if $\bX=\bY$ and $\bs f=\bs\id_\bY$ then we choose $\bs g\ci\bs\id_\bY=\bs g$ and $\Th_{jj'k}^{\bs g,\bs\id_\bY}=\bs G_{jj'}^k$, and if $\bZ=\bY$ and $\bs g=\bs\id_\bY$ then we choose $\bs\id_\bY\ci\bs f=\bs f$ and $\Th_{ijj'}^{\bs\id_\bY,\bs f}=\bs F_i^{jj'}$. This is allowed by Proposition~\ref{kt4prop1}(c).

The definition of a weak 2-category in \S\ref{ktA1} includes 2-isomorphisms $\bs\be_{\bs f}\colon \bs f\ci\bs\id_\bX\Ra\bs f$ and $\bs\ga_{\bs f}\colon \bs\id_\bY\ci\bs f\Ra\bs f$ in \eq{ktAeq7}, since one does not require $\bs f\ci\bs\id_\bX=\bs f$ and $\bs\id_\bY\ci\bs f=\bs f$ in a general weak 2-category. We define
\begin{equation*}
\bs\be_{\bs f}=\bs\id_{\bs f}\colon \bs f\ci\bs\id_\bX\Longra\bs f,\quad
\bs\ga_{\bs f}=\bs\id_{\bs f}\colon \bs\id_\bY\ci\bs f\Longra\bs f.
\end{equation*}

\label{kt4def4}
\end{dfn}

Since composition of 1-morphisms $\bs g\ci\bs f$ is natural only up to canonical 2-isomorphism, as in Proposition \ref{kt4prop1}(b), composition is associative only up to canonical 2-isomorphism. The next proposition comes from \cite[\S 4.3]{Joyc5}, \cite[\S 6.2]{Joyc6}.

\begin{prop} Let\/ $\bs e\colon \bW\ra\bX,$ $\bs f\colon \bX\ra\bY,$ $\bs g\colon \bY\ra\bZ$ be $1$-morphisms of Kuranishi spaces, and define composition of\/ $1$-morphisms as in Definition\/ {\rm\ref{kt4def4}}. Then using notation\/ {\rm\eq{kt4eq1}--\eq{kt4eq4},} there is a unique $2$-morphism
\begin{equation*}
\bs\al_{\bs g,\bs f,\bs e}\colon (\bs g\ci\bs f)\ci\bs e\Longra\bs g\ci(\bs f\ci\bs e)
\end{equation*}
with the property that for all\/ $h\in H,$ $i\in I,$ $j\in J$ and\/ $k\in K$ we have
\begin{equation*}
(\bs\al_{\bs g,\bs f,\bs e})_{hk}\od\Th^{\bs g\ci\bs f,\bs e}_{hik}\od (\Th_{ijk}^{\bs g,\bs f}*\id_{\bs e_{hi}})=
\Th_{hjk}^{\bs g,\bs f\ci\bs e}\od(\id_{\bs g_{jk}}*\Th_{hij}^{\bs f,\bs e})\od\bs\al_{\bs g_{jk},\bs f_{ij},\bs e_{hi}}.
\end{equation*}

\label{kt4prop2}
\end{prop}

We define vertical and horizontal composition of 2-morphisms:

\begin{dfn} Let $\bs f,\bs g,\bs h\colon \bX\ra\bY$ be 1-morphisms of Kuranishi spaces, using notation \eq{kt4eq2}--\eq{kt4eq3}, and $\bs\eta=(\bs\eta_{ij})\colon \bs f\Ra\bs g$, $\bs\ze=(\bs\ze_{ij})\colon \bs g\Ra\bs h$ be 2-mor\-phisms. Define a 2-morphism of Kuranishi spaces $\bs\ze\od\bs\eta\colon \bs f\Ra\bs h$ called the {\it vertical composition of\/ $2$-morphisms\/}  by
\begin{equation*}
\bs\ze\od\bs\eta=\bigl(\bs\ze_{ij}\od\bs\eta_{ij},\; i\in I,\; j\in J\bigr).
\end{equation*}

Next let $\bs e,\bs f\colon \bX\ra\bY$ and $\bs g,\bs h\colon \bY\ra\bZ$ be 1-morphisms of Kuranishi spaces, using notation \eq{kt4eq2}--\eq{kt4eq4}, and $\bs\eta=(\bs\eta_{ij})\colon \bs e\Ra\bs f$, $\bs\ze=(\bs\ze_{jk})\colon \bs g\Ra\bs h$ be 2-morphisms. As in \cite[\S 4.3]{Joyc5}, \cite[\S 6.2]{Joyc6}, using the stack property Theorem \ref{kt3thm2}(b) we can show that there is a unique 2-morphism of Kuranishi spaces $\bs\ze*\bs\eta\colon \bs g\ci\bs e\Ra\bs h\ci\bs f$ called the {\it horizontal composition of\/ $2$-morphisms}, such that for all $i\in I$, $j\in J$ and $k\in K$ we have
\begin{equation*}
(\bs\ze*\bs\eta)_{ik}\od\Th_{ijk}^{\bs g,\bs e}=
\Th_{ijk}^{\bs h,\bs f}\od(\bs\ze_{jk}*\bs\eta_{ij}).
\end{equation*}

\label{kt4def5}
\end{dfn}

We have now defined all the structures of a weak 2-category of Kuranishi spaces $\Kur$. As in \cite[\S 4.3]{Joyc5}, \cite[\S 6.2]{Joyc6}, the 2-category axioms are satisfied.

\begin{thm} The definitions and propositions above define the \begin{bfseries}weak\/ $2$-cat\-eg\-ory of Kuranishi spaces\end{bfseries}~$\Kur$.
\label{kt4thm1}
\end{thm}

\subsection{Manifolds, orbifolds, and m-Kuranishi spaces}
\label{kt43}

By imposing conditions on the Kuranishi neighbourhoods $(V_i,E_i,\Ga_i,s_i,\psi_i)$ on objects $\bX$ in $\Kur$, we can define interesting 2-subcategories of $\Kur$:
\begin{itemize}
\setlength{\itemsep}{0pt}
\setlength{\parsep}{0pt}
\item There is a full and faithful (2-)functor $F_\Man^\Kur\colon \Man\ra\Kur$ which identifies $\Man$ with the full 2-subcategory of objects $\bX=(X,\cK)$ in $\Kur$ in which $\cK$ has indexing set $I=\{0\}$, and one Kuranishi neighbourhood $(V_0,E_0,\Ga_0,s_0,\psi_0)$ with $E_0$ the zero vector bundle and~$\Ga_0=\{1\}$. 
\item The full 2-subcategory $\OrbKur$ of $\bX=(X,\cK)$ in $\Kur$ for which $E_i$ is the zero vector bundle for all Kuranishi neighbourhoods $(V_i,E_i,\Ga_i,s_i,\psi_i)$ in $\cK$ is equivalent as a 2-category to the 2-category of orbifolds~$\Orb$.
\item Write $\mKur\subset\Kur$ for the full 2-subcategory of $\bX\in\Kur$ for which $\Ga_i=\{1\}$ for all Kuranishi neighbourhoods $(V_i,E_i,\Ga_i,s_i,\psi_i)$ in $\cK$. We call such $\bX$ {\it m-Kuranishi spaces}, where `m-' is short for `manifold'. They are a kind of derived manifold, in the sense of \cite{Joyc2,Joyc3,Joyc4}, just as Kuranishi spaces are a kind of derived orbifold.
\item The full 2-subcategory of $\bX=(X,\cK)$ in $\Kur$ with only one Kuranishi neighbourhood $(V_i,E_i,\Ga_i,s_i,\psi_i)$ in $\cK$, so that $\md{I}=1$, is equivalent to the 2-category of global Kuranishi neighbourhoods $\GKN$ from Theorem~\ref{kt3thm1}.
\end{itemize}

We discuss the first three in more detail. Weak 2-functors are defined in~\S\ref{ktA2}. 

\begin{dfn} We will define a weak 2-functor $F_\Man^\Kur\colon \Man\ra\Kur$. If $X$ is a manifold, define a Kuranishi space $F_\Man^\Kur(X)=\bX=(X,\cK)$ with topological space $X$ and Kuranishi structure $\cK=\bigl(\{0\},(V_0,E_0,\Ga_0,s_0,\psi_0),\ab\Phi_{00},\ab\La_{000}\bigr)$, with indexing set $I=\{0\}$, one Kuranishi neighbourhood $(V_0,E_0,\Ga_0,s_0,\psi_0)$ with $V_0=X$, $E_0\ra V_0$ the zero vector bundle, $\Ga_0=\{1\}$, $s_0=0$, and $\psi_0=\id_X$, one coordinate change $\Phi_{00}=\id_{(V_0,E_0,\Ga_0,s_0,\psi_0)}$, and one 2-morphism~$\La_{000}=\id_{\Phi_{00}}$.

On (1-)morphisms, if $f\colon X\ra Y$ is a smooth map of manifolds and $\bX=F_\Man^\Kur(X)$, $\bY=F_\Man^\Kur(Y)$, define a 1-morphism $F_\Man^\Kur(f)=\bs f\colon \bX\ra\bY$ by $\bs f=(f,\bs f_{00},\bs F_{00}^0,\bs F_0^{00})$, where $\bs f_{00}=(P_{00},\pi_{00},f_{00},\hat f_{00})$ with $P_{00}=X$, $\pi_{00}=\id_X$, $f_{00}=f$, and $\hat f_{00}$ is the zero map on zero vector bundles, and~$\bs F_{00}^0=\bs F_0^{00}=\id_{\bs f_{00}}$.

On 2-morphisms, regarding $\Man$ as a 2-category, the only 2-morphisms are identity morphisms $\id_f\colon f\Ra f$ for (1-)morphisms $f\colon X\ra Y$ in $\Man$. We define~$F_\Man^\Kur(\id_f)=\bs\id_{F_\Man^\Kur(f)}$.

If $f\colon X\ra Y$, $g\colon Y\ra Z$ are (1-)morphisms in $\Man$, there is a unique 2-morphism in $\Kur$
\begin{equation*}
(F_\Man^\Kur)_{g,f}\colon F_\Man^\Kur(g)\ci F_\Man^\Kur(f)\Longra F_\Man^\Kur(g\ci f).
\end{equation*}
For any object $X$ in $\Man$, define
\begin{equation*}
(F_\Man^\Kur)_X:=\bs\id_{\bs\id_\bX}\colon F_\Man^\Kur(\id_X)=\bs\id_\bX\Longra\bs\id_\bX=\bs\id_{F_\Man^\Kur(X)}.
\end{equation*}
As in \cite[\S 4.3]{Joyc5}, this defines a full and faithful weak 2-functor $F_\Man^\Kur\colon \Man\ra\Kur$. Thus we can regard manifolds $X$ as special examples of Kuranishi spaces. 

We say that a Kuranishi space $\bX$ {\it is a manifold\/} if $\bX\simeq F_\Man^\Kur(X')$ in $\Kur$, for some manifold~$X'$.
\label{kt4def6}
\end{dfn}

The next remark reviews different definitions of orbifolds in the literature.

\begin{rem} Orbifolds are generalizations of manifolds locally modelled on $\R^n/G$, for $G$ a finite group acting linearly on $\R^n$. They were introduced by Satake \cite{Sata}, who called them `V-manifolds'. Later they were studied by Thurston \cite[Ch.~13]{Thur} who gave them the name `orbifold'.

As for Kuranishi spaces, defining orbifolds $\fX,\fY$ and smooth maps $\ff\colon \fX\ra\fY$ was initially problematic, and early definitions of ordinary categories of orbifolds \cite{Sata,Thur} had some bad differential-geometric behaviour  (e.g. for some definitions, one cannot define pullbacks $\ff^*(\mathfrak E)$ of orbifold vector bundles $\mathfrak E\ra\fY$). It is now generally agreed that it is best to define orbifolds to be a 2-category. See Lerman \cite{Lerm} for a good overview of ways to define orbifolds.

There are three main definitions of ordinary categories of orbifolds:
\begin{itemize}
\setlength{\itemsep}{0pt}
\setlength{\parsep}{0pt}
\item[(a)] Satake \cite{Sata} and Thurston \cite{Thur} defined an orbifold $\fX$ to be a Hausdorff topological space $X$ with an atlas $\bigl\{(V_i,\Ga_i,\psi_i)\colon i\in I\bigr\}$ of orbifold charts $(V_i,\Ga_i,\psi_i)$, where $V_i$ is a manifold, $\Ga_i$ a finite group acting smoothly (and locally effectively) on $V_i$, and $\psi_i\colon V_i/\Ga_i\ra X$ a homeomorphism with an open set in $X$, and pairs of charts $(V_i,\Ga_i,\phi_i),(V_j,\Ga_j,\phi_j)$ satisfy compatibility conditions on their overlaps in $X$. Smooth maps $\ff\colon \fX\ra\fY$ between orbifolds are continuous maps $f\colon X\ra Y$ of the underlying spaces, which lift locally to smooth maps on the charts, giving a category $\Orb_{\rm ST}$.
\item[(b)] Chen and Ruan \cite[\S 4]{ChRu} defined orbifolds $\fX$ in a similar way to \cite{Sata,Thur}, but using germs of orbifold charts $(V_p,\Ga_p,\psi_p)$ for $p\in X$. Their morphisms $\ff\colon \fX\ra\fY$ are called {\it good maps}, giving a category $\Orb_{\rm CR}$. 
\item[(c)] Moerdijk and Pronk \cite{Moer,MoPr} defined a category of orbifolds $\Orb_{\rm MP}$ as {\it proper \'etale Lie groupoids\/} in $\Man$. Their definition of smooth map $\ff\colon \fX\ra\fY$, called {\it strong maps\/} \cite[\S 5]{MoPr}, is complicated: it is an equivalence class of diagrams $\smash{\fX\,{\buildrel\phi\over \longleftarrow}\,\fX'\,{\buildrel\psi\over\longra}\,\fY}$, where $\fX'$ is a third orbifold, and $\phi,\psi$ are morphisms of groupoids with $\phi$ an equivalence (loosely, a diffeomorphism).
\end{itemize}
A book on orbifolds in the sense of \cite{ChRu,Moer,MoPr} is Adem, Leida and Ruan~\cite{ALR}.

There are four main definitions of 2-categories of orbifolds: 
\begin{itemize}
\setlength{\itemsep}{0pt}
\setlength{\parsep}{0pt}
\item[(i)] Pronk \cite{Pron} defines a strict 2-category $\bf LieGpd$ of Lie groupoids in $\Man$ as in (c), with the obvious 1-morphisms of groupoids, and localizes by a class of weak equivalences $\cW$ to get a weak 2-category $\Orb_{\rm Pr}={\bf LieGpd}[\cW^{-1}]$.
\item[(ii)] Lerman \cite[\S 3.3]{Lerm} defines a weak 2-category $\Orb_{\rm Le}$ of Lie groupoids in $\Man$ as in (c), with a non-obvious notion of 1-morphism called `Hilsum--Skandalis morphisms' involving `bibundles', and does not need to localize.

Henriques and Metzler \cite{HeMe} also use Hilsum--Skandalis morphisms.
\item[(iii)] Behrend and Xu \cite[\S 2]{BeXu}, Lerman \cite[\S 4]{Lerm} and Metzler \cite[\S 3.5]{Metz} define a strict 2-category of orbifolds $\Orb_{\rm ManSta}$ as a class of Deligne--Mumford stacks on the site $(\Man,{\mathcal J}_\Man)$ of manifolds with Grothendieck topology ${\mathcal J}_\Man$ coming from open covers.
\item[(iv)] The author \cite{Joyc1} defines a strict 2-category of orbifolds $\Orb_{C^\iy{\rm Sta}}$ as a class of Deligne--Mumford stacks on the site $(\CSch,{\mathcal J}_\CSch)$ of $C^\iy$-schemes.
\end{itemize}

As in Behrend and Xu \cite[\S 2.6]{BeXu}, Lerman \cite{Lerm}, Pronk \cite{Pron}, and the author \cite[Th.~9.26]{Joyc1}, approaches (i)--(iv) give equivalent weak 2-categories $\Orb_{\rm Pr},\ab\Orb_{\rm Le},\ab\Orb_{\rm ManSta},\ab\Orb_{C^\iy{\rm Sta}}$. As they are equivalent, the differences between them are not of mathematical importance, but more a matter of convenience or taste. Properties of localization also imply that $\Orb_{\rm MP}\simeq \Ho(\Orb_{\rm Pr})$. Thus, all of (c) and (i)--(iv) are equivalent at the level of homotopy categories.
\label{kt4rem1}
\end{rem}

We now give a fifth definition of a weak 2-category of orbifolds.

\begin{dfn} Define the weak 2-category $\OrbKur\subset\Kur$ of {\it Kuranishi orbifolds}, or just {\it orbifolds}, to be the full 2-subcategory of $\bX=(X,\cK)$ in $\Kur$ for which $E_i$ is the zero vector bundle for all Kuranishi neighbourhoods $(V_i,E_i,\Ga_i,s_i,\psi_i)$ in~$\cK$.

This allows us to simplify the notation a lot. Equations in \S\ref{kt31} involving error terms $O\bigl(\pi_{ij}^*(s_i)\bigr)$ or $O\bigl(\pi_{ij}^*(s_i)^2\bigr)$ become exact, as $s_i=0$. As $E_i,s_i$ are zero we can take `orbifold charts' to be $(V_i,\Ga_i,\psi_i)$. As $\hat\phi_{ij}=0$ we can take coordinate changes to be $\Phi_{ij}=(P_{ij},\pi_{ij},\phi_{ij})$, and we can also take $V_{ij}=\pi_{ij}(P_{ij})$ to be equal to $\bar\psi_i^{-1}(S)$, rather than just an open neighbourhood of $\bar\psi_i^{-1}(S)$ in $V_i$, since $\bar\psi_i^{-1}(S)$ is open in $V_i$ when $s_i=0$. For 2-morphisms $\La_{ij}=[\dot P_{ij},\la_{ij},\hat\la_{ij}]\colon \Phi_{ij}\Ra\Phi_{ij}'$ in \S\ref{kt31}, we have $\hat\la_{ij}=0$, and we are forced to take $\dot P_{ij}=P_{ij}$, and the equivalence relation $\approx$ in Definition \ref{kt3def4} becomes trivial, so we can take 2-morphisms to be just~$\la_{ij}$.

We say that a Kuranishi space $\bX$ {\it is an orbifold\/} if $\bX\simeq \bX'$ in $\Kur$ for some~$\bX'\in\OrbKur\subset\Kur$.
\label{kt4def7}
\end{dfn}

The next theorem is proved in \cite[\S 4.5]{Joyc5}, \cite[\S 6.6]{Joyc6}, by giving a full and faithful embedding $\OrbKur\hookra\Orb_{\rm Le}$, and showing it is surjective on equivalence classes.

\begin{thm} The weak\/ $2$-category of Kuranishi orbifolds $\OrbKur$ above is equivalent to the $2$-categories of orbifolds $\Orb_{\rm Pr},\ab\Orb_{\rm Le},\ab\Orb_{\rm ManSta},\ab\Orb_{C^\iy{\rm Sta}}$ in {\rm\cite{BeXu,Joyc1,Lerm,Metz,Pron}} described in Remark\/ {\rm\ref{kt4rem1}}. Also there is an equivalence of categories $\Ho(\OrbKur)\simeq\Orb_{\rm MP},$ for $\Orb_{\rm MP}$ the category of orbifolds from  Moerdijk and Pronk\/ {\rm\cite{Moer,MoPr}}.
\label{kt4thm2}
\end{thm}

Fukaya et al.\ \cite[\S 9]{FOOO4} and McDuff \cite{McDu} also define (effective) orbifolds as special examples of FOOO Kuranishi spaces and MW Kuranishi atlases.

\begin{dfn} As in \cite[\S 4.7]{Joyc5}, \cite[Ch.~4 \& \S 6.2]{Joyc6}, define the weak 2-category $\mKur\subset\Kur$ of {\it m-Kuranishi spaces\/} to be the full 2-subcategory of $\bX=(X,\cK)$ in $\Kur$ for which $\Ga_i=\{1\}$ for all Kuranishi neighbourhoods $(V_i,E_i,\Ga_i,s_i,\psi_i)$ in~$\cK$.

As in Definition \ref{kt4def7}, this allows us to simplify the notation a lot. We take `m-Kuranishi neighbourhoods' to be $(V_i,E_i,s_i,\psi_i)$ rather than $(V_i,E_i,\Ga_i,s_i,\psi_i)$. In coordinate changes, $\pi_{ij}\colon P_{ij}\ra V_{ij}$ in \S\ref{kt31} is a diffeomorphism, since it is a principal $\Ga_j$-bundle for $\Ga_j=\{1\}$, so we can replace $P_{ij}$ by $V_{ij}$ and take coordinate changes to be $(V_{ij},\phi_{ij},\hat\phi_{ij})$, and 2-morphisms to be $[\dot V_{ij},\hat\la_{ij}]$ for open $\dot V_{ij}\subseteq V_{ij}\subseteq V_i$. The m-Kuranishi space analogues $\mKN_S(X),\GmKN$ of $\KN_S(X),\GKN$ in \S\ref{kt31} are then {\it strict\/} 2-categories, as the canonical identification $\la_{ijkl}$ in \eq{kt3eq3} is replaced by the identity map.
\label{kt4def8}
\end{dfn}

\subsection{Relation to FOOO, MW, DY, polyfolds, and d-orbifolds}
\label{kt44}

The next four theorems are proved in \cite[\S 4.8]{Joyc5}, \cite[Ch.~7]{Joyc6}. Theorem \ref{kt4thm6} is a corollary of Theorem \ref{kt4thm5} and Theorem \ref{kt2thm5} from Yang~\cite[Th.~3.1.7]{Yang1}.

\begin{thm} Suppose $\bX=(X,\cK)$ is a FOOO Kuranishi space without boundary, as in Definition\/ {\rm\ref{kt2def3}}. Then we can construct a Kuranishi space $\bX'=(X,\cK')$ in the sense of\/ {\rm\S\ref{kt41}} with\/ $\vdim\bX'=\vdim\bX,$ with the same compact topological space $X,$ and\/ $\bX'$ is natural up to equivalence in the $2$-category~$\Kur$.

\label{kt4thm3}
\end{thm}

\begin{thm} Suppose $X$ is a compact, metrizable topological space with an MW weak Kuranishi atlas\/ $\cK$ without boundary, of virtual dimension $n\in\Z,$ in the sense of Definition\/ {\rm\ref{kt2def8}}. Then we can make $X$ into a Kuranishi space $\bX'=(X,\cK')$ in the sense of\/ {\rm\S\ref{kt41}} with\/ $\vdim\bX'=n,$ and\/ $\bX'$ is natural up to equivalence in the $2$-category $\Kur$. Commensurate MW weak Kuranishi atlases $\cK,\ti\cK$ on $\bX$ yield equivalent Kuranishi spaces $\bX{}',\ti\bX{}'$.
\label{kt4thm4}
\end{thm}

\begin{thm} Suppose $X$ is a compact, metrizable topological space with a DY Kuranishi structure $\cK$ without boundary in the sense of Definition\/ {\rm\ref{kt2def9}}. Then we can construct a Kuranishi space $\bX'=(X,\cK')$ in the sense of\/ {\rm\S\ref{kt41}} with\/ $\vdim\bX'=\vdim\bX,$ with the same topological space $X,$ and\/ $\bX'$ is natural up to equivalence in the $2$-category $\Kur$. R-equivalent DY Kuranishi structures $\cK,\ti\cK$ on $\bX$ yield equivalent Kuranishi spaces\/~$\bX{}',\ti\bX{}'$.
\label{kt4thm5}
\end{thm}

\begin{thm} Suppose we are given a `polyfold Fredholm structure'\/ $\cP$ on a compact metrizable topological space $X$ in the sense of\/ {\rm\cite{Hofe,HWZ1,HWZ2,HWZ3,HWZ4,HWZ5,HWZ6,HWZ7},} that is, we write $X$ as the zeroes of an sc-Fredholm section ${\mathfrak s}\colon \fV\ra\fE$ of a strong polyfold vector bundle $\fE\ra\fV$ over a polyfold $\fV,$ where $\mathfrak s$ has constant Fredholm index $n\in\Z$. Then we can make $X$ into a Kuranishi space $\bX=(X,\cK)$ in the sense of\/ {\rm\S\ref{kt41}} with\/ $\vdim\bX=n,$ and\/ $\bX$ is natural up to equivalence in the $2$-category\/~$\Kur$.
\label{kt4thm6}
\end{thm}

The constructions of Theorems \ref{kt4thm3}--\ref{kt4thm6} are compatible with other geometric structures. In particular, orientations on FOOO Kuranishi spaces or MW weak Kuranishi atlases or polyfold Fredholm structures correspond canonically to orientations on the corresponding Kuranishi spaces $\bX'$. Also smooth maps $\bs f\colon \bX\ra Y$ from a FOOO Kuranishi space $\bX$ to a manifold $Y$ correspond to 1-morphisms $\bs f'\colon \bX'\ra Y$ in $\Kur$, uniquely up to 2-isomorphism. Theorems \ref{kt4thm3}--\ref{kt4thm6} have applications to symplectic geometry.

The next theorem will be proved in \cite[Ch.~8]{Joyc6}.

\begin{thm} There are equivalences of weak\/ $2$-categories
\begin{equation*}
\smash{\mKur\simeq\dMan \qquad\text{and\/}\qquad \Kur\simeq\dOrb,}
\end{equation*}
where\/ $\dMan,\dOrb$ are the strict\/ $2$-categories of `d-manifolds' and `d-orbifolds', as in the author {\rm\cite{Joyc2,Joyc3,Joyc4}} and\/~{\rm\S\ref{kt36}}.
\label{kt4thm9}
\end{thm}

Thus, m-Kuranishi spaces and d-manifolds, and Kuranishi spaces and d-orbifolds, are for most purposes interchangeable.

To understand the relationship between Kuranishi spaces and d-orbifolds, consider the two equivalent ways of defining ordinary manifolds:
\begin{itemize}
\setlength{\itemsep}{0pt}
\setlength{\parsep}{0pt}
\item[(a)] As a Hausdorff, second countable topological space $X$ equipped with an atlas of charts $\bigl\{(V_i,\phi_i)\colon i\in I\bigr\}$; or
\item[(b)] As a Hausdorff, second countable topological space $X$ equipped with a sheaf of $\R$-algebras $\O_X^{C^\iy}$ which is embedded as a subsheaf $\O_X^{C^\iy}\subset\O_X^{C^0}$ of the sheaf $\O_X^{C^0}$ of continuous functions $f\colon X\ra\R$, such that $(X,\O_X^{C^\iy})$ is locally modelled on $(\R^n,\O_{\R^n}^{C^\iy})$, where $\O_{\R^n}^{C^\iy}$ is the sheaf of smooth functions~$f\colon \R^n\ra\R$.
\end{itemize}
If we try to define `derived manifolds' and `derived orbifolds' by generalizing definition (a), we end up with (m-)Kuranishi spaces (or something similar); if we try to do it by generalizing definition (b), we end up with d-manifolds and d-orbifolds (or something similar).

Kuranishi spaces are more convenient than d-orbifolds if we wish to define `derived' versions of categories of `manifolds' which are not ordinary smooth manifolds, for instance, manifolds with corners of various kinds, or some classes of singular manifolds. This is because in the Kuranishi neighbourhoods $(V_i,\ab E_i,\ab\Ga_i,\ab s_i,\ab\psi_i)$ on a Kuranishi space $\bX$, it is easy to take $V_i$ to be some different kind of `manifold', but smooth manifolds are built into the theory of d-manifolds and d-orbifolds at a deep level, in the notion of $C^\iy$-ring.

\section{Differential geometry of Kuranishi spaces}
\label{kt5}

\subsection{Isotropy groups, and tangent and obstruction spaces}
\label{kt51}

The next definition is taken from \cite[\S 4.6]{Joyc5}, \cite[\S 6.5]{Joyc6}.

\begin{dfn} Let $\bX=(X,\cK)$ be a Kuranishi space, with $\cK=\bigl(I,(V_i,\ab E_i,\ab\Ga_i,\ab s_i,\ab\psi_i)_{i\in I}$, $\Phi_{ij,\;i,j\in I}$, $\La_{ijk,\; i,j,k\in I}\bigr)$, and let $x\in\bX$. Choose an arbitrary $i\in I$ with $x\in\Im\psi_i$, and choose $v_i\in s_i^{-1}(0)\subseteq V_i$ with~$\bar\psi_i(v_i)=x$. 

Define a finite group $G_x\bX$ called the {\it isotropy group of $\bX$ at\/} $x$ by
\e
G_x\bX=\bigl\{\ga\in\Ga_i\colon \ga\cdot v_i=v_i\bigr\}=\Stab_{\Ga_i}(v_i),
\label{kt5eq1}
\e
as a subgroup of $\Ga_i$. Define finite-dimensional real vector spaces $T_x\bX$, the {\it tangent space of\/ $\bX$ at\/} $x$, and $O_x\bX$, the {\it obstruction space of\/ $\bX$ at\/} $x$, to be the kernel and cokernel of $\d s_i\vert_{v_i}$, so that they fit into an exact sequence
\begin{equation*}
\xymatrix@C=23pt{ 0 \ar[r] & T_x\bX \ar[r] & T_{v_i}V_i \ar[rr]^{\d s_i\vert_{v_i}} && E_i\vert_{v_i} \ar[r] & O_x\bX \ar[r] & 0. }
\end{equation*}
The actions of $\Ga_i$ on $V_i,E_i$ induce linear actions of $G_x\bX$ on $T_x\bX,O_x\bX$, by the commutative diagram for each $\ga\in G_x\bX$:
\begin{equation*}
\xymatrix@C=26pt@R=15pt{ 0 \ar[r] & T_x\bX \ar[r] \ar@{.>}[d]^{\ga\cdot} & T_{v_i}V_i \ar[rr]_{\d s_i\vert_{v_i}} \ar[d]^{\d(\ga\cdot)} && E_i\vert_{v_i} \ar[d]^{\ga\cdot} \ar[r] & O_x\bX \ar@{.>}[d]^{\ga\cdot} \ar[r] & 0 \\
0 \ar[r] & T_x\bX \ar[r] & T_{v_i}V_i \ar[rr]^{\d s_i\vert_{v_i}} && E_i\vert_{v_i} \ar[r] & O_x\bX \ar[r] & 0.\!\! }
\end{equation*}
This makes $T_x\bX,O_x\bX$ into representations of $G_x\bX$. Definition \ref{kt4def1}(b) yields
\begin{equation*}
\dim T_x\bX-\dim O_x\bX=\vdim\bX.
\end{equation*}

The triple $(G_x\bX,T_x\bX,O_x\bX)$ depends on the choice of $i,v_i$. In \cite[\S 4.6]{Joyc5} and \cite[\S 6.5]{Joyc6} we show two choices $i,v_i$ and $i',v'_{i'}$ yield triples $(G_x\bX,T_x\bX,O_x\bX)$, $(G_x\bX',\ab T_x\bX',\ab O_x\bX')$ which have a finite, nonempty class of natural isomorphisms. Two such isomorphisms $(G_x\bX,T_x\bX,O_x\bX)\ra(G_x'\bX,T_x'\bX,O_x'\bX)$ differ by conjugation by an element of $G_x'\bX$, and behave as expected under composition.

\label{kt5def1}
\end{dfn}

We discuss functoriality of the $G_x\bX,T_x\bX,O_x\bX$ under 1- and 2-morphisms:

\begin{dfn} Let $\bs f\colon \bX\ra\bY$ be a 1-morphism of Kuranishi spaces, with notation \eq{kt4eq2}--\eq{kt4eq3}, and let $x\in\bX$ with $\bs f(x)=y\in\bY$. Then Definition \ref{kt5def1} gives $G_x\bX,T_x\bX,O_x\bX$, defined using $i\in I$ and $u_i\in U_i$ with $\bar\chi_i(u_i)=x$, and $G_y\bY,T_y\bY,O_y\bY$, defined using $j\in J$ and $v_j\in V_j$ with $\bar\chi_j(v_j)=y$. In $\bs f$ we have a 1-morphism $\bs f_{ij}=(P_{ij},\pi_{ij},f_{ij},\hat f_{ij})$ over $f$. Define
\begin{equation*}
S_{x,\bs f}=\bigl\{p\in P_{ij}\colon \pi_{ij}(p)=u_i,\;\> f_{ij}(p)=v_j\bigr\}.
\end{equation*}
Then $S_{x,\bs f}$ is invariant under the commuting actions of $G_x\bX=\Stab_{\Be_i}(u_i)\subseteq\Be_i$ and $G_y\bY=\Stab_{\Ga_j}(v_j)\subseteq\Ga_j$ on $P_{ij}$ induced by the $\Be_i,\Ga_j$-actions on $P_{ij}$, where $G_y\bY$ acts freely and transitively on $S_{x,\bs f}$, but the action of $G_x\bX$ need not be free or transitive. Choose an arbitrary point $p_0\in S_{x,\bs f}$. Our definitions of $G_x\bs f,T_x\bs f,O_x\bs f$ will depend on this choice. 

Define morphisms $G_x\bs f\colon G_x\bX\ra G_y\bY$, $T_x\bs f\colon T_x\bX\ra T_y\bY$, $O_x\bs f\colon O_x\bX\ra O_y\bY$ by $G_x\bs f(\ga)=\ga'$ if $\ga\cdot p_0=(\ga')^{-1}\cdot p_0$ in $S_{x,\bs f}$, using the actions of $G_x\bX,G_y\bY$ on $S_{x,\bs f}$ with $G_y\bY$ free and transitive, and the commutative diagram
\begin{equation*}
\xymatrix@C=26pt@R=19pt{ 0 \ar[r] & T_x\bX \ar[r] \ar@{.>}[d]^{T_x\bs f} & T_{u_i}U_i \ar[rr]_{\d r_i\vert_{u_i}} \ar[d]^(0.49){\d f_{ij}\vert_{p_0}\ci(\d\pi_{ij}\vert_{p_0})^{-1}} && D_i\vert_{u_i} \ar[d]^{\hat f_{ij}\vert_{p_0}} \ar[r] & O_x\bX \ar@{.>}[d]^{O_x\bs f} \ar[r] & 0 \\
0 \ar[r] & T_y\bY \ar[r] & T_{v_j}V_j \ar[rr]^{\d s_j\vert_{v_j}} && E_j\vert_{v_j} \ar[r] & O_y\bY \ar[r] & 0.\!\! }
\end{equation*}
Then $T_x\bs f,O_x\bs f$ are $G_x\bs f$-equivariant linear maps.

If $p_0'\in S_{x,\bs f}$ is an alternative choice for $p_0$, yielding $G'_x\bs f,T'_x\bs f,O'_x\bs f$, there is a unique $\de\in G_y\bY$ with $\de\cdot p_0=p_0'$, and then $G'_x\bs f(\ga)=\de (G_x\bs f(\ga))\de^{-1}$, $T'_x\bs f(v)=\de\cdot T_x\bs f(v)$, $O'_x\bs f(w)=\de\cdot O_x\bs f(w)$ for all $\ga\in G_x\bX$, $v\in T_x\bX$, and $w\in O_x\bX$. That is, the triple $(G_x\bs f,T_x\bs f,O_x\bs f)$ is canonical up to conjugation by an element of~$G_y\bY$.

Continuing with the same notation, suppose $\bs g\colon \bX\ra\bY$ is another 1-mor\-ph\-ism and $\bs\eta\colon \bs f\Ra\bs g$ a 2-morphism in $\Kur$. Then above we define $G_x\bs g,T_x\bs g,O_x\bs g$ by choosing an arbitrary point $q_0\in S_{x,\bs g}$, where
\begin{equation*}
S_{x,\bs g}=\bigl\{q\in Q_{ij}\colon \pi_{ij}(q)=u_i,\;\> g_{ij}(q)=v_j\bigr\},
\end{equation*}
with $\bs g_{ij}=(Q_{ij},\pi_{ij},g_{ij},\hat g_{ij})$ in $\bs g$. In $\bs\eta$ we have $\bs\eta_{ij}=[\dot P_{ij},\la_{ij},\hat\la_{ij}]$ represented by $(\dot P_{ij},\la_{ij},\hat\la_{ij})$, where $\dot P_{ij}\subseteq P_{ij}$ and $\la_{ij}\colon \dot P_{ij}\ra Q_{ij}$. From the definitions we find that $S_{x,\bs f}\subseteq\dot P_{ij}$, and $\la_{ij}\vert_{S_{x,\bs f}}\colon S_{x,\bs f}\ra S_{x,\bs g}$ is a bijection. Since $G_y\bY$ acts freely and transitively on $S_{x,\bs g}$, there is a unique element $G_x\bs\eta\in G_y\bY$ with $G_x\bs\eta\cdot \la_{ij}(p_0)=q_0$. Then we have
\begin{align*}
G_x\bs g(\ga)&=(G_x\bs\eta)(G_x\bs f(\ga))(G_x\bs\eta)^{-1},\quad T_x\bs g(v)=G_x\bs\eta\cdot T_x\bs f(v),\quad\text{and} \\ 
O_x\bs g(w)&=G_x\bs\eta\cdot O_x\bs f(w)\quad\text{for all $\ga\in G_x\bX$, $v\in T_x\bX$, and $w\in O_x\bX$.} 
\end{align*}
That is, $(G_x\bs g,T_x\bs g,O_x\bs g)$ is conjugate to $(G_x\bs f,T_x\bs f,O_x\bs f)$ under $G_x\bs\eta\in G_y\bY$, the same indeterminacy as in the definition of~$(G_x\bs f,T_x\bs f,O_x\bs f)$.

Suppose instead that $\bs g\colon \bY\ra\bZ$ is another 1-morphism of Kuranishi spaces and $\bs g(y)=z\in\bZ$. Then in a similar way we can show there is a canonical element $G_{x,\bs g,\bs f}\in G_z\bZ$ such that for all $\ga\in G_x\bX$, $v\in T_x\bX$, $w\in O_x\bX$ we have
\begin{align*}
G_x(\bs g\ci\bs f)(\ga)&=(G_{x,\bs g,\bs f})((G_y\bs g\ci G_x\bs f)(\ga))(G_{x,\bs g,\bs f})^{-1},\\ 
T_x(\bs g\ci\bs f)(v)&=G_{x,\bs g,\bs f}\cdot (T_y\bs g\ci T_x\bs f)(v), \\ 
O_x(\bs g\ci\bs f)(w)&=G_{x,\bs g,\bs f}\cdot (O_y\bs g\ci O_x\bs f)(w). 
\end{align*}
That is, $(G_x(\bs g\ci\bs f),T_x(\bs g\ci\bs f),O_x(\bs g\ci\bs f))$ is conjugate to $(G_y\bs g,T_y\bs g,O_y\bs g)\ci (G_x\bs f,\ab T_x\bs f,\ab O_x\bs f)$ under~$G_{x,\bs g,\bs f}\in G_z\bZ$.

Since 2-morphisms $\bs\eta\colon \bs f\Ra\bs g$ relate triples $(G_x\bs f,T_x\bs f,O_x\bs f)$ and $(G_x\bs g,\ab T_x\bs g,\ab O_x\bs g)$ by isomorphisms, if $\bs f\colon \bX\ra\bY$ is an equivalence in $\Kur$ then $G_x\bs f,T_x\bs f,O_x\bs f$ are isomorphisms for all~$x\in\bX$.
\label{kt5def2}
\end{dfn}

We could use the Axiom of (Global) Choice to choose particular values of $G_x\bX,T_x\bX,O_x\bX,G_x\bs f,T_x\bs f,O_x\bs f,G_x\bs\eta$ in Definitions \ref{kt5def1}--\ref{kt5def2} for all $\bX,x,\bs f,\bs\eta$, but there seems no need for this. All the definitions we give using $G_x\bX,\ldots,G_x\bs\eta$ will be independent of the arbitrary choices involved.

The next theorem is proved in \cite[\S 10.4 \& \S 10.5]{Joyc6}. 

\begin{thm}{\bf(a)} A Kuranishi space $\bX$ is a manifold, in the sense of Definition\/ {\rm\ref{kt4def6},} if and only if\/ $G_x\bX=\{1\}$ and\/ $O_x\bX=0$ for all\/~$x\in\bX$.
\smallskip

\noindent{\bf(b)} A Kuranishi space $\bX$ is an orbifold, in the sense of Definition\/ {\rm\ref{kt4def7},} if and only if\/ $O_x\bX=0$ for all\/~$x\in\bX$. 
\smallskip

\noindent{\bf(c)} A $1$-morphism $\bs f\colon \bX\ra\bY$ in $\Kur$ is an equivalence in $\Kur$ if and only if\/ $G_x\bs f\colon G_x\bX\ra G_{\bs f(x)}\bY,$ $T_x\bs f\colon T_x\bX\ra T_{\bs f(x)}\bY$ and\/ $O_x\bs f\colon O_x\bX\ra O_{\bs f(x)}\bY$ are isomorphisms for all\/ $x\in\bX,$ and\/ $f\colon X\ra Y$ is a bijection.
\label{kt5thm1}
\end{thm}

\begin{dfn} Write $\KurtrG\subset\Kur$ for the full 2-subcategory of Kuranishi spaces $\bX$ with $G_x\bX=\{1\}$ for all $x\in\bX$. Note that the 2-subcategory $\mKur\subset\Kur$ in \S\ref{kt43} has $\mKur\subset\KurtrG$, since if $\bX=(X,\cK)\in\mKur$ then the Kuranishi neighbourhoods $(V_i,E_i,\Ga_i,s_i,\psi_i)$ in $\cK$ have $\Ga_i=\{1\}$, which implies that $G_x\bX=\{1\}$ as $G_x\bX\subseteq\Ga_i$ by~\eq{kt5eq1}.
\label{kt5def3}
\end{dfn}

The next theorem is proved in \cite[\S 4.7]{Joyc5} and~\cite[\S 6.5]{Joyc6}.

\begin{thm} The inclusion $\mKur\hookra\KurtrG$ is an equivalence of weak\/ $2$-categories.

\label{kt5thm2}
\end{thm}

Combined with Theorem \ref{kt4thm9} this gives an equivalence $\KurtrG\simeq\dMan$ of weak 2-categories. That is, $\dMan,\mKur,\KurtrG$ are 2-categories of `derived manifolds', and the theorems say that they are all essentially the same. The moral is that a derived orbifold (Kuranishi space) $\bX$ is a derived manifold if and only if the isotropy groups $G_x\bX$ for $x\in\bX$ are all trivial.

\subsection{W-transverse morphisms and fibre products}
\label{kt52}

Smooth maps of manifolds $g\colon X\ra Z$ and $h\colon Y\ra Z$ are called {\it transverse\/} if $T_xg\op T_yh\colon T_xX\op T_yY\ra T_zZ$ is surjective for all $x\in X$ and $y\in Y$ with $g(x)=h(y)=z\in Z$. It is well known that if $g,h$ are transverse then the fibre product $W=X\t_{g,Z,h}Y$ exists in the category $\Man$, in the sense of category theory, with $\dim W=\dim X+\dim Y-\dim Z$.

We now explain the analogue of this for Kuranishi spaces. We define two notions of transverse 1-morphisms in $\Kur$, a weak and a strong:

\begin{dfn} Let $\bs g\colon \bX\ra\bZ$ and $\bs h\colon \bY\ra\bZ$ be 1-morphisms of Kuranishi spaces. Call $\bs g,\bs h$ {\it w-transverse\/} if for all $x\in\bX$, $y\in\bY$ with $\bs g(x)=\bs h(y)=z$ in $\bZ$, and all $\ga\in G_z\bZ$, then $O_x\bs g\op (\ga\cdot O_y\bs h)\colon O_x\bX\op O_y\bY\ra O_z\bZ$ is surjective.

We call $\bs g,\bs h$ {\it transverse\/} if for all $x\in\bX$, $y\in\bY$ with $\bs g(x)=\bs h(y)=z\in\bZ$, and all $\ga\in G_z\bZ$, then $T_x\bs g\op(\ga\cdot T_y\bs h)\colon T_x\bX\op T_y\bY\ra T_z\bZ$ is surjective, and $O_x\bs g\op(\ga\cdot O_y\bs h)\colon O_x\bX\op O_y\bY\ra O_z\bZ$ is an isomorphism.
\label{kt5def4}
\end{dfn}

The next theorem is proved in~\cite[Ch.~11]{Joyc6}.

\begin{thm} Suppose $\bs g\colon \bX\ra\bZ,$ $\bs h\colon \bY\ra\bZ$ are w-transverse $1$-morphisms of Kuranishi spaces. Then the fibre product\/ $\bW=\bX\t_{\bs g,\bZ,\bs h}\bY$ exists in the $2$-category $\Kur,$ with\/ $\vdim\bW=\vdim\bX+\vdim\bY-\vdim\bZ$. It is unique up to canonical equivalence in\/~$\Kur$.

This\/ $\bW$ is an orbifold if and only if\/ $\bs g,\bs h$ are transverse.

The topological space $W$ of\/ $\bW$ is given as a set by
\begin{align*}
W=\bigl\{(x,y,C)\colon x\in X,\;\> y\in Y,\;\> \bs g(x)=\bs h(y)=z\in Z&,\\
C\in G_x\bs g(G_x\bX)\backslash G_z\bZ/G_y\bs h(G_y\bY)\bigr\}&.
\end{align*}
For\/ $(x,y,C)\in W$ with\/ $\ga\in C\subseteq G_z\bZ,$ there is an exact sequence
\e
\begin{gathered}
\xymatrix@C=19pt@R=15pt{ 0 \ar[r] & T_{(x,y,C)}\bW \ar[rrr] &&& T_x\bX\op T_y\bY \ar[rrr]_(0.54){T_x\bs g\op (\ga\cdot T_y\bs h)} &&& *+[l]{T_z\bZ} \ar[d] \\
0 & O_z\bZ \ar[l] &&& O_x\bX\op O_y\bY \ar[lll]_(0.54){O_x\bs g\op (\ga\cdot O_y\bs h)} &&& *+[l]{O_{(x,y,C)}\bW.\!} \ar[lll] }
\end{gathered}
\label{kt5eq2}
\e

\label{kt5thm3}
\end{thm}

Here fibre products in 2-categories are explained in \S\ref{ktA3}. They are characterized by a universal property involving 2-morphisms. In Theorem \ref{kt5thm3}, it is essential that Kuranishi spaces $\Kur$ are a 2-category. The $\bW$ described in Theorem \ref{kt5thm3} generally does not satisfy a universal property in the homotopy category $\Ho(\Kur)$, and is not a fibre product in $\Ho(\Kur)$ in the category-theoretic sense. This is an important reason for making Kuranishi spaces into a 2-category rather than an ordinary category.

Note that w-transversality is equivalent to exactness of \eq{kt5eq2} at $O_z\bZ$, and exactness of the rest of \eq{kt5eq2} determines $T_{(x,y,C)}\bW,O_{(x,y,C)}\bW$ up to isomorphism. Transversality is equivalent to \eq{kt5eq2} being exact from $T_z\bZ$ onwards with $O_{(x,y,C)}\bW\ab =0$, so the claim that $\bW$ is an orbifold if and only if $\bs g,\bs h$ are transverse follows from Theorem~\ref{kt5thm1}(b).

If $\bs g\colon \bX\ra\bZ$ and $\bs h\colon \bY\ra\bZ$ are 1-morphisms in $\Kur$ with $\bZ$ a manifold or orbifold then $O_z\bZ=0$ for all $z\in\bZ$ by Theorem \ref{kt5thm1}, so $\bs g,\bs h$ are w-transverse by Definition \ref{kt5def4}, and Theorem \ref{kt5thm3} yields:

\begin{cor} Suppose $\bs g\colon \bX\ra\bZ,$ $\bs h\colon \bY\ra\bZ$ are $1$-morphisms of Kuranishi spaces, with\/ $\bZ$ a manifold or orbifold. Then the fibre product\/ $\bW=\bX\t_{\bs g,\bZ,\bs h}\bY$ exists in $\Kur,$ with\/~$\vdim\bW=\vdim\bX+\vdim\bY-\dim\bZ$.
\label{kt5cor1}
\end{cor}

In symplectic geometry it is often important to consider fibre products of Kuranishi spaces over manifolds or orbifolds, so Corollary \ref{kt5cor1} is very useful.

\begin{ex} If $\bX,\bY$ are Kuranishi spaces, and $\bZ=\bs *$ is the point in $\Kur$, then $\bX\t_{\bs *}\bY$ is the product $\bX\t\bY$, with $\vdim(\bX\t\bY)=\vdim\bX+\vdim\bY$. We can define a model for products of Kuranishi spaces explicitly: if $\bX=(X,\cI)$, $\bY=(Y,\cJ)$ with notation \eq{kt4eq2}--\eq{kt4eq3}, then we may write $\bX\t\bY=(X\t Y,\cK)$, where $\cK$ has indexing set $I\t J$ and Kuranishi neighbourhoods
\begin{align*}
&(W_{(i,j)},F_{(i,j)},\De_{(i,j)},t_{(i,j)},\om_{(i,j)})=\\
&\bigl(U_i\t V_j,\pi_{U_i}^*(D_i)\op \pi_{V_j}^*(E_j),\Be_i\t\Ga_j,\pi_{U_i}^*(r_i)\op \pi_{V_j}^*(s_j),\chi_i\t\psi_j\bigr).
\end{align*}
The remaining data in $\cK$ is also easy to write down explicitly. 

\label{kt5ex1}
\end{ex}

\subsection{Submersions and w-submersions}
\label{kt53}

A smooth map of manifolds $f\colon X\ra Y$ is called a {\it submersion\/} if $T_xf\colon T_xX\ra T_yY$ is surjective for all $x\in X$ with $f(x)=y\in Y$. If $g\colon X\ra Z$ and $h\colon Y\ra Z$ are smooth with $g$ a submersion then $g,h$ are transverse, so the fibre product $X\t_{g,Z,h}Y$ exists in $\Man$. Here is the analogue of this for Kuranishi spaces. We define two notions of submersion in $\Kur$, a weak and a strong.

\begin{dfn} A 1-morphism of Kuranishi spaces $\bs f\colon \bX\ra\bY$ is called a {\it weak submersion}, or {\it w-submersion}, if $O_x\bs f\colon O_x\bX\ra O_{\bs f(x)}\bY$ is surjective for all $x\in\bX$. It is called a {\it submersion\/} if $T_x\bs f\colon T_x\bX\ra T_{\bs f(x)}\bY$ is surjective and $O_x\bs f\colon O_x\bX\ra O_{\bs f(x)}\bY$ is an isomorphism for all $x\in\bX$. If $\bY$ is a manifold or orbifold then any 1-morphism $\bs f\colon \bX\ra\bY$ is a w-submersion.
\label{kt5def5}
\end{dfn}

If $\bs g\colon \bX\ra\bZ$ and $\bs h\colon \bY\ra\bZ$ are 1-morphisms in $\Kur$ with $\bs g$ a w-submersion then $\bs g,\bs h$ are w-transverse. If $\bs g$ is a submersion and $\bY$ is an orbifold then $\bs g,\bs h$ are transverse. Thus Theorem \ref{kt5thm3} implies:

\begin{cor} Suppose $\bs g\colon \bX\ra\bZ,$ $\bs h\colon \bY\ra\bZ$ are $1$-morphisms in $\Kur$ with\/ $\bs g$ a w-submersion. Then the fibre product\/ $\bW=\bX\t_{\bs g,\bZ,\bs h}\bY$ exists in $\Kur,$ with\/~$\vdim\bW=\vdim\bX+\vdim\bY-\vdim\bZ$.

If\/ $\bs g$ is a submersion and\/ $\bY$ is an orbifold, then\/ $\bW$ is an orbifold.
\label{kt5cor2}
\end{cor}

\appendix

\section{Background from Category Theory and Algebraic Geometry}
\label{ktA}

Sections \ref{ktA1}--\ref{ktA3} discuss 2-categories, and \S\ref{ktA4} defines sheaves and stacks on topological spaces.

\subsection{Basics of 2-categories}
\label{ktA1}

We discuss 2-categories, both strict and weak. Some references are Borceux \cite[\S 7]{Borc}, Kelly and Street \cite{KeSt}, and Behrend et al.~\cite[App.~B]{BEFF}.

\begin{dfn} A {\it strict\/} 2-{\it category\/} $\bs\cC$ consists of
a class of {\it objects\/} $\Obj(\bs\cC)$, for all $X,Y\in\Obj(\bs\cC)$ an essentially small category $\Hom(X,Y)$, for all $X,Y,Z$ in $\Obj(\bs\cC)$ a functor $\mu_{X,Y,Z}\colon \Hom(X,Y)\t\Hom(Y,Z)\ra\Hom(X,Z)$ called {\it composition}, and for all $X$ in $\Obj(\bs\cC)$ an object $\id_X$ in $\Hom(X,X)$ called the {\it identity $1$-morphism}. These must satisfy
the {\it associativity property}, that
\e
\mu_{W,Y,Z}\ci(\mu_{W,X,Y}\t\id_{\Hom(Y,Z)})
=\mu_{W,X,Z}\ab\ci\ab(\id_{\Hom(W,X)}\ab\t\mu_{X,Y,Z})
\label{ktAeq1}
\e
as functors $\Hom(W,X)\t\Hom(X,Y)\t\Hom(Y,Z)\ra\Hom(W,X)$, and the
{\it identity property}, that
\e
\mu_{X,X,Y}(\id_X,-)=\mu_{X,Y,Y}(-,\id_Y)=\id_{\Hom(X,Y)}
\label{ktAeq2}
\e
as functors $\Hom(X,Y)\ra\Hom(X,Y)$.

Objects $f$ of $\Hom(X,Y)$ are called 1-{\it morphisms}, written
$f\colon X\ra Y$. For 1-morphisms $f,g\colon X\ra Y$, morphisms $\eta\in
\Hom_{\Hom(X,Y)}(f,g)$ are called 2-{\it mor\-ph\-isms}, written
$\eta\colon f\Ra g$. Thus, a 2-category has objects $X$, and two kinds of
morphisms: 1-morphisms $f\colon X\ra Y$ between objects, and 2-morphisms
$\eta\colon f\Ra g$ between 1-morphisms.

A {\it weak\/ $2$-category}, or {\it bicategory}, is like a strict 2-category, except that the equations of functors \eq{ktAeq1}, \eq{ktAeq2} are required to hold only up to specified natural isomorphisms. That is, a weak 2-category $\bs\cC$ consists of data $\Obj(\bs\cC),\Hom(X,Y),\ab\mu_{X,Y,Z},\ab\id_X$ as above, but in place of \eq{ktAeq1}, a natural isomorphism of functors 
\e
\al\colon \mu_{W,Y,Z}\ci(\mu_{W,X,Y}\t\id_{\Hom(Y,Z)})
\Longra\mu_{W,X,Z}\ci(\id_{\Hom(W,X)}\t\mu_{X,Y,Z}),
\label{ktAeq3}
\e
and in place of \eq{ktAeq2}, natural isomorphisms
\e
\be\colon \mu_{X,X,Y}(\id_X,-)\Longra \id_{\Hom(X,Y)},\;\>
\ga\colon \mu_{X,Y,Y}(-,\id_Y)\Longra \id_{\Hom(X,Y)}.
\label{ktAeq4}
\e
These $\al,\be,\ga$ must satisfy identities which we give below in \eq{ktAeq6} and~\eq{ktAeq8}.

A strict 2-category $\bs\cC$ can be regarded as an example of a weak 2-category, in which the natural isomorphisms $\al,\be,\ga$ in \eq{ktAeq3}--\eq{ktAeq4} are the identities.
\label{ktAdef1}
\end{dfn}

We now unpack Definition \ref{ktAdef1}, making it more explicit.

There are three kinds of composition in a 2-category, satisfying
various associativity relations. If $f\colon X\ra Y$ and $g\colon Y\ra Z$ are
1-morphisms then $\mu_{X,Y,Z}(f,g)$ is the {\it composition of\/ $1$-morphisms}, written $g\ci f\colon X\ra Z$. If $f,g,h\colon X\ra Y$ are 1-morphisms and $\eta\colon f\Ra g$, $\ze\colon g\Ra h$ are 2-morphisms then composition of $\eta,\ze$ in $\Hom(X,Y)$ gives the {\it vertical composition of\/ $2$-morphisms} of $\eta,\ze$, written $\ze\od\eta\colon f\Ra h$, as a diagram\vskip -15pt\begin{equation*}
\xymatrix@C=25pt{ X \rruppertwocell^f{\eta} \rrlowertwocell_h{\ze}
\ar[rr]_(0.35)g && Y & \ar@{~>}[r] && X
\rrtwocell^f_h{{}\,\,\,\,\ze\od\eta\!\!\!\!\!} && Y.}
\end{equation*}\vskip -10pt
\noindent Vertical composition is associative.

If $f,\dot f\colon X\ra Y$ and $g,\dot g\colon Y\ra Z$ are 1-morphisms and
$\eta\colon f\Ra\dot f$, $\ze\colon g\Ra\dot g$ are 2-morphisms then
$\mu_{X,Y,Z}(\eta,\ze)$ is the {\it horizontal composition of\/
$2$-morphisms}, written $\ze*\eta\colon g\ci f\Ra\dot g\ci\dot f$, as a
diagram\vskip -15pt\begin{equation*}
\xymatrix@C=20pt{ X \rrtwocell^f_{\dot f}{\eta} && Y
\rrtwocell^g_{\dot g}{\ze} && Z & \ar@{~>}[r] && X \rrtwocell^{g\ci
f}_{\dot g\ci\dot f}{{}\,\,\,\ze*\eta\!\!\!\!\!} && Z. }
\end{equation*}\vskip -5pt
\noindent As $\mu_{X,Y,Z}$ is a functor, these satisfy {\it compatibility of vertical and horizontal composition\/}: given a diagram of 1- and 2-morphisms\vskip -15pt\begin{equation*}
\xymatrix@C=20pt{ X \rruppertwocell^f{\eta} \rrlowertwocell_{\ddot f}{\dot\eta} \ar[rr]^(0.3){\dot f} && Y \rruppertwocell^g{\ze} \rrlowertwocell_{\ddot g}{\dot\ze} \ar[rr]^(0.3){\dot g} && Z, }
\end{equation*}\vskip -10pt
\noindent we have
\begin{equation*}
(\dot\ze\od\ze)*(\dot\eta\od\eta)=(\dot\ze*\dot\eta)\od(\ze*\eta)\colon g\ci f\Longra \ddot g\ci\ddot f.
\end{equation*}
There are also two kinds of identity: {\it identity\/
$1$-morphisms\/} $\id_X\colon X\ra X$ and {\it identity\/
$2$-morphisms\/}~$\id_f\colon f\Ra f$.

In a strict 2-category $\bs\cC$, composition of 1-morphisms is strictly associative, $(g\ci f)\ci e=g\ci(f\ci e)$, and horizontal composition of 2-morphisms is strictly associative, $(\ze*\eta)*\ep=\ze*(\eta*\ep)$. In a weak 2-category $\bs\cC$, composition of 1-morphisms is associative up to specified 2-isomorphisms. That is, if $e\colon W\ra X$, $f\colon X\ra Y$, $g\colon Y\ra Z$ are 1-morphisms in $\bs\cC$ then the natural isomorphism $\al$ in \eq{ktAeq3} gives a 2-isomorphism
\e
\al_{g,f,e}\colon (g\ci f)\ci e\Longra g\ci(f\ci e).
\label{ktAeq5}
\e
As $\al$ is a natural isomorphism, given 1-morphisms $e,\dot e\colon W\ra X$, $f,\dot f\colon X\ra Y$, $g,\dot g\colon Y\ra Z$ and 2-morphisms $\ep\colon e\Ra\dot e$, $\eta\colon f\Ra\dot f$, $\ze\colon g\Ra\dot g$ in $\bs\cC$, the following diagram of 2-morphisms must commute:
\begin{equation*}
\xymatrix@C=150pt@R=14pt{ *+[r]{(g\ci f)\ci e} \ar@{=>}[r]_{\al_{g,f,e}} \ar@{=>}[d]^{(\ze*\eta)*\ep} & *+[l]{g\ci (f\ci e)} \ar@{=>}[d]_{\ze*(\eta*\ep)} \\
*+[r]{(\dot g\ci\dot f)\ci\dot e} \ar@{=>}[r]^{\al_{\dot g,\dot f,\dot e}} & *+[l]{\dot g\ci(\dot f\ci\dot e).\!\!} }
\end{equation*}

The $\al_{g,f,e}$ must satisfy the {\it associativity coherence axiom\/}: if $d\colon V\ra W$ is another 1-morphism, then the following diagram of 2-morphisms must commute:
\e
\begin{gathered}
\xymatrix@C=76pt@R=14pt{ *+[r]{((g\ci f)\ci e)\ci d} \ar@{=>}[r]_(0.6){\al_{g,f,e}*\id_d} \ar@{=>}[d]^{\al_{g\ci f,e,d}} & (g\ci (f\ci e)) \ci d \ar@{=>}[r]_(0.4){\al_{g,f\ci e,d}} & *+[l]{g\ci ((f\ci e) \ci d)} \ar@{=>}[d]_{\id_g*\al_{f,e,d}} \\
*+[r]{(g\ci f)\ci (e \ci d)} \ar@{=>}[rr]^{\al_{g,f,d\ci e}} && *+[l]{g\ci (f\ci (e \ci d)).\!\!}}
\end{gathered}
\label{ktAeq6}
\e

In a strict 2-category $\bs\cC$, given a 1-morphism $f\colon X\ra Y$, the identity 1-morphisms $\id_X,\id_Y$ satisfy $f\ci\id_X=\id_Y\ci f=f$. In a weak 2-category $\bs\cC$, the natural isomorphisms $\be,\ga$ in \eq{ktAeq4} give 2-isomorphisms
\e
\be_f\colon f\ci\id_X\Longra f,\qquad \ga_f\colon \id_Y\ci f\Longra f.
\label{ktAeq7}
\e
As $\be,\ga$ are natural isomorphisms, if $\eta\colon f\Ra\dot f$ is a 2-morphism we must have
\begin{align*}
\eta\od\be_f&=\be_{\dot f}\od(\eta*\id_{\id_X})\colon f\ci\id_X\Ra\dot f,\\
\eta\od\ga_f&=\ga_{\dot f}\od(\id_{\id_Y}*\eta)\colon \id_Y\ci f\Ra\dot f.
\end{align*}

The $\be_f,\ga_f$ must satisfy the {\it identity coherence axiom\/}: if $g\colon Y\ra Z$ is another 1-morphism, then the following diagram of 2-morphisms must commute:
\e
\begin{gathered}
\xymatrix@C=120pt@R=-1pt{ *+[r]{(g\ci\id_Y)\ci f} \ar@{=>}[dr]^{\be_g*\id_f} \ar@{=>}[dd]^{\al_{g,\id_Y,f}} \\
& g\ci f. \\
*+[r]{g\ci(\id_Y\ci f)} \ar@{=>}[ur]_{\id_g*\al_f} }
\end{gathered}
\label{ktAeq8}
\e

A basic example of a strict 2-category is the {\it $2$-category of categories\/} $\mathfrak{Cat}$, with objects small categories $\cC$, 1-morphisms
functors $F\colon \cC\ra\cD$, and 2-morphisms natural transformations
$\eta\colon F\Ra G$ for functors $F,G\colon \cC\ra\cD$. Orbifolds naturally form
a 2-category (strict or weak, depending on the definition), as in \S\ref{kt43}, and so do stacks in algebraic geometry.

In a 2-category $\bs\cC$, there are three notions of when objects $X,Y$
in $\bs\cC$ are `the same': {\it equality\/} $X=Y$, and 1-{\it
isomorphism}, that is we have 1-morphisms $f\colon X\ra Y$, $g\colon Y\ra X$
with $g\ci f=\id_X$ and $f\ci g=\id_Y$, and {\it equivalence}, that
is we have 1-morphisms $f\colon X\ra Y$, $g\colon Y\ra X$ and 2-isomorphisms
$\eta\colon g\ci f\Ra\id_X$ and $\ze\colon f\ci g\Ra\id_Y$. Usually equivalence
is the correct notion. 

When we say that objects $X,Y$ in a 2-category $\bs\cC$ are {\it canonically equivalent}, we mean that there is a nonempty distinguished class of equivalences $f\colon X\ra Y$ in $\bs\cC$, and given any two such equivalences $f,g\colon X\ra Y$ there is a 2-isomorphism $\eta\colon f\Ra g$. Often there is a distinguished choice of such~$\eta$.

{\it Commutative diagrams\/} in 2-categories should in general only
commute {\it up to (specified)\/ $2$-isomorphisms}, rather than
strictly. A simple example of a commutative diagram in a 2-category
$\bs\cC$ is
\begin{equation*}
\xymatrix@C=50pt@R=8pt{ & Y \ar[dr]^g \ar@{=>}[d]^\eta \\
X \ar[ur]^f \ar[rr]_h && Z, }
\end{equation*}
which means that $X,Y,Z$ are objects of $\bs\cC$, $f\colon X\ra Y$, $g\colon Y\ra
Z$ and $h\colon X\ra Z$ are 1-morphisms in $\bs\cC$, and $\eta\colon g\ci f\Ra h$
is a 2-isomorphism.

Let $\bs\cC$ be a 2-category. The {\it homotopy category\/} $\Ho(\bs\cC)$ of $\bs\cC$ is the category whose objects are objects of $\bs\cC$, and whose morphisms $[f]\colon X\ra Y$ are 2-isomorphism classes $[f]$ of 1-morphisms $f\colon X\ra Y$ in $\bs\cC$. Then equivalences in $\bs\cC$ become isomorphisms in $\Ho(\bs\cC)$, 2-commutative diagrams in $\bs\cC$ become commutative diagrams in $\Ho(\bs\cC)$, and so on.

\subsection{2-functors between 2-categories}
\label{ktA2}

Next we discuss 2-functors between 2-categories, following Borceux \cite[\S 7.2, \S 7.5]{Borc} and Behrend et al.~\cite[\S B.4]{BEFF}.

\begin{dfn} Let $\bs\cC,\bs\cD$ be strict 2-categories. A {\it strict\/ $2$-functor\/} $F\colon\bs\cC\ra\bs\cD$ assigns an object $F(X)$ in $\bs\cD$ for each object $X$ in $\bs\cC$, a 1-morphism $F(f):F(X)\ra F(Y)$ in $\bs\cD$ for each 1-morphism $f\colon X\ra Y$ in $\bs\cC$, and a 2-morphism $F(\eta)\colon F(f)\Ra F(g)$ in $\bs\cD$ for each 2-morphism $\eta\colon f\Ra g$ in $\bs\cC$, such that $F$ preserves all the structures on $\bs\cC,\bs\cD$, that is,
\ea
F(g\ci f)&=F(g)\ci F(f), & F(\id_X)&=\id_{F(X)}, & F(\ze*\eta)&=F(\ze)* F(\eta), 
\label{ktAeq9}\\
F(\ze\od\eta)&=F(\ze)\od F(\eta), & F(\id_f)&=\id_{F(f)}.
\label{ktAeq10}
\ea

Now let $\bs\cC,\bs\cD$ be weak 2-categories. Then strict 2-functors $F\colon \bs\cC\ra\bs\cD$ are not well-behaved. To fix this, we need to relax \eq{ktAeq9} to hold only up to specified 2-isomorphisms. A {\it weak\/ $2$-functor\/} (or {\it pseudofunctor\/}) $F\colon \bs\cC\ra\bs\cD$ assigns an object $F(X)$ in $\bs\cD$ for each object $X$ in $\bs\cC$, a 1-morphism $F(f)\colon F(X)\ra F(Y)$ in $\bs\cD$ for each 1-morphism $f\colon X\ra Y$ in $\bs\cC$, a 2-morphism $F(\eta)\colon F(f)\Ra F(g)$ in $\bs\cD$ for each 2-morphism $\eta\colon f\Ra g$ in $\bs\cC$, a 2-isomorphism $F_{g,f}\colon F(g)\ci F(f)\Ra F(g\ci f)$ in $\bs\cD$ for all 1-morphisms $f\colon X\ra Y$, $g\colon Y\ra Z$ in $\bs\cC$, and a 2-isomorphism $F_X\colon F(\id_X)\Ra \id_{F(X)}$ in $\bs\cD$ for all objects $X$ in $\bs\cC$ such that \eq{ktAeq10} holds, and for all $e\colon W\ra X$, $f\colon X\ra Y$, $g\colon Y\ra Z$ in $\bs\cC$ the following diagram of 2-isomorphisms commutes in $\bs\cD$:
\begin{equation*}
\xymatrix@C=87pt@R=15pt{ *+[r]{(F(g)\ci F(f))\ci F(e)} \ar@{=>}[d]^{\al_{F(g),F(f),F(e)}} \ar@{=>}[r]_(0.65){F_{g,f}*\id_{F(e)}} & F(g\ci f)\ci F(e) \ar@{=>}[r]_(0.4){F_{g\ci f,e}} & *+[l]{F((g\ci f)\ci e)} \ar@{=>}[d]_{F(\al_{g,f,e})} \\
*+[r]{F(g)\ci (F(f)\ci F(e))} \ar@{=>}[r]^(0.65){\id_{F(g)}*F_{f,e}} & F(g)\ci F(f\ci e) \ar@{=>}[r]^(0.4){F_{g,f\ci e}} &
*+[l]{F(g\ci (f\ci e)),\!\!} } 
\end{equation*}
and for all 1-morphisms $f\colon X\ra Y$ in $\bs\cC$, the following commute in $\bs\cD$:
\begin{equation*}
\xymatrix@!0@C=25pt@R=30pt{
*+[r]{F(f)\ci F(\id_X)} \ar@{=>}[rrrrrr]_(0.58){F_{f,\id_X}} \ar@{=>}[d]^{\id_{F(f)}*F_X} &&&&&& *+[l]{F(f\ci\id_X)} \ar@{=>}[d]_{F(\be_f)} & *+[r]{F(\id_Y)\ci F(f)} \ar@{=>}[rrrrrr]_(0.58){F_{\id_Y,f} } \ar@{=>}[d]^{F_Y*\id_{F(f)}} &&&&&& *+[l]{F(\id_Y\ci f)} \ar@{=>}[d]_{F(\ga_f)} \\
*+[r]{F(f)\ci \id_{F(X)}} \ar@{=>}[rrrrrr]^(0.6){\be_{F(f)}} &&&&&& *+[l]{F(f),\!\!} &
*+[r]{\id_{F(Y)}\ci F(f)} \ar@{=>}[rrrrrr]^(0.6){\ga_{F(f)}} &&&&&& *+[l]{F(f),\!\!} } 
\end{equation*}
and if $f,\dot f\colon X\ra Y$ and $g,\dot g\colon Y\ra Z$ are 1-morphisms and
$\eta\colon f\Ra\dot f$, $\ze\colon g\Ra\dot g$ are 2-morphisms in $\bs\cC$ then the following commutes in $\bs\cD$:
\begin{equation*}
\xymatrix@C=140pt@R=15pt{ *+[r]{F(g)\ci F(f)} \ar@{=>}[d]^{F(\ze)*F(\eta)} \ar@{=>}[r]_{F_{g,f}} & *+[l]{F(g\ci f)} \ar@{=>}[d]_{F(\ze*\eta)} \\
*+[r]{F(\dot g)\ci F(\dot f)} \ar@{=>}[r]^{F_{\dot g,\dot f}} & *+[l]{F(\dot g\ci\dot f).\!\!} } 
\end{equation*}

There are obvious notions of {\it composition\/} $G\ci F$ of strict and weak 2-functors $F\colon \bs\cC\ra\bs\cD$, $G\colon \bs\cD\ra\bs\cE$, {\it identity\/ $2$-functors\/} $\id_{\bs\cC}$, and so on.

If $\bs\cC,\bs\cD$ are strict 2-categories, then a strict $2$-functor $F\colon \bs\cC\ra\bs\cD$ can be made into a weak 2-functor by taking all $F_{g,f},F_X$ to be identity 2-morphisms. 

We define a weak 2-functor $F\colon \bs\cC\ra\bs\cD$ to be an {\it equivalence of weak\/ $2$-categories}, if for all objects $X,Y$ in $\bs\cC$, the functor $F_{X,Y}\colon \Hom_{\bs\cC}(X,Y)\ra\Hom_{\bs\cD}(F(X),F(Y))$ is an equivalence of categories, and the map induced by $F$ from equivalence classes of objects in $\bs\cC$ to equivalence classes of objects in $\bs\cD$ is surjective (and hence a bijection). There is an alternative definition using 2-natural transformations of 2-functors, which we will not explain.
\label{ktAdef2}
\end{dfn}

\subsection{Fibre products in 2-categories}
\label{ktA3}

We define fibre products in 2-categories, as in Behrend et al.~\cite[Def.~B.13]{BEFF}.

\begin{dfn} Let $\bs\cC$ be a strict 2-category and $g\colon X\ra Z$, $h\colon Y\ra Z$ be 1-morphisms in $\bs\cC$. A {\it fibre product\/} $X\t_ZY$ in $\bs\cC$
consists of an object $W$, 1-morphisms $\pi_X\colon W\ra X$ and $\pi_Y\colon W\ra Y$ and a 2-isomorphism $\eta\colon g\ci\pi_X\Ra h\ci\pi_Y$ in $\bs\cC$ with the
following universal property: suppose $\pi_X'\colon W'\ra X$ and
$\pi_Y'\colon W'\ra Y$ are 1-morphisms and $\eta'\colon g\ci\pi_X'\Ra
h\ci\pi_Y'$ is a 2-isomorphism in $\bs\cC$. Then there should exist a
1-morphism $b\colon W'\ra W$ and 2-isomorphisms $\ze_X\colon \pi_X\ci
b\Ra\pi_X'$, $\ze_Y\colon \pi_Y\ci b\Ra\pi_Y'$ such that the following
diagram of 2-isomorphisms commutes:
\e
\begin{gathered}
\xymatrix@C=100pt@R=14pt{ *+[r]{g\ci\pi_X\ci b} \ar@{=>}[r]_{\eta*\id_b}
\ar@{=>}[d]^{\,\id_g*\ze_X} & *+[l]{h\ci\pi_Y\ci b}
\ar@{=>}[d]_{\id_h*\ze_Y} \\ *+[r]{g\ci\pi_X'}
\ar@{=>}[r]^{\eta'} & *+[l]{h\ci\pi_Y'.\!\!} }
\end{gathered}
\label{ktAeq11}
\e
Furthermore, if $\ti b,\ti\ze_X,\ti\ze_Y$ are alternative choices of
$b,\ze_X,\ze_Y$ then there should exist a unique 2-isomorphism
$\th\colon \ti b\Ra b$ with
\begin{equation*}
\ti\ze_X=\ze_X\od(\id_{\pi_X}*\th)\quad\text{and}\quad
\ti\ze_Y=\ze_Y\od(\id_{\pi_Y}*\th).
\end{equation*}
If a fibre product $X\t_ZY$ in $\bs\cC$ exists then it is unique up to
canonical equivalence. That is, if $W,W'$ are two choices for $X\t_ZY$ then there is an equivalence $b\colon W\ra W'$ in $\bs\cC$, which is canonical up to 2-isomorphism.

If instead $\bs\cC$ is a weak 2-category, we must replace \eq{ktAeq11} by
\begin{equation*}
\xymatrix@C=80pt@R=14pt{ 
*+[r]{(g\ci\pi_X)\ci b} \ar@{=>}[d]^{\al_{g,\pi_X,b}}
\ar@{=>}[r]_(0.6){\eta*\id_b} & {(h\ci\pi_Y)\ci b}
\ar@{=>}[r]_(0.4){\al_{h,\pi_Y,b}} & *+[l]{h\ci(\pi_Y\ci b)\!\!}
\ar@{=>}[d]_{\id_h*\ze_Y} \\ 
*+[r]{g\ci(\pi_X\ci b)} \ar@{=>}[r]^(0.6){\id_g*\ze_X} &
g\ci\pi_X' \ar@{=>}[r]^(0.4){\eta'} & *+[l]{h\ci\pi_Y'.\!\!} }
\end{equation*}

\label{ktAdef3}
\end{dfn}

Orbifolds, and stacks in algebraic geometry, form 2-categories, and
Definition \ref{ktAdef3} is the right way to define fibre products
of orbifolds or stacks.

\subsection{Sheaves and stacks on topological spaces}
\label{ktA4}

The next definition of sheaves on a topological space, as in Hartshorne \cite[\S II.1]{Hart}, will not actually be used in this paper, but we include it as motivation for the following definition of stacks on a topological space.

\begin{dfn} Let $X$ be a topological space. A {\it presheaf of sets\/} $\cE$ on $X$ consists of the data of a set $\cE(S)$ for every open set $S\subseteq X$, and a map $\rho_{ST}\colon \cE(S)\ra\cE(T)$ called the {\it restriction map\/} for every inclusion $T\subseteq S\subseteq X$ of open sets, satisfying the conditions that
\begin{itemize}
\setlength{\itemsep}{0pt}
\setlength{\parsep}{0pt}
\item[(i)] $\rho_{SS}=\id_{\cE(S)}\colon \cE(S)\ra\cE(S)$ for all open
$S\subseteq X$; and
\item[(ii)] $\rho_{SU}=\rho_{TU}\ci\rho_{ST}\colon \cE(S)\ra\cE(U)$ for all
open~$U\subseteq T\subseteq S\subseteq X$.
\end{itemize}

A presheaf of sets $\cE$ on $X$ is called a {\it sheaf\/} if it also satisfies
\begin{itemize}
\setlength{\itemsep}{0pt}
\setlength{\parsep}{0pt}
\item[(iii)] If $S\subseteq X$ is open, $\{T_i\colon i\in I\}$ is an open
cover of $S$, and $s,t\in\cE(S)$ have $\rho_{ST_i}(s)=\rho_{ST_i}(t)$ in
$\cE(T_i)$ for all $i\in I$, then $s=t$ in $\cE(S)$; and
\item[(iv)] If $S\subseteq X$ is open, $\{T_i\colon i\in I\}$ is an open cover of
$S$, and we are given elements $s_i\in\cE(T_i)$ for all $i\in I$
such that $\rho_{T_i(T_i\cap T_j)}(s_i)=\rho_{T_j(T_i\cap
T_j)}(s_j)$ in $\cE(T_i\cap T_j)$ for all $i,j\in I$, then there
exists $s\in\cE(S)$ with $\rho_{ST_i}(s)=s_i$ for all $i\in I$.
This $s$ is unique by~(iii).
\end{itemize}

Suppose $\cE,\cF$ are presheaves or sheaves of sets on $X$. A {\it morphism\/} $\phi\colon \cE\ra\cF$ consists of a map $\phi(S)\colon \cE(S)\ra\cF(S)$ for all open $S\subseteq
X$, such that the following diagram commutes for all open $T\subseteq S\subseteq X$
\begin{equation*}
\xymatrix@C=95pt@R=15pt{
*+[r]{\cE(S)} \ar[r]_{\phi(S)} \ar[d]^{\rho_{ST}} & *+[l]{\cF(S)}
\ar[d]_{\rho_{ST}'} \\ *+[r]{\cE(T)} \ar[r]^{\phi(T)} & *+[l]{\cF(T),\!} }
\end{equation*}
where $\rho_{ST}$ is the restriction map for $\cE$, and $\rho_{ST}'$
the restriction map for~$\cF$.
\label{ktAdef4}
\end{dfn}

Sheaves are basic objects in algebraic geometry. Informally, something forms a sheaf if it can be defined locally and glued. For example, if $X,Y$ are manifolds then smooth maps $f\colon X\ra Y$ form a sheaf on $X$, because an arbitrary map $f\colon X\ra Y$ is smooth if and only if it is a smooth on an open neighbourhood of each point in $X$. One uses the sheaf property of smooth maps all the time in differential geometry without noticing. 

Finally we define stacks on topological spaces. Stacks are the 2-category version of sheaves, roughly, they are sheaves of groupoids rather than sheaves of sets. Informally, something forms a stack if it can be defined locally and glued up to 2-isomorphisms. For example, if $X,Y$ are orbifolds then smooth maps $f\colon X\ra Y$ (i.e. 1-morphisms of orbifolds) and their 2-morphisms form a stack on $X$, but in general they do not form a sheaf on $X$.

\begin{dfn} Let $X$ be a topological space. A {\it prestack\/} (or {\it prestack in groupoids}, or 2-{\it presheaf\/}) $\bcE$ on $X$, consists of the data of a groupoid $\bcE(S)$ for every open set $S\subseteq X$, and a functor $\rho_{ST}\colon \bcE(S)\ra\bcE(T)$ called the {\it restriction map\/} for every inclusion $T\subseteq S\subseteq X$ of open sets, and a natural isomorphism of functors $\eta_{STU}\colon \rho_{TU}\ci\rho_{ST}\Ra \rho_{SU}$ for all inclusions $U\subseteq T\subseteq S\subseteq X$ of open sets, satisfying the conditions that:
\begin{itemize}
\setlength{\itemsep}{0pt}
\setlength{\parsep}{0pt}
\item[(i)] $\rho_{SS}=\id_{\bcE(S)}\colon \bcE(S)\ra\bcE(S)$ for all open
$S\subseteq X$, and $\eta_{SST}=\eta_{STT}=\id_{\rho_{ST}}$ for all open $T\subseteq S\subseteq X$; and
\item[(ii)] $\eta_{SUV}\od(\id_{\rho_{UV}}*\eta_{STU})=\eta_{STV}\od(\eta_{TUV}*\id_{\rho_{ST}})\colon \rho_{UV}\ci\rho_{TU}\ci\rho_{ST}\Longra\rho_{SV}$ for all open~$V\subseteq U\subseteq T\subseteq S\subseteq X$.
\end{itemize}

A prestack $\bcE$ on $X$ is called a {\it stack\/} (or {\it stack in groupoids}, or 2-{\it sheaf\/}) on $X$ if whenever $S\subseteq X$ is open and $\{T_i\colon i\in I\}$ is an open cover of $S$, then:
\begin{itemize}
\setlength{\itemsep}{0pt}
\setlength{\parsep}{0pt}
\item[(iii)] If $\al,\be\colon A\ra B$ are morphisms in $\bcE(S)$ and $\rho_{ST_i}(\al)=\rho_{ST_i}(\be)\colon \rho_{ST_i}(A)\ra \rho_{ST_i}(B)$ in $\bcE(T_i)$ for all $i\in I$, then $\al=\be$.
\item[(iv)] If $A,B$ are objects of $\bcE(S)$ and $\al_i\colon \rho_{ST_i}(A)\ra \rho_{ST_i}(B)$ are morphisms in $\bcE(T_i)$ for all $i\in I$ with
\begin{align*}
&\eta_{ST_i(T_i\cap T_j)}(B)\ci\rho_{T_i(T_i\cap T_j)}(\al_i)\ci\eta_{ST_i(T_i\cap T_j)}(A)^{-1}\\
&\quad=\eta_{ST_j(T_i\cap T_j)}(B)\ci\rho_{T_j(T_i\cap T_j)}(\al_j)\ci\eta_{ST_j(T_i\cap T_j)}(A)^{-1}
\end{align*}
in $\bcE(T_i\cap T_j)$ for all $i,j\in I$, then there exists $\al\colon A\ra B$ in $\bcE(S)$ (necessarily unique by (iii)) with $\rho_{ST_i}(\al)=\al_i$ for all~$i\in I$.
\item[(v)] If $A_i\in\bcE(T_i)$ for $i\in I$ and $\al_{ij}\colon \rho_{T_i(T_i\cap T_j)}(A_i)\ra \rho_{T_j(T_i\cap T_j)}(A_j)$ are morphisms in $\bcE(T_i\cap T_j)$ for all $i,j\in I$ satisfying
\end{itemize}
\begin{align*}
& \\[-22pt]
&\eta_{T_k(T_j\cap T_k)(T_i\cap T_j\cap T_k)}(A_k)\ci\rho_{(T_j\cap T_k)(T_i\cap T_j\cap T_k)}(\al_{jk})\ci\eta_{T_j(T_j\cap T_k)(T_i\cap T_j\cap T_k)}(A_j)^{-1}\\
&\ci\eta_{T_j(T_i\cap T_j)(T_i\cap T_j\cap T_k)}(A_j)\ci \rho_{(T_i\cap T_j)(T_i\cap T_j\cap T_k)}(\al_{ij})\ci\eta_{T_i(T_i\cap T_j)(T_i\cap T_j\cap T_k)}(A_i)^{-1}\\
&=\eta_{T_k(T_i\cap T_k)(T_i\cap T_j\cap T_k)}(A_k)\ci\rho_{(T_i\cap T_k)(T_i\cap T_j\cap T_k)}(\al_{ik})\ci\eta_{T_i(T_i\cap T_k)(T_i\cap T_j\cap T_k)}(A_i)^{-1}
\end{align*}
\begin{itemize}
\setlength{\itemsep}{0pt}
\setlength{\parsep}{0pt}
\item[]for all $i,j,k\in I$, then there exist an object $A$ in $\bcE(S)$ and morphisms $\be_i\colon A_i\ra\rho_{ST_i}(A)$ for $i\in I$ such that for all $i,j\in I$ we have
\begin{equation*}
\eta_{ST_i(T_i\cap T_j)}(A)\ci\rho_{T_i(T_i\cap T_j)}(\be_i)=\eta_{ST_j(T_i\cap T_j)}(A)\ci\rho_{T_j(T_i\cap T_j)}(\be_j)\ci\al_{ij}.
\end{equation*}

If $\ti A,\ti\be_i$ are alternative choices then (iii),(iv) imply there is a unique isomorphism $\ga\colon A\ra\ti A$ in $\bcE(S)$ with $\rho_{ST_i}(\ga)=\ti\be_i\ci\be_i^{-1}$ for all~$i\in I$.
\end{itemize}
\label{ktAdef5}
\end{dfn}

In the examples of stacks on topological spaces that will be important to us, we will have $\rho_{TU}\ci\rho_{ST}=\rho_{SU}$ and $\eta_{STU}=\id_{\rho_{SU}}$ for all open $U\subseteq T\subseteq S\subseteq X$. So (ii) is automatic, and all the $\eta_{\cdots}(\cdots)$ terms in (iv),(v) can be omitted.

\medskip

\noindent{\sc The Mathematical Institute, Radcliffe
Observatory Quarter, Woodstock Road, Oxford, OX2 6GG, U.K.

\noindent E-mail: {\tt joyce@maths.ox.ac.uk.}}

\end{document}